\documentclass{article}

\usepackage{arxiv}

\usepackage[T1]{fontenc}    
\usepackage{hyperref}       
\usepackage{url}            
\usepackage{booktabs}       
\usepackage{amsfonts}       
\usepackage{nicefrac}       
\usepackage{microtype}      
\usepackage{amsmath}
\usepackage{amssymb}
\usepackage{amsmath, amsthm}
\usepackage{graphicx}
\usepackage[colorinlistoftodos]{todonotes}
\usepackage{datetime}
\usepackage{tikz}
\usetikzlibrary{calc,trees,positioning,arrows,chains,shapes.geometric,%
decorations.pathreplacing,decorations.pathmorphing,shapes,%
matrix,shapes.symbols,plotmarks,decorations.markings,shadows}

\definecolor{mygray}{gray}{0.8}
\definecolor{mmygray}{gray}{0.6}    
    
\usepackage{tikz}
\usepackage{verbatim}
\usepackage{subfig}

\theoremstyle{definition}
\newtheorem{definition}{Definition}[section]
\newtheorem{theorem}{Theorem}[section]
\newtheorem{lemma}[theorem]{Lemma}
\newtheorem{proposition}[theorem]{Proposition}
\newtheorem{corollary}[theorem]{Corollary}
\newtheorem{example}[theorem]{Example}
\newtheorem{remark}[theorem]{Remark}

\newcommand\norm[1]{\left\lVert#1\right\rVert}
\newcommand{\turnnw}[1]{\rotatebox[origin=c]{270}{\ensuremath#1}}

\DeclareUnicodeCharacter{391}{$\Alpha$}
\DeclareUnicodeCharacter{395}{$\Epsilon$}

\title{On hybrid  order dimensions}

\author{
  Athanasios Andrikopoulos\thanks{Associate professor  (https://www.ceid.upatras.gr/webpages/faculty/aandriko/)} \\
  Dept. of Computer Engineering and Informatics\\
  University of Patras\\
  Patras, 26504, Greece \\
  \texttt{aandriko@ceid.upatras.gr} \\
}

\begin{document}
\maketitle

\begin{abstract}
The notion of interval order was introduced by Norbert Wiener \cite{wie}
in order to clarify the relation between the notion of an instant of time and that of a period of time. This was a problem
on which 
Bertrand
Russell \cite{rus} worked at the time.
Interval orders play
an important role in many areas of pure and applied mathematics, graph theory, computer science and engineering.
Special cases of interval order are the semiorder and linear order.
All of these notions 
are especially important in the study of linear-interval and linear-semiorder dimension of a binary relation. 
This kind of dimension, which we call {\it hybrid order dimension},
gives a common generalization of linear order and
interval order (semiorder) dimension and
is
arguably the most important measure of ordered set complexity. 
In this paper,
we present three main results of the theory of hybrid order dimension. More specifically, we obtain necessary and sufficient conditions for a binary relation to have an interval order (resp. linear-interval order, linear-simiorder) extension, 
as well as an interval order realizer of interval orders (resp. linear-interval orders, linear-simiorders).
We also obtain a characterization of the interval order (resp. linear-interval order, linear-simiorder) dimension.
Because a binary relation's hybrid order dimension is less than its (linear) order dimension, these results will be able to 
improve known results in graph theory and computer science by identifying more efficient algorithms.
\end{abstract}

\keywords{Graph\and Extension theorems\and Linear-interval order\and Linear-semiorder\and  Realizer\and Dimension}

\section{Introduction} 
Zeno’s paradox posed for first time the question of whether time should be represented by a discrete or a continuous variable.
Since any experience has some duration, 
we have become accustomed to thinking of times as either durationless instants or else collections of such instants. In this direction, it is standard to take instant time points and the precedence relation between them to define time intervals (sets of instant time points).
Russell, however, proposed to go the other way around: temporal instants should be constructed from what he calls events.
He wanted especially to derive, an instant of time (or of a point on a line)
from a period of time (or from an interval on this line). 
In his paper \cite{wie}, Wiener provides an axiomatic frame for Russell's problem in which instants
can be defined.
To do that, he defines a precedence relation $R$ defined on a set of events $X$
satisfying the following condition:
\begin{center}
$\forall x, y, z, w\in X,\ (x,y)\in R,\ (z,y)\notin R$ and $(z,w)\in R$ imply $(x,w)\in R$\ \ \ \ \ \
\ \ \ \ \ ($\star$)
\end{center}
where $(x,y)\in R$ means that $x, y\in X$ and $x$ temporally wholly precedes $y$,
i.e., every time at which $x$ exists is temporally precedent to any time at which $y$ exists. 
Russell and Wiener postulate that for each $x\in X$, $(x,x)\notin R$ holds.
We shall call statement ($\star$) the {\it Russell - Wiener axiom}.
Intuitively the formula states that if $x$ precedes $y$ and $y$ is 
simultaneous with $z$, and $z$ procedes $w$, then $x$ precedes $w$.
Fishburn \cite{fis1}, \cite{fis2} was the first to establish the name interval order for these relations.
Interval orders are important subclasses of strict partial orders that arise in graph theory, computer science, economics, psychology, biology, and scheduling problems, among others.
For example, interval orders and graph theory, along with their incomparability graphs, also known as interval graphs, serve as a natural model for the study of scheduling and preference models.
Interval orders are used in distributed computing (vector clocks and global predicate detection), concurrency theory (pomsets and 
occurrence nets), programming language semantics (fixed-point semantics), data mining (concept analysis), and other fields.
In general, the precise time of each event occurrence is usually not required 
for many applications in computer science, but what really matters is the precedence relation.
In most of these cases, the precedence relation holds for events $x$ and $y$ if $x$ ends before $y$ begin, 
and thus according to this logic, we can construct a time model in which each event corresponds to an interval representing its duration.
In this case,
two events are incomparable if their temporal durations overlap. 
By using
the Russell-Wiener axiom,
the transitivity of the precedences   
axiom and the notion of overlapping intervals allow us to infer information regarding the sequence of events. 
Using the Russell-Wiener axiom, the transitivity of precedences axiom, 
and the concept of overlapping intervals, we can infer information about the sequence of events.
To make it more understandable, let's see the following example from the field of computer science:
In scheduling 
modeled by precedence constraints,
we have a number of tasks, say,  $t_{_1}, t_{_2},...,t_{_n}$ that must be executed by a number of parallel 
processors $p_{_1}, p_{_2},...,p_{_n}$.  
We assume that all processors are identical and that all tasks are known 
in advance and can be performed independently of one another.
Each assignment of tasks to processors is called a {\it schedule}. 
The sum of the processing times of the
tasks, assigned to a processor, is the {\it load} of this processor and,  the maximum load of any processor is the {\it length}
 of the schedule.
Our
strategy is an optimal schedule, that is, a schedule of minimal length.
When the precedence constraint is an interval order,
Papadimitriou and Yannakakis \cite{PY} showed that if tasks are placed in a list sorted by 
non-increasing size of successor sets,\footnote{The successor set of a task $t_{_i}$, $i\in \{1,2,...,n\}$ be the set of tasks that can not start before $t_{_i}$ is finished.}
and whenever a processor becomes idle, it executes the leftmost unscheduled task in the list that is ready for execution, then 
a schedule of minimal length is obtained (see also \cite[Page 3]{sta}).
Finally, recall that if an interval order $R$ represents the time intervals for a given set of tasks, the width 
(the maximum size of an antichain in $R$) gives an upper bound on how many tasks can be executed concurrently. 
This has applications, such as the allocation of a registry to a computer CPU.
Finally,
recall that, if an interval order $R$ represents the time intrvals for a given set of tasks, the width (the maximal size of an antichain in $R$) gives an upper bound on how many tasks are executed at the same time. This has applications, 
for example, in the register allocation on a computer CPU. 

On the other hand, it is well known that graph is a very useful tool to model problems in all
areas of our life.
A graph $G=(V,E)$ is the intersection graph for a non-empty family $\mathcal{F}$ of geometric objects
if there is a one-to-one correspondence between $\mathcal{F}$ and $V$ 
such that two
geometric objects in $\mathcal{F}$ have non-empty intersection if and only if the corresponding vertices 
are adjacent in $V$. 
Such a family of geometric objects is called an intersection representation of the graph. 
One of the most important intersection graphs are that of intervals on the real line and that of
triangles defined by a point on a horizontal line and an interval or a unit interval on another horizontal line. 
Intersection graphs have natural applications in several fields, including bioinformatics and involving the 
physical mapping of DNA and the genome reconstruction.

A partially ordered set or poset, $(X,\prec)$, consists of a set $X$ together with an irreflexive and transitive binary relation $\prec$ on it. A realizer of a poset $(X,\prec)$ is a family of linear orders on $X$ whose intersection is
the binary relation $R$. Szprilrajn \cite{szp} first proved that a realizer for a partial order $R$ always exists.
Dushnik and Miller \cite{DM} defined the order dimension $dim(R)$ of a poset $(X,R)$ to be the minimum cardinality of a realizer. 
The concept of order dimension plays a role which in many instances is analogous to the chromatic number for graphs.
Spinrad \cite{spi} is of the opinion that order dimension is a parameter that, in some sense, measures 
the complexity of a partial order.
In fact, various problems may be easier to solve when restricted to partial orders of small order dimension.
There are efficient algorithms to test if a partial order has order dimension at most 2. In 1982, 
Yannakakis \cite {yan} showed that for $k\geq 3$ to test if a partial order has order dimension $\leq k$ is $NP$-complete.
Generally speaking,
dimension seems to be a particularly hard $NP$-complete problem. 
This is indicated by the lack of heuristic or approximation algorithms for producing realizers with 
reasonable size for partial orders (for details see  \cite{FMP}, \cite{EPL}, \cite{MS}, \cite{spi}, \cite{sta} , \cite{tro1}, \cite{yan}).
Interval order dimension and semiorder dimension of a poset $(X,R)$, denoted
$idim(R)$ and $sdim(R)$,
are defined analogously to the order dimension but with interval orders and semiorders instead of linear orders.
Since strict linear orders are semiorders and semiorders are interval orders, we trivially obtain that order dimension 
is an upper bound and interval dimension is a lower bound for semiorder dimension. The dimension of 
acyclic binary relations which are the intersection of orders from the same class have been 
extensively investigated.
In contrast, not much is known for dimension of acyclic binary relations 
that are the intersection of orders from different classes. Two main examples 
in this area are
linear-interval orders (resp. linear-semiorders), i.e., acyclic binary relations $R$, where $R=R_{_1}\cap R_{_2}$, 
with $R_{_1}$ being a linear order and $R_{_2}$ being an interval order (resp. semiorder).
The linear-interval (resp. linear-semiorder) dimension
is defined analogously to the order dimension but with linear-interval orders (resp. linear-semiorders) instead 
of linear orders (see \cite{mer}, \cite{rav}, \cite{tak0} and \cite{tak1}).

In this paper, we give 
three main results on: (i) the (linear-) interval order and (linear-) semiorder extensions of a binary relation; (ii) the existence of
a realizer of a (linear-) interval order and (linear-) semiordern of a binary relation; and 
(iii) the characterization of the 
(linear-) interval order and (linear-) semiorder
dimension of a binary relation. 
These results 
give an analogue of the: ($\mathfrak{i}$) Szpilran extension theory for posets\cite{szp}, 
($\mathfrak{ii}$) Dushnik and Miller \cite{DM} measure of poset complexity (order dimension)
and ($\mathfrak{iii}$) Hiraguchi \cite{hir}, Ore \cite{ore} and Milner and Pouzet \cite{MP} characterization of 
order dimension for posets, in the hybrid order case.

\section{Notations and definitions}

Let $X$ be a non-empty universal set of alternatives and
$R\subseteq X\times X$ be a binary relation on $X$.
We sometimes
abbreviate $(x,y)\in R$ as $xRy$. 
An abstract system \cite{VM} is a pair $(X, R)$, where $X$ is a set and $R$ is a binary
relation where for $x, y\in X$, $xRy$ means that $x$ dominates $y$.
We say that $R$ on $X$ is (i) {\it reflexive} if for each $x\in X$, $(x,x)\in R$; (ii)
{\it irreflexive} if we never have $(x,x)\in R$;
(iii) {\it asymmetric} if for all $x, y\in X$, $(x,y)\in R \Rightarrow (y,x)\notin R$;
(iv) {\it transitive} if for all $x,y,z\in X$, [$(x,z)\in R$ and
$(z,y)\in R$] $\Rightarrow (x,y)\in R$; 
(v)
{\it antisymmetric} if for each $x,y\in X$,
[$(x,y)\in R$ and
$(y,x)\in R$] $\Rightarrow x=y$; 
(vi) {\it total} if for each $x,y\in X$,
$x\neq y$ we have $xRy$ or $yRx$.
Let $\mathcal{B}$ be the set of binary relations on $X$.
The {\it diagonal relation} $\Delta$ on $X$ is defined by $\Delta=\{(x,x)\vert x\in X\}$. 
A {\it unary operator} $\rho$ is a mapping from $\mathcal{B}$ to $\mathcal{B}$.
Thus, given a binary relation $R$, $\rho(R)\in \mathcal{B}$  is a binary relation. 
We first define the basic unary operator for binary relations. 
Given a binary relation $R$, the {\it asymmetric part} $P(R)$ of $R$ is defined as follows:
\begin{center}
$P(R)=\{(x,y)\in X\times X\vert (x,y)\in R$ and $(y,x)\notin R\}$.
\end{center}
A {\it closure operator} is a unary operator $\varphi$ from $\mathcal{B}$ to $\mathcal{B}$ that satisfies the following three properties: for all $R, R^{\prime}\in \mathcal{B}$,
($\mathfrak{a})$ $R\subseteq \varphi(R)$ (extensiveness); ($\mathfrak{b})$ 
$R\subseteq R^{\prime}\Rightarrow \varphi(R)\subseteq \varphi(R^{\prime})$ (monotonicity) and
($\mathfrak{c})$
$\varphi(\varphi(R))=\varphi(R)$ (idempotence).
For a particular property $\mathcal{P}$, a closure operation
of $R$ is defined to be the smallest relation $R_{_0}$
that contains $R$ and has the
desired property $\mathcal{P}$.
Now, we provide two examples of closure operations. First, the 
{\it transitive closure} of a relation $R$
is denoted by $\overline{R}$, that is, for all 
$x, y\in X, (x,y)\in \overline{R}$ if there exists $m\in \mathbb{N}$ and $z_{_0},...,z_{_m}\in X$
such that $x=z_{_0}, (z_{_k},z_{_{k+1}})\in R$ for all $k\in \{0,...,m-1\}$ and $z_{_m}=y$. Clearly, 
$\overline{R}$ is transitive and, because the case $m=1$ is included, it follows that 
$R\subseteq \overline{R}$. 
Secontly, the {\it reflexive closure} of $R$ is defined as follows:
\begin{center}
$rc(R)=R\cup \Delta$.
\end{center}

The following combinations of properties are considered in the next
theorems. A binary relation $R$ on $X$ is: 
(1) {\it complete} if for all $x, y\in X$, $xRy$ or $yRx$; 
(2) {\it total} if for all $x, y\in X$ with $x\neq y$, $xRy$ or $yRx$ 
($R$ is complete if and only if it is reflexive and total);
(3) {\it a strict partial order} if
it is irreflexive and transitive; (4) a {\it partial order} if
it is reflexive, transitive and antisymmetric;
(5)  a {\it strict linear order} if
it is a total strict partial order and (6) a {\it  linear order} if
it is a total partial order;
(6) an {\it interval order} if it is a strict partial order which satisfies the
Russell-Wiener axiom; (7)  a {\it strong interval order} (see \cite[Definition 3]{JPP}) if it is 
the reflexive closure of an interval order $Q$ ($R=rc(Q)$;
(8) a {\it semitransitive order} if it is irreflexive and for all $x, y, z, w\in X$, 
$(x,y)\in R$ and $(y,z)\in R$ implies $(x,w)\in R$ or $(w,z)\in R$; 
(9) a {\it semiorder} if it is an interval order which is also a semitransitive order;
A subset $Y\subseteq X$ is an $R$-cycle if, for all $x, y\in Y$, we have $(x, y) \in \overline{R}$ and 
$(y,x) \in \overline{R}$.
We say that $R$ is {\it acyclic} if there does not exist an $R$-cycle.
A binary relation $R^{\ast}$
is an {\it extension} of a binary relation $R$ if and only if
$R\subseteq R^{\ast}$ and $P(R)\subseteq P(R^{\ast})$. 
The (interval) order dimension of an a partially ordered set $(X,\prec)$ is the least $\lambda$ such that there are $\lambda$ (interval order) linear order extensions of $\prec$ whose intersection is $\prec$.
Since 
a strict linear order is a special case of an
interval order and of a semiorder respectively,
we conclude that 
a strict linear order
extension of a binary relation $R$
is also an
interval order as well as a semiorder extension of $R$.
The converse is not true.
In the simple example which follows this can be confirmed.
\begin{example}{\rm Let $X=\{x_{_1}, x_{_2}, x_{_3}, x_{_4}\}$ be a set
and let $R_{_1}=\{(x_{_1}, x_{_2}), (x_{_3}, x_{_4})\}$ and 
$R_{_2}=\{(x_{_1}, x_{_2}), (x_{_1}, x_{_3}),$
$(x_{_2}, x_{_3})\}$
be two relations on $X$.
Then, $\widetilde{R}_{_1}=\{(x_{_1}, x_{_2}), (x_{_3}, x_{_4}),(x_{_1}, x_{_4})\}$
is an interval order extension of $R$ which is not a linear order
and 
$\widetilde{R}_{_2}=\{(x_{_1}, x_{_2}), (x_{_1}, x_{_3}),(x_{_2}, x_{_3}),(x_{_4}, x_{_3})\}$
is a semiorder extension of $R$ which is not a linear order as well.}
\end{example}
Cerioli, Oliveira and Szwarcfiter in \cite{COS} gave a common generalization of
interval order dimension and (linear) order dimension of partial order $\precsim$. 
We extend this generalization in acyclic binary relations as follows: 
An acyclic binary relation $R$
is called a {\it linear-interval order} if there exist a linear order $L$ and an interval order $Q$ such that $\overline{R}=L\cap Q$ . 
In this direction,
we call an acyclic binary relation a {\it linear-semiorder} if its 
transitive closure is the intersection of a linear order and a semiorder (see \cite{tak2}).
Suppose $\mathcal{S}=\{S_{_i}\vert i\in I\}$ be a family of geometric objects.
A graph $G=(V,E)$ is an {\it intersection graph} if we can associate $\mathcal{S}$ to $G$
such that each $S_{_i}$ is corresponded to a vertex in $V$ and $(x, y)\in E$ if and only if the $S_{_i}$ corresponding to $x$ and $y$ have non-empty intersection. 
That is, there is a one-to-one correspondence between $\mathcal{S}$ and $G$ such that two sets in $\mathcal{S}$ have non-empty intersection if and only if their corresponding vertices in $G$ are adjacent.  
Intersection graphs are very important in both theoretical as well as application point of view. An interval graph is the intersection graph of a family of intervals of the real line, called an interval model. 
Let $L_{_1}$ and $L_{_2}$ be two distinct parallel lines. A {\it permutation graph} 
is the intersection graph of a family of line segments
whose endpoints lie on two parallel lines $L_{_1}$ and $L_{_2}$.
A trapezoid graph is the intersection graph of a family of trapezoids $ABCD$, such that $AB$ is on 
$L_{_1}$ and $CD$ on $L_{_2}$.
A point-interval graph (or $PI$ graph) is the intersection graph of 
a family of triangles $ABC$, such that $A$ is on $L_{_1}$ and $BC$ is on $L_{_2}$. 
Figure 1 illustrates a PI graph.
Point-interval graphs
generalize both permutation and interval graphs and
lie between permutation and trapezoid graphs.
In fact, an acyclic binary relation $R$ is called a {\it linear-interval order} if 
for each $x\in X$ there exists a 
triangle $T(x)$ such that
\begin{center}
$x\overline{R}y$ if and only if $T(x)$ lies completely to the left of $T(y)$.
\end{center}
In fact, the ordering of the apices of the triangles gives the linear order $L$, and the bases of the triangles give an interval representation of the interval order $I$.
Let $\mathcal{K}$ be a family of geometric objects on $X$ and
let $L_{_1}$ and $L_{_2}$ be two horizontal lines in the $xy$-plane with 
$L_{_1}$ above $L_{_2}$. 
Generally speaking, a binary relation $R$ on a set $X$ is $\mathcal{K}$-{\it order} 
if for each element $x\in X$, there is a geometric object $\mathcal{K}$
between $L_{_1}$ and $L_{_2}$ so that for any two elements $x, y\in X$, 
we have $x\prec y$ in $R$ if and only if $\mathcal{K}(x)$ lies completely to the left of 
$\mathcal{K}(y)$.
The set $\{\mathcal{K}(x)\vert x\in X\}$ is called a $\mathcal{K}$ representation of $R$. 
Linear-interval orders have an triangle representation and 
Linear-semiorders have a unit triangle representation.

We say that $\mathcal{R}$ is a $(p,q)$-{\it linear}-{\it interval realizer} of $R$, if $\mathcal{R}$ is an interval realizer of $\overline{R}$ ($\overline{R}=\displaystyle\bigcap\mathcal{R}$) with $p$ elements and precisely $q$ of them are non-linear. 
In this case we say that $\mathcal{R}$ $(p,q)$-realize $R$.
We define 
$(p,q)\leq (p^{\prime},q^{\prime})$ if $(p,q)$ is lexicographically smaller than or equal to $(p^{\prime},q^{\prime})$. A {\it linear}-{\it interval dimension} of an order $R$, denoted by $lidim(R)$, is the lexicographically smallest ordered pair $(p,q)$ such that there exists a 
$(p,q)$-linear-interval realizer of $R$ (see \cite[Page 113]{COS}). 
Similarly we define the notion $(p,q)$-{\it linear}-{\it semiorder} realizer of $R$.

 \begin{figure}

\centering

 \begin{tikzpicture}
    \tikzstyle{every node}=[draw,circle,fill=white,minimum size=4pt,
                            inner sep=0pt]

    \draw (0,0)(2,0) node (1234) [label=left:$a$] {}
        -- ++(00:1.25cm) node (3214) [label=right:$b$] {}
        -- ++(270:1.25cm) node (2314) [label=right:$c$] {}
         -- ++(180:1.25cm) node (2314) [label=left:$d$] {}
      -- (1234);     

\draw (2314) -- ++(-90:1.25cm) node (4132) [label=left:$g$] {};
\draw (4132) -- ++(00:1.25cm) node (5132) [label=right:$f$] {};
\draw (5132) -- ++(90:1.25cm) node (6132) [label=right:$$] {};
\draw (1234) -- ++(-45:1.767cm) node (2314) [label=left:$$] {};

\put(68,-84.8) {\small{$(a)$}}

\put(268,-84.8) {\small{$(b)$}}

\draw  [->][ultra thick] (5,0) -- (14,0);
\draw  [->][ultra thick] (5,-2.5) -- (14,-2.5);

\path [fill=mygray] (5.5,0) -- (8,-2.48) to 
(8.5,-2.48);     
\path [fill=mygray] (8.7,0) -- (5,-2.48) to 
(12,-2.48);    
\begin{scope} 
\clip (5.5,0) -- (8,-2.48) to 
(8.5,-2.48);  
\draw [draw=none, fill=black] (8.7,0) -- (5,-2.48) to (12,-2.48); 
\end{scope} %

 \path [fill=mygray] (7,0) -- (5.3,-2.48) to 
(5.5,-2.48);     
 \begin{scope} 
\clip (5.5,0) -- (8,-2.48) to 
(8.5,-2.48);  
\draw [draw=none, fill=black] (7,0) -- (5.3,-2.48) to 
(5.5,-2.48); 
\end{scope} %

 \begin{scope} 
\clip (8.7,0) -- (5,-2.48) to 
(12,-2.48);
 \draw [draw=none, fill=black] (7,0) -- (5.3,-2.48) to 
(5.5,-2.48); 
\end{scope} %

\path [fill=mygray] (10.3,0) -- (13,-2.48) to 
(13.5,-2.48);     

\path [fill=mygray] (13.395,0) -- (9.4,-2.48) to 
(9.9,-2.48);     
 
 \begin{scope} 
	\clip (13.395,0) -- (9.4,-2.48) to 
(9.9,-2.48);     
 \draw [draw=none, fill=black]  (8.7,0) -- (5,-2.48) to (12,-2.48); 
\end{scope} %

\path [fill=mygray] (11.73,0) -- (6.4,-2.48) to 
(7.1,-2.48);    
  
   \begin{scope} 
\clip (11.73,0) -- (6.4,-2.48) to 
(7.1,-2.48);
\draw [draw=none, fill=black]  (8.7,0) -- (5,-2.48) to (12,-2.48); 
\end{scope} %

\path [fill=mygray] (10.3,0) -- (13,-2.48) to 
(13.5,-2.48);     
   \begin{scope} 
	\clip (10.3,0) -- (13,-2.48) to 
(13.5,-2.48); 
	 \draw [draw=none, fill=black] (13.395,0) -- (9.4,-2.48) to 
(9.9,-2.48);     \end{scope} %

  \put(306.3,-12.7) {$\bf{.}$}

  \put(307.5,-13.7) {$\bf{.}$}

  \put(307,-13.9) {$\bf{.}$}

  \put(306,-13) {$\bf{.}$}

  \put(306.3,-13.4) {$\bf{.}$}

  \put(306.9,-13.1) {$\bf{.}$}

 \put (152,3) {$a$};

 \put (198,3) {$b$};

 \put (244,3) {$c$};

 \put (290,3) {$g$};

 \put (336,3) {$d$};

 \put (380,3) {$f$};

 \end{tikzpicture}
\par\bigskip\par
\caption{\small {(a) A simple-triangle graph A. (b) An intersection representation of A.} } \label{fig1: }
\end{figure}

\section{Main result}

The Szpilrajn's extension theorem shows that any irreflexive and transitive binary relation has an irreflexive, transitive and total (strict linear order)
extension (see Szpilrajn \cite{szp}). 
A general result of Szpilrajn's extension theorem is the following corollary.

\begin{corollary}\label{87}{\rm A binary relation $R$ on a set $X$ has a strict linear order extension 
if and only if
$R$ is an acyclic binary relation.}
\end{corollary}
\begin{proof} To prove the necessity of the corollary, we assume that $R$ is acyclic.
Then, $\overline{R}$ is irreflexive and transitive. 
By Szpilrajn’s extension theorem $\overline{R}$ has a
strict linear order extension $R^{\ast}$. Since $R\subseteq \overline{R}$ we have that
$R^{\ast}$ is a
strict linear order extension of $R$.
To prove the sufficiency, let us assume that $R$ has a strict linear order extension $Q^{\ast}$. Then, $R$ is acyclic. Indeed, suppose to the contrary that there exist 
$x, y\in X$ such that $x\overline{R}y$ and $y\overline{R}x$. It follows that
$x Q^{\ast} x$, a contradiction to irreflexivity of $Q^{\ast}$. The last conclusion completes the proof. 
\end{proof}

Szpilrajn's result remains true if asymmetry is replaced with reflexivity and antisymmetric (see \cite[Page 64]{arr}, \cite{han}),
that is, every reflexive, transitive and antisymmetric binary relation has a linear order extension.
We generalize this result as follows:

\begin{definition}{\rm A binary relation $R$ on a set $X$ is {\it transitively antisymmetric}
if and only if $\overline{R}$ is antisymmetric.}
\end{definition}

\begin{proposition}\label{876}{\rm A binary relation $R$ on a set $X$ has a linear order extension 
if and only if
$R$ is transitively antisymmetric.}
\end{proposition}
\begin{proof} To prove the necessity of the corollary, we assume that $R$ is transitively antisymmetric.
Then, $\overline{R}$ is transitive and antisymmetric. Then, by Arrow \cite[Page 64]{arr} and Hansson \cite{han}, 
$\overline{R}$ has a linear order extension. Therefore, $R$ has a linear order extension.
To prove the sufficiency suppose that $R$ has a linear order extension.
If $R$ is not transitively antisymmetric, then there are $x, y\in X$ such that
$(x,y)\in \overline{R}$, $(y,x)\in \overline{R}$ and $x\neq y$. But then,
$(x,y)\in Q$, $(y,x)\in Q$ and $x\neq y$ which
is impossible by the antisymmetry of $Q$. The last contradiction shows that $R$ is transitively antisymmetric.
\end{proof}

\par
To continue the study on the interval order dimension
let us make the following assumption.
\par\noindent
{\it Negative interval order assumption}. 
Let a binary relation $R$ on $X$ be given. Then, there exists
$x,y, a, b\in X$ such that
$(x,a)\in R$, $(b,y)\in R$, $(b,a)\notin \overline{R}$ and $(x,y)\notin R$ hold.
The set 
\begin{center}
$\mathfrak{D}_{_R}=\{((x,y),(a,b))\in X^2\times X^2\vert \ 
(x,a)\in R, (b,y)\in R,\ (b,a)\notin \overline{R} \ {\rm and}\  (x,y)\notin R\}$
\end{center}
is called the {\it negative interval order assumption set with respect to} $R$.

\par\smallskip\par\noindent
{\it Negative semiorder assumption}. 
Let a binary relation $R$ on $X$ be given. Then, there exists
$x,y, z, w\in X$ such that
$(x,y)\in R$, $(y,z)\in R$, $(x,w)\notin R$ and $(w,y)\notin R$ hold.

\begin{remark}\label{slem}{\rm If a binary relation $R$
is assumed to satisfy the negative interval order assumption
generalizes the 2+2 rule and if
it is assumed to satisfy the semiorder assumption
is equivalent to fulfil the 3+1 rule.
In this paper, we use the first notation which is more convenient for presentation of proofs.
}
\end{remark}

 \begin{figure}

\centering

 \begin{tikzpicture}
[xscale=1,yscale=1]

    \tikzstyle{every node}=[draw,circle,fill=white,minimum size=4pt,
                            inner sep=-2pt]

    \draw (0,0)(-3,0) node (1234) [label=left:\!\!$b$] {$\bullet$}
        -- ++(90:1.25cm) node (3214) [label=left:\!\!$a$] {$\bullet$} 
   -- (1234);

\draw (0,0)(-2,0) node (3334) [label=right:\ $d$] {$\bullet$}
        -- ++(90:1.25cm) node (3214) [label=right:\ $c$] {$\bullet$} 
   -- (3334);

\draw (2,2)(2,0) node (3334) [label=right:\ $z$] {$\bullet$}
        -- ++(90:1.25cm) node (3214) [label=right:\ $y$] {$\bullet$} 
         -- ++(90:1.25cm) node (3214) [label=right:\ $x$] {$\bullet$}
   -- (3334);

\put(14,33) {\small{$\bullet$}}

\put(5,33) {\small{$w$}}

\put(-77,-20) {\small{$(a)$}}

\put(32,-20) {\small{$(b)$}}

 \end{tikzpicture}
\par\bigskip\bigskip\par
\caption{\small {$R$  satisfies the negative interval order assumption (or the 2+2 rule) iff a restriction of it is isomorphic to (a) and it satisfies the negative semiorder assumption (or 3+1 rule)
iff a restriction of it is isomorphic to either (a)  or (b).} } \label{fig2: }
\end{figure}

\begin{lemma}\label{43} {\rm Let $R$ be an acyclic binary relation on a set $X$, which does not satisfy the 
negative interval order assumption. Then, $\overline{R}$ is an
interval
order extension of $R$ (not necessarily strict linear order).}
\end{lemma}
\begin{proof} By definition, $R\subseteq \overline{R}$ and
$\overline{R}$ is transitive. 
Since $R$ is acyclic, we also have that
$\overline{R}$ is irreflexive. 
To complete the proof, we have only to verify that
$\overline{R}$
satisfies the Russell-Wiener axiom. Indeed, since
$R$ does not satisfy the 
negative interval order assumption, we are led to conclude that
for all $x, y, a, b\in X$, which satisfy
$xRa, bRy$ and $(b,a)\notin \overline{R}$, we have $(x,y)\in R$.
Let now
$z, w, c, d\in X$ such that 
$z\overline{R}c, d\overline{R}w$ and $(d,c)\notin \overline{R}$.
Then, there exist natural numbers $\mu, \nu$ and alternatives $s_{_1}, s_{_2},...,s_{_{\mu}}$,
$t_{_1}, t_{_2},...,t_{_{\nu}}$
such that 
\begin{center}
$zRs_{_1}Rs_{_2}...Rs_{_{\mu}}Rc$\ {\rm and}\ 
$dRt_{_1}Rt_{_2}...Rt_{_{\nu}}Rw$.
\end{center}
But then, $s_{_{\mu}}Rc$, $dRt_{_1}$ and $(d,c)\notin \overline{R}$ imply that 
$(s_{_{\mu}},t_{_{1}})\in R$. It follows that $(z,w)\in\overline{R}$. Hence,
$\overline{R}$ is an interval order extension of $R$.
\end{proof}

\begin{theorem}\label{a11} A binary relation $R$ on a set $X$ has an
interval
order extension (not necessarily a strict linear order)
if and only if $R$ is acyclic.
\end{theorem}
\begin{proof} 
Let us prove the necessity of the theorem.
We assume that $R$ is an acyclic binary relation defined on a set $X$.
If $R$ is an interval order (if $x, y, a, b\in X$ such that 
$xRa, bRy, (b,a)\notin R=\overline{R}$, then $(x,y)\in R$), then there is nothing to prove.
Otherwise, 
$\mathfrak{D}_{_R}\neq\emptyset$. That is,
there exists $x, y, a, b\in X$ such that 
$xRa, bRy, (b,a)\notin \overline{R}$ and $(x,y)\notin R$.
We put
\begin{center}
$R^{\prime}=R\cup\{(x,y)\in X\times X \vert \ \exists 
\ a, b\in X$  such that
$xRa$,\ $bRy$\ {\rm and}\  $(b,a)\notin \overline{R}\ \}$.
\end{center}
Clearly, $R^{\prime}$ is irreflexive and $R\subset R^{\prime}$.
To verify that $R^{\prime}$
is acyclic, take any $z\in X$ and suppose that $(z,z)\in \overline{R}$. Then, there exists a natural number $m$ and alternatives $x_{_1}, x_{_2},...,x_{_m}$ such that 
\begin{center}
$z=x_{_1}R^{\prime}x_{_2}...R^{\prime}x_{_{m-1}}R^{\prime}x_{_{m}}=z$.
\end{center}
Since $R$ is acyclic, 
there is at least one $k\in \{1,...,m-1\}$ such that $(x_{_k},x_{_{k+1}})=(x,y)$ with $(x,y)\in R^{\prime}\setminus R$.
Let $x_{_{k^{\ast}}}$ be the first occurrence of $x$ and 
let $x_{_{l^{\ast}}}$ be the last occurrence of $y$. 
Clearly, for all
$k\in \{1,...,m-1\}$,
if $(x_{_k},x_{_{k+1}})\neq (x,y)$, then $(x_{_k},x_{_{k+1}})\in R$.
Then,
\begin{center}
$y=x_{_{l^{\ast}}}Rx_{_{l^{\ast}+1}}...RzRx_{_1}...Rx_{_{k^{\ast}}}=x$
\end{center}
It follows that $(y,x)\in \overline{R}$ which jointly to 
$(x,a)\in R$ and $(b,y)\in R$ implies that $(b,a)\in \overline{R}$, yielding an absurdity.
Therefore, $R^{\prime}$ is acyclic. On the other hand, if $R^{\prime}$ does not satisfies the 
negative interval order assumption, then Lemma \ref{43} implies that $\overline{R}$ is
an interval order extension of $R$, which ends the proof of the necessity of the theorem.
Otherwise, we proceed by assuming that $R^{\prime}$ satisfies 
the 
negative interval order assumption. 

Now let 
\begin{center}
$\mathcal{E}=\{Q\subseteq X\times X\vert\ Q$ 
is an acyclic extension of $R$ which satisfies the 
negative interval order assumption$\}$.
\end{center}
We have that
$R^{\prime}\in\mathcal{E}$, so this class is non-empty.
Let ${{\mathcal{C}}}=(Q_{_\theta})_{_{\theta\in \Theta}}$ be a chain in $\mathcal{E}$ and let
$\widehat{Q}=\displaystyle\bigcup_{\theta\in \Theta}Q_{_\theta}$. Then, $\widehat{Q}\in \mathcal{E}$.
To prove it we first show that $\widehat{Q}$ is acyclic (resp. irreflexive). Take $(x,x)\in\overline{\widehat{Q}}$ (resp. $(x,x)\in\widehat{Q}$)
for some $x\in X$.
Then, since ${{\mathcal{C}}}$ is a chain, there exists an $Q_{_{\theta^{\ast}}}\in {{\mathcal{C}}}, \theta^{\ast}\in \Theta$ such that $(x,x)\in\overline{Q_{_{\theta^{\ast}}}}$ (resp. $(x,x)\in Q_{_{\theta^{\ast}}}$). This is impossible by acyclicity (irreflexivity) of $Q_{_{\theta^{\ast}}}$. Therefore, $\widehat{Q}$ is irreflexive and acyclic.
On the other hand, we assume that $\widehat{Q}$
satisfies the 
negative interval order assumption, because otherwise, 
Lemma \ref{43} implies that $\overline{\widehat{Q}}$ is
an interval order extension of $R$, which ends the proof of the necessity of the theorem.
Since $R\subset \widehat{Q}$ we have that $\widehat{Q}\in \mathcal{E}$.
Therefore,
any chain 
in $\mathcal{E}$ has an upper bound in $\mathcal{E}$ (with respect to set inclusion).
By Zorn's lemma, there is a maximal element $Q^{\ast}$ in $\mathcal{E}$.
We prove that $\overline{Q^{\ast}}$ is an interval order extension of $R$. 
Clearly, $\overline{Q^{\ast}}$ is an irreflexive and transitive extension of $R$. It remains to prove that 
$\overline{Q^{\ast}}$
satisfies the Russell-Wiener axiom.
We proceed by way of contradiction. Suppose there are
$x, y, a, b\in X$ such that 
$(x,a)\in \overline{Q^{\ast}}, (b,y)\in \overline{Q^{\ast}}, (b,a)\notin \overline{Q^{\ast}}$
and $(x,y)\notin \overline{Q^{\ast}}$.
Then, $\overline{Q^{\ast}}\supset  Q^{\ast}$ is an acyclic extension of $R$ which satisfies
the negative interval order assumption, a contradiction to the maximal character of $Q^{\ast}$. Clearly, 
in any case of the proof, the extension of the interval $R$ is not required to be a linear order.
Thus,
the last contradiction completes the necessity of the theorem.

To prove the sufficiency, let us assume that $R$ has a
not necessarily linear interval
order extension 
$Q^{\ast}$.  Then, $R$ is acyclic. Indeed, suppose to the contrary that there exist 
$x\in X$, 
a natural number $m$ and alternatives $x_{_1}, x_{_2},...,x_{_m}$ such that 
\begin{center}
$xRx_{_1}Rx_{_2}...Rx_{_{m}}Rx$.
\end{center}
Since $Q^{\ast}$ is transitive and $R\subseteq Q^{\ast}$, we have that
$x Q^{\ast} x$, a contradiction to irreflexivity of $Q^{\ast}$. The last conclusion completes the proof. 
\end{proof}

\begin{corollary}\label{117}{\rm A binary relation $R$ on a set $X$ has a strong interval order extension ((not necessarily a linear order)
if and only if
$R$ is transitively antisymmetric.}
\end{corollary}
\begin{proof} To prove the necessity of the corollary, we assume that $R$ is transitively antisymmetric.
Then, $\overline{R}\setminus \Delta$ is acyclic.
By Theorem \ref{a11}, $\overline{R}\setminus \Delta$ has a
interval order extension $R^{\ast}$. Then, we have
\begin{center}
$R\subseteq \overline{R}=(\overline{R}\setminus \Delta)\cup (\overline{R}\cap \Delta)\subseteq R^{\ast}\cup
(\overline{R}\cap \Delta)\subseteq R^{\ast}\cup \Delta$.
\end{center}
Therefore,
$Q=rc(R^{\ast})=
R^{\ast}\cup \Delta$ is a
strong interval order extension of $R$.
To prove the sufficiency, let us assume that $R$ has a strong interval order extension $Q$.
Suppose on the contrary that there are $x, y\in X$ such that 
$(x,y)\in \overline{R}$, $(y,x)\in \overline{R}$ and $x\neq y$. It follows that
$(x,y)\in Q$, $(y,x)\in Q$ and $x\neq y$ which
is impossible by the asymmetry of $Q\setminus \Delta$. The last contradiction completes the proof.
\end{proof}

\begin{theorem}\label{awer}{\rm 
A binary relation $R$ on a set $X$ is a
linear-interval
order 
if and only if $R$ is acyclic.}
\end{theorem}
\begin{proof}
To prove the necessity of the theorem let us suppose
that $R$ is an acyclic binary relation defined on a set $X$.
By Theorem \ref{a11} there exists an 
interval order extension $Q$ of $R$ ($Q$ is not necessarily a strict linear order).
Then, $\overline{R}\subseteq Q$ which implies that $Q$ is 
an 
interval order extension of $\overline{R}$.

We put
\begin{center}
$R^{\ast}=\overline{R}\cup\{(x,y)\in X\times X \vert
(y,x)\in {Q\setminus \overline{R}}\}.$
\end{center}
Since $R$ is acyclic and $Q$ is irreflexive, we have that $Q^{\ast}$ is irreflexive.
If $Q=\overline{R}$, then $R^{\ast}=\overline{R}$.
By the theorem of Szpilrajn, $\overline{R}$ has a strict linear order extension $L$.
It follows that $\overline{R}=Q\cap L$ which implies that $R$ is a 
linear-interval order.
Now suppose $Q\setminus \overline{R}\neq \emptyset$. 
We now prove that
$R^{\ast}$ is acyclic and thus it is an acyclic extension
of $\overline{R}$.
Indeed, suppose to the contrary that 
there are alternatives
$\nu, z_{_0},z_{_1},z_{_2},...,z_{_m}\in X$ such that
\begin{center}
$\nu=z_{_0}R^{\ast} z_{_1}R^{\ast}z_{_2}...R^{\ast}z_{_m}=\nu.$
\end{center}
Since $R$ is 
acyclic, there is at least one $\kappa\in \{0,1,...,m-1\}$ such that
$(z_{_\kappa},z_{_{\kappa+1}})=(x,y)$.
Let $z_{_{\kappa^{\ast}}}$ be the first occurrence of $x$ and let $z_{_{\lambda^{\ast}}}$ be the last occurrence of $y$. 
Then,
\begin{center}
$y=z_{_{\lambda^{\ast}}}R z_{_{{\lambda^{\ast}+1}}}...Rz_{_m}=\nu=z_{_0}Rz_{_1}R...Rz_{_{\kappa^{\ast}}}=x$.
\end{center}
It follows that $(y,x)\in\overline{R}$, a contradiction to $(y,x)\in Q\setminus \overline{R}$.

Suppose that
$\widetilde{{\mathcal{R}}}=\{\widetilde{R}_i\vert i\in I\}$ denote the
set of acyclic extensions of $\overline{R}$ such that 
$(x,y)\in \widetilde{R}_i\setminus \overline{R}$ if and only if $(y,x)\in Q\setminus \overline{R}$.
Since $R^{\ast}\in \widetilde{{\mathcal{R}}}$ we have that $\widetilde{{\mathcal{R}}}\neq \emptyset$.
Let ${{\mathcal{C}}}=(C_{_i})_{_{i\in
I}}$ be a chain in $\widetilde{{\mathcal{R}}}$, and let
$\widehat{C}=\displaystyle\bigcup_{i\in I}C_{_i}$. 
We prove that $\widehat{C}\in \widetilde{{\mathcal{R}}}$.
To prove that $\widehat{C}$ is acyclic suppose to the contrary that
there exixts 
$\mu, s_{_0},s_{_1},s_{_2},...,s_{_n}\in X$ such that
\begin{center}
$\mu=s_{_0}\widehat{C} s_{_1}\widehat{C} s_{_2}...\widehat{C} s_{_n}=\mu.$
\end{center}
Since ${{\mathcal{C}}}$ is a chain, there exists $i^{\ast}\in I$ such that
\begin{center}
$\mu=s_{_0}C_{_{i^{\ast}}} s_{_1}C_{_{i^{\ast}}} s_{_2}...C_{_{i^{\ast}}} s_{_n}=\mu,$
\end{center}
contradicting the acyclisity of $C_{_{i^{\ast}}}$. On the other hand, it is easy to chec that
$(x,y)\in \widehat{C}\setminus \overline{R}$ implies $(y,x)\in Q\setminus \overline{R}$.

By Zorn's lemma $\widetilde{{\mathcal{R}}}$  possesses an
element, say $\widehat{R}$, that is maximal with respect to set
inclusion.
We have two cases to consider:
$\widehat{R}$ is total or not.
If $\widehat{R}$ is total then $\widehat{R}$
is a strict linear order extension of
$\overline{R}$. Then, $\overline{R}=Q\cap \widehat{R}$.
Indeed, since $\overline{R}\subseteq Q\cap \widehat{R}$,
one needs only to prove that $Q\cap \widehat{R}\subseteq \overline{R}$.
Let to the contrary that $(x,y)\in Q\cap \widehat{R}$ and $(x,y)\notin \overline{R}$.
The $(x,y)\in Q\setminus \overline{R}$ which implies that $(y,x)\in \widehat{R}$, 
a contrsdiction to asymmetry of $\widehat{R}$ (irreflexive and transitive).
Therefore, $\overline{R}=Q\cap \widehat{R}$.

If $\widehat{R}$ is not total, then there exists $x, y\in X$ such that $(x,y)\notin \widehat{R}$
and $(y,x)\notin \widehat{R}$. It follows that 
 $(x,y)\notin \overline{R}$
and $(y,x)\notin \overline{R}$. But then, 
 $(x,y)\notin Q$
and $(y,x)\notin Q$, because otherwise  $(x,y)\notin Q\setminus \overline{R}$
or $(y,x)\notin Q\setminus \overline{R}$ which implies that 
$(y,x)\notin \widehat{R}$
or $(x,y)\notin \widehat{R}$ which is impossible.
Since $\widehat{R}$ and transitive, by the theorem of Szpilrajn there exists 
a strict linear order extension $\widehat{\widehat{R}}$ of $\widehat{R}$.
Since $(\widehat{\widehat{R}}\setminus\widehat{R})\cap Q=\emptyset$
we conclude that 
 $Q\cap \widehat{\widehat{R}}=\overline{R}$. The last conclusion shows that
 $R$ is a linear-interval binary relation.

The converse is similar to the proof of the converse of Theorem \ref{a11}.
\end{proof}

\begin{theorem}\label{awerr}{\rm 
A binary relation $R$ on a set $X$ has a
semiorder extension
if and only if $R$ is acyclic.}
\end{theorem}
\begin{proof} Let $R$ be an acyclic binary relation on $X$. By Theorem \ref{a11} has an interval
order extension $Q$ of $R$. Put
\begin{center}
$Q^{\ast}=Q\cup\{(x,w)\in X\times X\setminus \Delta \vert\ {\rm there\ exist}\ 
y, z\in X\ {\rm such\ that}$
\end{center}
\begin{center}
$\ (x,y)\in Q, (y,z)\in Q, (x,w)\notin Q\ {\rm and}\ 
(w,z)\notin Q\}=Q\cup T.$
\end{center}
Clearly, $Q^{\ast}$ is irreflexive. We prove that $Q^{\ast}$ is transitive.
Indeed, let $a, b, c\in X$ such that $(a,b)\in Q^{\ast}$ and $(c,d)\in Q^{\ast}$. Then,
we have four cases
to consider:
\par\smallskip\par\noindent
{\it Case 1.} $(a,b)\in Q$ and $(b,c)\in Q$. Then, $(a,c)\in Q\subseteq Q^{\ast}$. 
\par\smallskip\par\noindent
{\it Case 2.} $(a,b)\in Q$ and $(b,c)\in T$.
Therefore, $(a,b)\in Q$ and there exists $\kappa, m\lambda \in X$ such that 
$(b,\kappa)\in Q$, $(\kappa,\lambda)\in Q$, $(b,c)\notin Q$ and $(c,\lambda)\notin Q$.
From $(a,b)\in Q$ and $(b,\kappa)\in Q$ we have that $(a,\kappa)\in Q$.
If $(a,c)\in Q\subseteq Q^{\ast}$ we have nothing to prove. We suppose that $(a,c)\notin Q$.
Then, from $(a,\kappa)\in Q$, $(\kappa,\lambda)\in Q$, $(a,c)\notin Q$ and $(c,\lambda)\notin Q$ we conclude that $(a,c)\in Q^{\ast}$.
\par\smallskip\par\noindent
{\it Case 3.} $(a,b)\in T$ and $(b,c)\in Q$. In this case, we have $(b,c)\in Q$ and 
there exists 
$\kappa, \lambda \in X$ such that 
$(a,\kappa)\in Q$, $(\kappa,\lambda)\in Q$, $(a,b)\notin Q$ and $(b,\lambda)\notin Q$.
Since, $(\kappa,\lambda)\in Q$, $(b,c)\in Q$ and $(b,\lambda)\notin Q$ we conclude that 
$(\kappa,c)\in Q$ which jointly to $(a,\kappa)\in Q$ implies that $(a,c)\in Q\subseteq Q^{\ast}$.
\par\smallskip\par\noindent
{\it Case 4.} $(a,b)\in T$ and $(b,c)\in T$. In this case, there are $\kappa, \lambda, \nu, \mu\in X$
such that 
$(a,\kappa)\in Q$, $(\kappa,\lambda)\in Q$, $(a,b)\notin Q$, $(b,\lambda)\notin Q$
and
$(b,\mu)\in Q$, $(\mu,\nu)\in Q$, $(b,c)\notin Q$ and $(c,\nu)\notin Q$.
If $(a,c)\in Q\subseteq Q^{\ast}$, then we have nothing to prove.
Suppose that $(a,c)\notin Q$.
If $(c,\lambda)\notin Q$, then from $(a,\kappa)\in Q$, $(\kappa,\lambda)\in Q$ and
$(a,c)\notin Q$ we conclude that $(a,c)\in T\subseteq Q^{\ast}$.
Otherwise, if 
$(c,\lambda)\in Q$ the we have two subcases to consider when $(a,\mu)\in Q$ or not.
If $(a,\mu)\in Q$, then from $(m,\nu)\in Q$, $(a,c)\notin Q$ and $(c,\nu)\notin Q$ we have
$(a,c)\in T\subseteq Q^{\ast}$. On the other hand, if $(a,\mu)\notin Q$, then
$(b,\mu)\in Q$, $(a,\kappa)\in Q$ implies that $(b,\kappa)\in Q$ which jointly to $(\kappa,\lambda)\in Q$ implies that $(b,\lambda)\in Q$ which is impossible.
Therefore, in all possible cases $(a,c)\in Q^{\ast}$ which implies that $Q^{\ast}$ is transitive.

To prove that $Q^{\ast}$ is an interval order we have four cases to consider.
\par\smallskip\par\noindent
{\it Case 1.} $(a,b)\in Q$, $(c,d)\in Q$ and $(c,b)\notin Q^{\ast}\supseteq Q$.
Since $Q$ is an interval order, in this case It is clear that $(a,d)\in Q\subseteq Q^{\ast}$.
\par\smallskip\par\noindent
{\it Case 2.} $(a,b)\in Q$, $(c,d)\in T$ and $(c,b)\notin Q^{\ast}\supseteq Q$.
In this case, there are $\kappa, \lambda\in X$ such that
$(c,\kappa)\in Q$, $(\kappa,\lambda)\in Q$, $(c,d)\notin Q$ and $(d,\lambda)\notin Q$.
Then, from $(a,b)\in Q$, $(c,k)\in Q$ and $(c,b)\notin Q$ we conclude that $(a,k)\in Q$.
If $(a,d)\in Q\subseteq Q^{\ast}$, then we have nothing to prove. If $(a,d)\notin Q$, then from 
$(a,\kappa)\in Q$, $(\kappa,\lambda)\in Q$, $(a,d)\notin Q$ and $(d,\lambda)\notin Q$ we have that $(a,d)\in T\subseteq Q^{\ast}$.
\par\smallskip\par\noindent
{\it Case 3.} $(a,b)\in T$, $(c,d)\in Q$ and $(c,b)\notin Q^{\ast}\supseteq Q$.
In this case,
we have $(c,d)\in Q$ and 
there exists
$\kappa, \lambda \in X$ such that 
$(a,\kappa)\in Q$, $(\kappa,\lambda)\in Q$, $(a,b)\notin Q$, $(b,\lambda)\notin Q$
and  $(c,b)\notin Q^{\ast}\supseteq Q$.
If $(a,d)\in Q\subseteq Q^{\ast}$, then we have nothing to prove.
Let $(a,d)\notin Q$. If $(d,\lambda)\notin Q$, then 
$(a,\kappa)\in Q$, $(\kappa,\lambda)\in Q$ implies $(a,d)\in T\subseteq Q^{\ast}$.
Otherwise, $(d,\lambda)\in Q$ which jointly to $(c,d)\in Q$, $(c,b)\notin Q$, $(b,\lambda)\notin Q$ imply that $(c,b)\in T\subseteq Q^{\ast}$, a contradiction. Therefore, 
$(a,d)\in T\subseteq Q^{\ast}$.
\par\smallskip\par\noindent
{\it Case 4.} $(a,b)\in T$ and $(b,c)\in T$. In this case, there are $\kappa, \lambda, \nu, \mu\in X$
such that 
$(a,\kappa)\in Q$, $(\kappa,\lambda)\in Q$, $(a,b)\notin Q$, $(b,\lambda)\notin Q$,
$(c,\mu)\in Q$, $(\mu,\nu)\in Q$, $(c,d)\notin Q$ and $(c,\nu)\notin Q$ and $(c,b)\notin Q^{\ast}\supseteq Q$.
If $(a,d)\in Q\subseteq Q^{\ast}$, then we have nothing to prove.
Let $(a,d)\notin Q$. 
If $(d,\lambda)\notin Q$, then from $(a,\kappa)\in Q$, $(\kappa,\lambda)\in Q$ and $(a,d)\notin Q$ we conclude that $(a,d)\in T\subseteq Q^{\ast}$.
If $(d,\lambda)\in Q$, then we have two subcases to consider:
($4_{_{\mathfrak{a}}}$) $(a,\mu)\in Q$ and ($4_{_{\mathfrak{b}}}$) $(a,\mu)\notin Q$.
If $(a,\mu)\in Q$, then from $(\mu,\nu)\in Q$, $(a,d)\notin Q$ and $(d,\lambda)\notin Q$
wwe conclude that $(a,d)\in T\subseteq Q^{\ast}$.
If $(a,\mu)\notin Q$, then from $(c,m)\in Q$ and $(a,\kappa)\in Q$ we conclude that $(c,k)\in Q$. But then, $(c,k)\in Q$, $(\kappa,\lambda)\in Q$, $(c,b)\notin Q$ and $(b,\lambda)\notin Q$ implies that $(c,b)\notin Q^{\ast}$, an absurdity. Hence, $(a,d)\in  Q^{\ast}$.
Therefore, $Q^{\ast}$ is an inderval order.
If $Q^{\ast}$ does not satisy the 
negative semiorder assumption, then $Q^{\ast}$ is a semiorder extension of $R$ and the proof is over.
Otherwise, $Q^{\ast}$ satisfies the 
negative semiorder assumption.
Now, let 
\begin{center}
$\mathcal{E}=\{Q\subseteq X\times X\vert\ Q$ 
is an interval order extension of $R$ which satisfies the 
negative semiorder assumption$\}$.
\end{center}
We have that
$Q^{\ast}\in\mathcal{E}$, so this class is non-empty.
Let $\mathfrak{C}=\{({\mathcal{C}}^i_{_\eta})_{_{\eta\in H_{_i}}})\vert i\in I\}$ be the family of chains in 
$\mathcal{E}$. If 
${{\mathcal{C}^{i^{\ast}}}}\!\!\!\!\!_{_\eta}=(Q_{_j})_{_{j\in J}}$ is a chain in $\mathcal{E}$ 
such that 
$\widehat{Q}=\displaystyle\bigcup_{j\in J}Q_{_j}$
does not satisfy the negative semiorder assumption then 
$\widehat{Q}$ is a semiorder extension of $R$. Otherwise, for each $i\in I$,
$\displaystyle\bigcup_{\eta\in H_{_i}}{\mathcal{C}}^i_{_\eta}\in \mathcal{E}$ holds.
By Zorn's lemma, there is a maximal element ${\widehat{Q}}^{\ast}$ in $\mathcal{E}$.
We prove that 
${\widehat{Q}}^{\ast}$
is a semiorder extension of $R$. 
Indeed, suppose to the contrary that ${\widehat{Q}}^{\ast}$ is not a semiorder.
Then, there axist$x, y, w, z\in X$ such that 
$(x,y)\in {\widehat{Q}}^{\ast}, (y,z)\in {\widehat{Q}}^{\ast}, (x,w)\notin {\widehat{Q}}^{\ast}$
and
$(w,z)\notin {\widehat{Q}}^{\ast}$.
But then, the relation
\begin{center}
${\widetilde{Q}}^{\ast}={\widehat{Q}}^{\ast}\cup\{(x,w)\in X\times X\setminus \Delta \vert\ {\rm there\ exist}\ 
y, z\in X\ {\rm such\ that}
\ (x,y)\in {\widehat{Q}}^{\ast}, (y,z)\in {\widehat{Q}}^{\ast}, (x,w)\notin {\widehat{Q}}^{\ast}\ {\rm and}\ 
(w,z)\notin {\widehat{Q}}$
\end{center}
belongs to $\mathcal{E}$, a contradiction to the maximal character of 
${\widehat{Q}}^{\ast}$.
Therefore, ${\widehat{Q}}^{\ast}$ is a semiorder extension of $R$.
The converse is evident.
\end{proof}

The following theorem is proved in a similar way to the proof of Theorem
\ref{awer}.

\begin{theorem}\label{awer1}{\rm 
A binary relation $R$ on a set $X$ is a
linear-semiorder
if and only if $R$ is acyclic.}
\end{theorem}

\section{Hybrid order dimension}

Nowadays, dimension theory is a strong brance in the graph theory and computer science. 
This is documented  by the recent book of Trotter \cite{tro}, which gives a comprehensive survey.
The notion of dimension of a poset $(X,\prec)$ was introduced in a seminal paper by Dushnik and Miller \cite{DM} as the least $\lambda$ such that there are $\lambda$ linear extensions of $\prec$ whose intersection is $\prec$.
Equivalently, the dimension of $\prec$
is the dimension of the Euclidean space $\mathbb{R}^{\lambda}$ in which $(X,\prec)$ can be embedded
in such a way that $x\prec y$ if and only if the point of $x$ is below the point of $y$ with respect to 
componentwise order (see Ore \cite{ore}).
In a more general context, we often have a class $\mathcal{R}$ of objects 
e.g., acyclic binary relations, graphs, digraphs, specific kinds of them, etc.-
and
a subclass $\mathcal{C}$ of $\mathcal{R}$ such that every $R\in \mathcal{R}$ 
is either equivalent to the intersection of a number of $C_{i}\in \mathcal{C}$
or
can be embedded into a 
product $\displaystyle\prod_{i<\lambda} C_{i}$ with $C_{i}\in \mathcal{C}$ and 
$\lambda$ being a cardinal number. 
It is then natural to regard the necessary number of the $C_{\lambda}$ as a measure of complexity of $R$, called the dimension of $R$ with respect to $\mathcal{C}$ and $\mathcal{R}$.

The following theorem is a generalized result to that of Dushnik and Miller and it is a 
key result for the study of the interval order dimension.

\begin{theorem}\label{a12} {\rm Let $(X,R)$ be an abstract system. Then,
$\overline{R}$ has as realizer the set of interval order extensions of $R$ if and only if $R$ is acyclic. } 
\end{theorem}
\begin{proof}
To prove necessity, let R be an acyclic binary relation on $X$ and let $\mathcal{Q}$ be the set of all interval order extensions of $R$. 
By Theorem \ref{a11}, the family of such extensions is
non-empty.
We show that $\overline{R}=\displaystyle\bigcap_{Q\in \mathcal{Q}}Q$. Clearly, 
$\overline{R}\subseteq \displaystyle\bigcap_{Q\in \mathcal{Q}}Q$. 
Therefore, we have only to show that $\displaystyle\bigcap_{Q\in \mathcal{Q}}Q\subseteq \overline{R}$. 
Suppose to the contrary that there exists a pair $(a,b)\in \displaystyle\bigcap_{Q\in \mathcal{Q}}Q$ but $(a,b)\notin \overline{R}$. We first show that $(b,a)\notin \overline{R}$. Indeed, if we suppose, for the sake of contradiction, that $(b,a)\in \overline{R}$, then we have $(b,a)\in \overline{Q}=Q$.
This contradicts the fact that $Q$ is asymmetric (irreflexive and transitive).
Therefore, $a, b\in X$ are non-comparable with respect to $\overline{R}$. Put
\begin{center}
$R^{\prime}=\overline{R}\cup\{(b,a)\}$
\end{center}
It is easy to check that $R^{\prime}$ is acyclic ($(a,b)\notin \overline{R}$).
By theorem \ref{a11}, $R^{\prime}$ has an interval order extension
$Q^{\ast}$. Therefore, $R$ has an interval order extension
$Q^{\ast}$ such that $(b,a)\in Q^{\ast}$, a contradiction to
asymmetry of $Q^{\ast}$ ($(a,b)\in \displaystyle\bigcap_{Q\in \mathcal{Q}}Q\subseteq Q^{\ast}).$
The last contradiction proves that $\overline{R}=\displaystyle\bigcap_{Q\in \mathcal{Q}}Q$.

To prove 
the sufficiency of the theorem, let $\overline{R}=\displaystyle\bigcap_{Q\in \mathcal{Q}}Q$, where 
$\mathcal{Q}$ is a family of interval order extensions of $R$. Then, 
$R$ is acyclic. Indeed, suppose to the contrary that there are alternatives
$x, x_{_0}, x_{_1},..., x_{_n}\in X$ such that
\begin{center}
$x=x_{_0}Rx_{_1}R...Rx_{_n}=x$.
\end{center}
Since $Q$ is a transitive extension of $R$ we have $xQx$, a contradiction to irreflexivity of $Q$.
Therefore, $R$ is acyclic. The last conclusion completes the proof.
\end{proof}

The following corollary is a consequence of Theorem \ref{a12}.

\begin{corollary} \label{a212} {\rm Let $(X,R)$ be an abstract system. Then,
$\overline{R}$ has as realizer the set of strong interval order extensions of $R$ if and only if $R$ is reflexive and transitively antisymmetric.}
\end{corollary}
\begin{proof} To prove the necessity, let $R$ be reflexive and transitively antisymmetric.
Then, $\overline{R}\setminus \Delta$ is acyclic.
By Theorem \ref{a12}, we have that $\overline{R}\setminus \Delta=\displaystyle\bigcap_{Q\in \mathcal{Q}}Q$,
where $Q$ is an interval order. Therefore, $\overline{R}=\displaystyle\bigcap_{Q\in \mathcal{Q}}rc(Q)$
where $rc(Q)$ is a strong interval order.
Conversely, suppose that $\overline{R}$ has as realizer the set $\mathcal{Q^{\ast}}$ of strong interval order extensions of $R$.
If $Q^{\ast}\in \mathcal{Q^{\ast}}$, then $Q^{\ast}\setminus \Delta$ is an interval order.
If we suppose that $R$ is not transitively antisymmetric,
then we conclude that $Q^{\ast}\setminus \Delta$ is not asymmetric which is a contradiction.
Therefore, $R$ is transitively antisymmetric. On the other hand, since
$\Delta\subseteq \displaystyle\bigcap_{Q^{\ast}\in \mathcal{Q^{\ast}}}Q^{\ast}=\overline{R}$,
we have that for all $x\in X$ there holds $(x,x)\in \overline{R}$.
Thus, here are alternatives
$x, x_{_0}, x_{_1},..., x_{_n}\in X$ such that
\begin{center}
$x=x_{_0}Rx_{_1}R...Rx_{_n}=x$.
\end{center}
Since
$R$ is transitively antisymmetric, we conclude that $x=x_{_0}=x_{_1}=...=x_{_n}$ which implies that
$(x,x)\in R$. Hence, $R$ is reflexive.
\end{proof}

Moreover, if $R$ is transitive, then
as immediate consequences of Theorem \ref{a12} and Corollary \ref{a212} we have the following results.

\begin{corollary}\label{a72}
A binary relation $R$ has as realizer the set of its interval order extensions if and only if $R$ is a strict partial order.  
\end{corollary}

\begin{corollary}\label{a721}
A binary relation $R$ has as realizer the set of its strong interval order extensions if and only if $R$ is a partial order.  
\end{corollary}

The following result is 
a generalization of 
the theorem of Dushnik and Miller \cite{DM}.

\begin{theorem}\label{a572} Let $(X,R)$ be an abstract system. Then,
$\overline{R}$ has as realizer the set of strict linear order extensions of $R$ if and only if $R$ is acyclic.   
\end{theorem}
\begin{proof}
Let R be an acyclic binary relation on $X$.
Then, $(X,\overline{R})$ is a poset. By (\cite[Theorem 2.32]{DM}
we have that the family
$\mathcal{Q}$ of
strict linear order extensions of $\overline{R}$ is a realizer of $\overline{R}$. That is,
$\overline{R}=\displaystyle\bigcap_{Q\in \mathcal{Q}}Q$. 
Since $R\subseteq \overline{R}$ and $R\subseteq Q$ imply $\overline{R}\subseteq \overline{Q}=Q$,
we have that the family of strict linear order extension of $\overline{R}$ coincides with
the family of strict linear order extension of $R$.

Conversely, suppose that $\overline{R}$ has as realizer the set of strict linear order extensions
of $R$, $\mathcal{Q}$. Then,
$\overline{R}=\displaystyle\bigcap_{Q\in \mathcal{Q}}Q$. Since 
$\displaystyle\bigcap_{Q\in \mathcal{Q}}Q$ is irreflexive, we conclude that $\overline{R}$ is acyclic.
\end{proof}

By analogy of the proof of Corollary \ref{a212} from Theorem \ref{a12}, we can prove the following
corollary from Theorem \ref{a572}.

\begin{corollary}\label{b572} Let $(X,R)$ be an abstract system. Then,
$\overline{R}$ has as realizer the set of linear order extensions of $R$ if and only if $R$ is 
reflexive and transitively antisymmetric.   
\end{corollary}

Moreover, if $R$ is transitive, then
as immediate consequences of Theorem \ref{a572} and Corollary \ref{b572} we have the following results.

\begin{corollary}\label{d572} Let $(X,R)$ be an abstract system. Then,
$R$ has as realizer the set of strict linear order extensions of $R$ if and only if $R$ is transitive and asymmetric.   
\end{corollary}

\begin{corollary}\label{e572} Let $(X,R)$ be an abstract system. Then,
$R$ has as realizer the set of linear order extensions of $R$ if and only if $R$ is 
reflexive, thansitive and antisymmetric.   
\end{corollary}

The following two theorems are proved in a similar way to the proof of Theorem \ref{a12}.

\begin{theorem}\label{aa122} {\rm Let $(X,R)$ be an abstract system. Then,
$\overline{R}$ has as realizer the set of linear-interval order extensions of $R$ if and only if $R$ is acyclic. } 
\end{theorem}

\begin{theorem}\label{aa12} {\rm Let $(X,R)$ be an abstract system. Then,
$\overline{R}$ has as realizer the set of linear-semiorder extensions of $R$ if and only if $R$ is acyclic. } 
\end{theorem}

As we mentioned above,
Ore \cite{ore} defined order dimension of a poset $\mathcal{P}=(X,\prec)$ as the least cardinal $\lambda$ (see also Hiraguchi \cite{hir}) such that there is an order preserving embedding  of $(X,\prec)$ into a direct product
\begin{center}
$dpc(\mathcal{P})=\otimes\{(X,\leq_{_i})\vert i<\lambda\}=(\displaystyle\prod_{i<\lambda}X^i,<_{_{Q}})$ 
\end{center}
of $\lambda$ linear orders
$\leq_i\ (i<\lambda)$, where
$<_{_Q}$ is defined by
\begin{center}
$(x_{i})_{_{i<\lambda}}\leq_{_Q} (y_{i})_{_{i<\lambda}}$ if and only if $x_{i}\leq_i y_{i}$ holds 
for all $i<\lambda$.
\end{center}
On the other hand, Milner and Pouzet \cite{MP}
proved that the dimension of a poset $\mathcal{P}$ is equal to the least cardinal $\lambda$
such that there is an order preserving embedding  of $(X,\prec)$ into a strict direct product
\begin{center}
$spc(\mathcal{P})=\odot\{(X,<_{_i})\vert i<\lambda\}=(\displaystyle\prod_{i<\lambda}X^i,<_{_S})$ 
\end{center}
of $\lambda$ strict linear orders
$\mathcal{<}_{_i}\ (i<\lambda)$, where $<_{_S}$ is defined by
\begin{center}
$(x_{i})_{_{i<\lambda}}<_{_S} (y_{i})_{_{i<\lambda}}$ if and only if $x_{i}<_i y_{i}$ holds for all $i<\lambda$.
\end{center}

 In order to give general results concerning those of (interval) order dimension, we extend
 the notions of order preserving embedding, componentwise order and (strict) direct product of a partial order to arbitrary binary relations.

In the following,
for the sake of maintaining uniformity of notations, for any  abstract system $(X,R)$ we denote 
$<_{_R}=P(R)$ and $\leq_{_{R}}=P(R)\cup \Delta=rc(P(R))$. Clearly, if $R$ is acyclic, then
$<_{_R}=R$ and $\leq_{_{R}}=rc(R)$.

\begin{definition}\label{a124}
A mapping from an abstract system $(X,R)$ to an abstract system $(X^{\prime},R^{\prime})$ is called an {\it dominance-preserving embedding} if it respects the dominance relation, that is, all $x, y\in X$ are mapped to 
$x^{\prime}, y^{\prime}\in R^{\prime}$ such that $xRy$ if and only if $x^{\prime}R^{\prime}y^{\prime}$.
Let $\lambda\in \aleph$ be a cardinal number and let $\mathfrak{R}=\{(X_i,R_{_i})\vert i<\lambda\}$ be a family of abstract systems.
The {\it strict componentwise dominance relation} of $\mathfrak{R}$ is a binary relation $S(\mathfrak{R})$ on the cartesian product 
$\displaystyle\prod_{i<\lambda}X^{i}$
such that given $(x_{_i})_{_{i<\lambda}}$,
$(y_{_i})_{_{i<\lambda}}\in
\displaystyle\prod_{i<\lambda}X^{i}$, we have
\begin{center}
$(x_{_i})_{_{i<\lambda}}<_{_{S(\mathfrak{R})}}
(y_{_i})_{_{i<\lambda}}$
if and only if $x_{_i}<_{_{R_{_i}}} y_{_i}$ for all
$i<\lambda$.
\end{center}
The {\it componentwise dominance relation} 
of $\mathfrak{R}$ is a binary relation $Q(\mathfrak{R})$ on 
the cartesian product 
$\displaystyle\prod_{i<\lambda}X^{i}$
such that given $(x_{_i})_{_{i<\lambda}}$,
$(y_{_i})_{_{i<\lambda}}\in
\displaystyle\prod_{i<\lambda}X^{i}$, we have
\begin{center}
$(x_{_i})_{_{i<\lambda}}\leq_{_{Q(\mathfrak{R})}}
(y_{_i})_{_{i<\lambda}}$ if and only if  $x_{_i}\leq_{_{R_{_i}}} y_{_i}$
for each $i<\lambda$.
\end{center}
The {\it strict direct product} of a family $\mathfrak{R}=\{(X_i,R_{_i})\vert i<\lambda\}$
of abstract systems,
denoted by $\odot\{(X,\mathcal{R}_{_i})\vert i<\lambda\}$,
is the Cartesian product 
$\displaystyle\prod_{i<\lambda}X^i$ equipped with the strict
componentwise dominance relation $<_{_{S(\mathfrak{R})}}$. In this case, we write 
$(\widetilde{X},<_{_{S(\mathfrak{R})}})=\odot\{(X,\mathcal{R}_{_i})\vert i<\lambda\}$ where
$\widetilde{X}=\displaystyle\prod_{i<\lambda}X^i$.
The {\it direct product} 
of a family $\mathfrak{R}=\{(X_i,R_{_i})\vert i<\lambda\}$
of abstract systems,
denoted by $\otimes\{(X,\mathcal{R}_{_i})\vert i<\lambda\}$,
is the Cartesian product 
$\displaystyle\prod_{i<\lambda}X^i$ equipped with the
componentwise dominance relation
$<_{_{Q(\mathfrak{R})}}$. 
In this case, we write 
$(\widetilde{X},<_{_{Q(\mathfrak{R})}})=\otimes\{(X,\mathcal{R}_{_i})\vert i<\lambda\}$ where
$\widetilde{X}=\displaystyle\prod_{i<\lambda}X^i$.
\end{definition}
In case of (strict) partial orders, the notions of dominance-preserving embedding, componentwise dominance relation and  (strict) direct product of an abstract system coincide with
the notions of order-preserving embedding, componentwise order and  (strict) direct product of linearly
ordered sets, respectively.

We now extend the notion of order dimension in order to study the problem of
(interval) order dimension in a general form.

\begin{definition} \label{a125} Let $\mathfrak{R}=(X,R)$ be an abstract system.
The ({\it interval order dimension}) {\it order dimension} ($idim(\mathfrak{R})$) $dim(\mathfrak{R})$ of $(X,R)$ is
the least cardinal $\lambda$ such that there are $\lambda$ (interval order) strict linear order extensions of $R$ 
whose intersection is the transitive closure $\overline{R}$ of $R$.
\end{definition}

Note that this definition coincides with the classical one when $R$ is transitive.

The following theorem generalizes the well known results of Hiraguchi \cite{hir}, Ore \cite{ore} and Milner and Pouzet
\cite{MP}.

\begin{theorem}\label{a126}Let $\mathfrak{R}=(X,R)$ be an abstract system where $R$ is acyclic. Then the following statements are equivalent.
\par\noindent
(a) The order dimension of $\mathfrak{R}$ is the least cardinal $\lambda$ such that $\overline{R}$ is the intersection of $\lambda$ strict linear orders.
\par\noindent
(b) The order dimension of $\mathfrak{R}$
is the least cardinal $\lambda$ such that there is an dominance-preserving embedding of $(X,\overline{R})$ 

\! \!\!
into a strict direct product 
of $\lambda$ strict linear orders.
\par\noindent
(c) The order dimension of $\mathfrak{R}$
is the least cardinal $\lambda$ such that there is a dominance-preserving embedding of $(X,\overline{R})$ 

\! \!\!
into a 
direct product 
of $\lambda$ linear orders.
\end{theorem}
\begin{proof}{\bf Step 1} ($dim(\mathfrak{R})\geq spc(\mathfrak{R})$).
Suppose that $\mathfrak{R}=(X,R)$ has order dimension $\lambda$. Therefore, 
$\overline{R}=\displaystyle\bigcap_{i<\lambda}\widehat{\mathcal{L}}_{_i}$ where $\widehat{\mathcal{L}}_{_i}$ are strict linear orders
on $X$. Let $\widehat{\mathfrak{L}}=\{(X_i,\widehat{\mathcal{L}}_{_i})\vert i<\lambda\}$.
We define the map $f:(X,\overline{R})\longrightarrow (\widetilde{X},<_{_{S(\widehat{\mathfrak{L}})}})=\odot\{(X,\widehat{\mathcal{L}}_{_i})\vert i<\lambda\}$ by $f(x)=(x_{_i})_{_{i<\lambda}}$
where $x_{_i}=x$ for all $i<\lambda$.
Since the ordering $<_{_{S(\widehat{\mathfrak{L}})}}$ is defined on
$\widetilde{X}$ by
\begin{center}
$(x_{i})_{_{i<\lambda}} <_{_{S(\widehat{\mathfrak{L}})}}) (y_{i})_{_{i<\lambda}}$ if and only if $x_{i}\widehat{\mathcal{L}}_{_i} y_{i}$ holds for all $i<\lambda$,
\end{center}
we have
\begin{center}
$x\overline{R}y\Leftrightarrow {\rm (}\forall i{\rm )} \ x_{_i}\widehat{\mathcal{L}}_{_i}y_{_i}
\Leftrightarrow
{\rm (}\forall i{\rm )} \ x\widehat{\mathcal{L}}_{_i}y 
\Leftrightarrow 
f(x) <_{_{S(\widehat{\mathfrak{L}})}} f(y)$.
\end{center}
\par\smallskip\par\noindent
{\bf Step 2} ($spc(\mathfrak{R})\geq dpc(\mathfrak{R}))$.
To prove this fact, it suffices to show that the strict direct product
$(\widetilde{X},<_{_{S(\widehat{\mathfrak{L}})}})=\odot\{(X,\widehat{\mathcal{L}}_{_i})\vert i<\lambda\}$ of the strict linear orders
$\widehat{\mathcal{L}}_{_i}$
 can be embedded into a direct
product of linear orders. Indeed, let for each $i<\lambda$, $\mathcal{L}_{_i}$ denote the ordering
on 
$\widetilde{X}=\displaystyle\prod_{i<\lambda}X^i$
defined by
\begin{center}
$(x_{_j})_{_{j<\lambda}} \mathcal{L}_{_i} (y_{_j})_{_{j<\lambda}}$ if and only if either $x_i \widehat{\mathcal{L}}_{_i} y_i$ or
$x_i=y_i$ and $y_{_\gamma}\widehat{\mathcal{L}}_{_\gamma}x_{_\gamma}$ where $\gamma=\min\{\beta\vert
x_{_\beta}\neq y_{_\beta}\}$.
\end{center}
Clearly, $\mathcal{L}_{_i}$ is reflexive ($\{\beta\vert
x_{_\beta}\neq y_{_\beta}\}=\emptyset$), antisymmetric and transitive on 
$\widetilde{X}$.
We prove that $\mathcal{L}_{_i}$ is also total on 
$\widetilde{X}$. Suppose that $x_{_i} \mathcal{L}_{_i} y_{_i}$ is false.
Since $\widehat{\mathcal{L}}_{_i}$ is total, it follows that 
\begin{center}
$[(x_{_i},y_{_i})\notin \widehat{\mathcal{L}}_{_i}\wedge x_{_i}\neq y_{_i}]\vee [(x_{_i},y_{_i})\notin 
\widehat{\mathcal{L}}_{_i}
\wedge(y_{_\gamma},x_{_\gamma})\notin \widehat{\mathcal{L}}_{_\gamma}$, $\gamma=\min\{\beta\vert
x_{_\beta}\neq y_{_\beta}\}]=$\\ 
$[(y_{_i}\widehat{\mathcal{L}}_{_i} x_{_i})]\vee [(y_{_i}\widehat{\mathcal{L}}_{_i} x_{_i})\wedge (x_{_\gamma}\widehat{\mathcal{L}}_{_\gamma} x_{_\gamma})]\vee [(y_{_i}=x_{_i})\wedge (x_{_\gamma}
\widehat{\mathcal{L}}_{_\gamma} x_{_\gamma})]=
A\vee B\vee C$.
\end{center}
In all cases $A, B$ and $C$
we have $y_{_i} \mathcal{L}_{_i} x_{_i}$.
 It follows that 
$\mathcal{L}_{_i}$ is a linear order extension of $<_{_{S(\widehat{\mathfrak{L}})}}$.
We prove that $(\widetilde{X},<_{_{S(\widehat{\mathfrak{L}})}})$ 
is embedded in the direct product
$(\widetilde{\widetilde{X}},<_{_{Q(\mathfrak{L})}})=\otimes\{(\widetilde{X},\mathcal{L}_{_i})\vert i<\lambda\}$,
where $\mathfrak{L}=\{(X_i, \mathcal{L}_{_i})\vert i<\lambda\}$ and $\widetilde{\widetilde{X}}=
\displaystyle\prod_{i<\lambda}\widetilde{X}^{^i}$.
Let $\widetilde{x}=(x_{_i})_{_{i<\lambda}}$
where $x_{_i}=x$ for all $i<\lambda$.
We claim that $(\widetilde{X},<_{_{S(\widehat{\mathfrak{L}})}})$
is embedded in the direct product $(\widetilde{\widetilde{X}},<_{_{Q(\mathfrak{L})}})$
by the mapping
$f(\widetilde{x})=(\widetilde{x}_{_i})_{_{i<\lambda}}$
where $\widetilde{x}_{_i}=\widetilde{x}$ for all $i<\lambda$. 
Indded, if $\widetilde{x} <_{_{S(\widehat{\mathfrak{L}})}} \widetilde{y}$, then $x_{_i}
\widehat{\mathcal{L}}_{_i}
y_{_i}$ 
and so $\widetilde{x} \mathcal{L}_{_i} \widetilde{y}$ for all $i<\lambda$. 
Therefore, $f(\widetilde{x}) <_{_{Q(\mathfrak{L})}} f(\widetilde{y})$. Conversely, if $f(\widetilde{x}) <_{_{Q(\mathfrak{L})}} f(\widetilde{y})$, then
$\widetilde{x}\neq \widetilde{y}$ and $\widetilde{x} \mathcal{L}_{_i} \widetilde{y}$ for all $i<\lambda$.
Therefore, either 
$x_{_i}\widehat{\mathcal{L}}_{_i}y_{_i}$ or $x_{_i}=y_{_i}$
for all $i<\lambda$.
If $x_{_i}=y_{_i}$, then there is some $\gamma<\lambda$ such that
$y_{_\gamma}\widehat{\mathcal{L}}_{_\gamma} x_{_\gamma}$
and $\gamma=\min\{\beta\vert
x_{_\beta}\neq y_{_\beta}\}$.
On the other hand, 
$\widetilde{x} \mathcal{L}_{_i} \widetilde{y}$ for all $i<\lambda$ implies that
$\widetilde{x} \mathcal{L}_{_\gamma} \widetilde{y}$ and thus
$x_{_\gamma}\widehat{\mathcal{L}}_{_\gamma} y_{_\gamma}$.
Since $\widehat{\mathcal{L}}_{_\gamma}$ is transitive, 
$x_{_\gamma}\widehat{\mathcal{L}}_{_\gamma} y_{_\gamma}$ and
$y_{_\gamma}\widehat{\mathcal{L}}_{_\gamma} x_{_\gamma}$
imply that $x_{_\gamma}\widehat{\mathcal{L}}_{_\gamma} x_{_\gamma}$,
a contradiction to irreflexivity of 
$\widehat{\mathcal{L}}_{_\gamma}$.
Therefore, $x_{_i}\widehat{\mathcal{L}}_{_i}y_{_i}$ for all $i<\lambda$. It follows that 
$\widetilde{x} <_{_{S(\widehat{\mathfrak{L}})}} \widetilde{y}$.
The last conclusion shows that
\begin{center}
$\widetilde{x}<_{_{S(\widehat{\mathfrak{L}})}}\widetilde{y}\Leftrightarrow 
f(\widetilde{x}) <_{_{Q(\mathfrak{L})}} f(\widetilde{y})$
\end{center}
\par\smallskip\par\noindent
{\bf Step 3} ($spc(\mathfrak{R})\geq dim(\mathfrak{R}))$.
Suppose that $spc(\mathfrak{R})=\lambda$. By definition,
$\lambda$ is
the least cardinal such that there is a dominance-preserving embedding $f$ of $(X,\overline{R})$ into a direct product 
$(\widetilde{X},<_{_{Q(\mathfrak{M})}})=\otimes\{(X,\mathcal{M}_{_i})\vert i<\lambda\}$,
where each $\mathcal{M}_{_i}$ is a linear order, 
$\mathfrak{M}=\{\mathcal{M}_{_i}\vert i<\lambda\}$ 
and
$<_{_{Q(\mathfrak{M})}}$ is defined by
\begin{center}
$(x_{i})_{_{i<\lambda}}\leq_{_{Q(\mathfrak{M})}}(y_{i})_{_{i<\lambda}}$ if and only 
if $x_{i}\leq_{_{\mathcal{M}_{_i}}} y_{i}$ hold for all $i<\lambda$.
\end{center}
Then, by supposition we have 
\begin{center}
$x\overline{R}y\Leftrightarrow f(x)<_{_{Q(\mathfrak{M})}}f(y)$.
\end{center}
If $f(x)=(x_{_i})_{_{i<\lambda}}$, we write $f_{_i}(x)=x_{_i}$. Then,
for each $i<\lambda$ define a linear order $\mathcal{C}_{_i}$ on $X$ by 
\begin{center}
$x \mathcal{C}_{_i} y$ if and only if either $f_{_i}(x)\neq f_{_i}(y)$ and 
$f_{_i}(x)\leq_{_{\mathcal{M}_{_i}}}f_{_i}(y)$ hold or 
$f_{_i}(x)=f_{_i}(y)$ and $f_{_j}(y)\leq_{_{\mathcal{M}_{_j}}}f_{_j}(x)$, where \\
$j=min\{k<\lambda\vert f_{_k}(x)\neq f_{_k}(y)\}$. 
\end{center}
We prove that
\begin{center}
 $f(x)<_{_{Q(\mathfrak{M})}}f(y)\Leftrightarrow (\forall i<\lambda) (x\mathcal{C}_{_i}y)$. 
 \end{center}
 Indeed, let $f(x)<_{_{Q(\mathfrak{M})}}f(y)$, then $x\neq y$ and thus for all $i<\lambda$, we have 
$f_{_i}(x)\leq_{_{\mathcal{M}_{_i}}} f_{_i}(y)$ and $f_{_i}(x)\neq f_{_i}(y)$. 
Therefore, for all $i<\lambda$ there holds $x \mathcal{C}_{_i} y$. Hence, 
\begin{center}
$f(x)<_{_{Q(\mathfrak{M})}}f(y)\Rightarrow (\forall i<\lambda) (x\mathcal{C}_{_i}y)$.
\end{center}
Converselly, let $x \mathcal{C}_{_i} y$ for all $i<\lambda$. 
Then, either
($\mathfrak{a}$) $f_{_i}(x)\neq f_{_i}(y)$ and $f_{_i}(x)\leq_{_{\mathcal{M}_{_i}}}f_{_i}(y)$ or ($\mathfrak{b}$)
$f_{_i}(x)=f_{_i}(y)$ and $f_{_j}(y)\leq_{_{\mathcal{M}_{_j}}}f_{_j}(x)$, where
$j=min\{k<\lambda\vert f_{_k}(x)\neq f_{_k}(y)\}$ hold.
Suppose that $f_{_i}(x)=f_{_i}(y)$ for some $i<\lambda$. Then,
$f_{_j}(y)\leq_{_{\mathcal{M}_{_j}}}f_{_j}(x)$, where
$j=min\{k<\lambda\vert f_{_k}(x)\neq f_{_k}(y)\}$.
Since $x \mathcal{C}_{_j} y$ and $f_{_j}(x)\neq f_{_j}(y)$ we have 
$f_{_j}(x)\leq_{_{\mathcal{M}_{_j}}}f_{_j}(y)$. By antisymmetry of $\leq_{_{\mathcal{M}_{_j}}}$ we have
$f_{_j}(x)=f_{_j}(y)$ which is impossible by the definition of $j$. The last contradiction shows that
for all $i<\lambda$ we have $f_{_i}(x)\neq f_{_i}(y)$ and $f_{_i}(x)\leq_{_{\mathcal{M}_{_i}}}f_{_i}(y)$.
It follows that 
\begin{center}
$(\forall i<\lambda) (x\mathcal{C}_{_i}y)\Rightarrow
(\forall i<\lambda)
[(f_{_i}(x)\neq f_{_i}(y))\wedge (f_{_i}(x)\leq_{_{\mathcal{M}_{_i}}}f_{_i}(y))]
\Rightarrow 
(\forall i<\lambda) (f_{_i}(x)<_{_{\mathcal{M}_{_i}}}f_{_i}(y))
\Rightarrow 
(\forall i<\lambda) (x_{_i}<_{_{\mathcal{M}_{_i}}}y_{_i})
\Rightarrow
(x_{_i})_{_{i<\lambda}}
<_{_{Q(\mathfrak{M})}}
(y_{_i})_{_{i<\lambda}}
\Rightarrow
f(x)<_{_{Q(\mathfrak{M})}}f(y)$.
\end{center}
The last conclusion implies that,
\begin{center}
 $f(x)<_{_{Q(\mathfrak{M})}}f(y)\Leftrightarrow (\forall i<\lambda) (x\mathcal{C}_{_i}y)$. 
 \end{center}
 Therefore,
 \begin{center}
 $x\overline{R}y\Leftrightarrow f(x)<_{_{Q(\mathfrak{M})}}f(y)\Leftrightarrow (\forall i<\lambda) (x\mathcal{C}_{_i}y)$. 
 \end{center}
Since $R$ is acyclic, the last implication implies that
$\overline{R}=
\displaystyle\bigcap_{i<\lambda}(\mathcal{C}_{_i}\setminus \Delta)$,
where for all $i<\lambda$, $\mathcal{C}_{_i}\setminus \Delta$ is a strict linear order.
As a consequence of the three steps above we conclude that
$dim(\mathfrak{R})=spc(\mathfrak{R})=dpc(\mathfrak{R})$, and the proof is complete.
\end{proof}

As an immediate consequence of Theorem \ref{a126}, we have the following 
corollary which is 
the main result of \cite{MP}.

\begin{corollary}\label{a1291} Let $\mathfrak{F}=(X,\prec)$ be a poset. Then the following statements are equivalent.
\par\noindent
(a) The order dimension of $\mathfrak{F}$ is the least cardinal $\lambda$ such that $\prec$ is the intersection of $\lambda$ strict linear orders.
\par\noindent
(b) The order dimension of $\mathfrak{F}$
is the least cardinal $\lambda$ such that there is an embedding of $(X,\prec)$ into a strict  

\! \!\!direct product of $\lambda$ strict linear orders.
\par\noindent
(c) The order dimension of $\mathfrak{F}$
is the least cardinal $\lambda$ such that there is an embedding of $(X,\prec)$ into a direct 

\! \!\!product 
of $\lambda$ linear orders.
\end{corollary}

An alternative definition of the interval order $\prec$ defined in $X$ can be made by assigning to each element $x\in X$ an open interval $I_x=(a_x,b_x)$ of the real line, such that $x\prec y$ in $X$ if and only if $b_x\leq a_y$.
Such a collection of intervals is called an {\it interval representation of} $\prec$.
Let $\lambda\in \aleph$ be a cardinal number and
let $\mathcal{I}=(I_{_i})_{_{i<\lambda}}$ be a
family of interval orders.
We denote by $\widetilde{I}_{_i}$ the interval order representation of each interval order
$I_{_i}$. Let $(a_{_x}^{i},b_{_x}^{i})$ be an interval corresponding
to $x\in X$ in the representation of $\widetilde{I}_{_i}$. With $x\in X$ we associate the box $\displaystyle\prod_{i<\lambda}(a_{_x}^{i},b_{_x}^{i})\subseteq \mathbb{R}^\lambda$. Each of these boxes is uniquely determined by its {\it upper extreme corner}
$u_x=(b_x^{i})_{_{i<\lambda}}$ and its 
{\it lower extreme corner} $l_x=(a_x^{i})_{_{i<\lambda}}$.
Such an assignment is called a {\it box embedding} of $X$.
For interval order dimension, the box embedding play the role of the point embedding into $\mathbb{R}^\lambda$
introduced by Ore. The projections of  a box embedding on each coordinate yields an 
interval order (see \cite{FHM}).

In order to approach the interval orders analogue of the 
Hiraguchi \cite{hir}, Ore \cite{ore} and Milner and Pouzet \cite{MP} results for posets, in a first step the concepts of direct product and strict direct product
have to be generalized from linear orders to interval orders
on $X$. 
The {\it direct product} of a family $\mathfrak{G}=\{(X_i,\preceq_i)\vert i<\lambda\}$
of
strong interval orders is the Cartesian product 
$\displaystyle\prod_{i<\lambda}X^i$ equipped with 
the ordering $\preceq_{_{Q(\mathfrak{G})}}$ defined by 
\begin{center}
$x\preceq_{_{Q(\mathfrak{G})}} y$ if and only if 
either
$b_{x}^{i}\leq a_{y}^{i}$ or $a_{x}^{i}=a_{y}^{i}, b_{x}^{i}=b_{y}^{i}$
holds for all $i<\lambda$.
\end{center} 
The {\it strict direct product} of a family $\mathfrak{G}=\{(X_i,\prec_i)\vert i<\lambda\}$
of
interval orders is the Cartesian product 
$\displaystyle\prod_{i<\lambda}X^i$
equipped with 
the ordering $\prec_{_{S(\mathfrak{G})}}$ defined by 
\begin{center}
$x\prec_{_{s(\mathfrak{G})}} y$ if and only if 
$b_{x}^{i}\leq a_{y}^{i}$
holds for all $i<\lambda$.
\end{center}

\begin{definition}\label{a127}{\rm Let $\mathfrak{P} = (X, R)$ be an abstract system. (i)We call $idpc(\mathfrak{P})$, the least cardinal $\lambda$ such that there is a box embedding of $(X,\overline{R})$ into a direct product of $\lambda$ strong interval orders.
(ii) 
We call $ispc(\mathfrak{P})$), the least cardinal $\lambda$ such that there is a box embedding of $(X,\overline{R})$ into a direct product of $\lambda$ interval orders.
}
\end{definition}

\begin{theorem}\label{than} Let $\mathfrak{P}=(X,R)$ be an abstract system where $R$ is acyclic. Then the following statements are equivalent.
\par\noindent
(a) The interval order dimension of $\mathfrak{P}$ is the least cardinal $\lambda$ such that $\overline{R}$ is the intersection of $\lambda$ interval orders.
\par\noindent
(b) The interval order dimension of $\mathfrak{P}$
is the least cardinal $\lambda$ such that there is a box embedding of $(X,\overline{R})$ into a strict 

\! \!\!direct product of $\lambda$ interval orders.
\par\noindent
(c) The interval order dimension of $\mathfrak{P}$
is the least cardinal $\lambda$ such that there is an box embedding of $(X,\overline{R})$ into a direct 

\! \!\!product of $\lambda$ strong interval orders.
\end{theorem}
\begin{proof}{\bf Step 1} ($idim(\mathfrak{P})\geq ispc(\mathfrak{P})$). Suppose that $\mathfrak{P}=(X,R)$ has interval order dimension $\lambda$. Therefore, 
$\overline{R}=\displaystyle\bigcap_{i<\lambda}\prec_i$ where $\prec_i$ are interval orders
on $X$. Let $\mathcal{I}=\{I_x^i \vert x\in X\}$, where $I_x^i=(a_{x}^{i},b_{x}^{i})$
be an interval representation of $\prec_{_i}$. Let also
$\widetilde{X}=\displaystyle\prod_{i<\lambda}X^i$ and $\widehat{\mathfrak{O}}=\{\prec_{i}\vert i<\lambda\}$.
We define the map $f:(X,\overline{R})\longrightarrow (\widetilde{X},<_{_{S(\widehat{\mathfrak{O}})}})=\odot\{(X,\prec_{_i})\vert i<\lambda\}$ by 
$f(x)=\displaystyle\prod_{_{i<\lambda}} (a_{x}^{i},b_{x}^{i})$.
The ordering $<_{_{S(\widehat{\mathfrak{O}})}})$ is defined by
\begin{center}
$f(x)<_{_{S(\widehat{\mathfrak{O}})}}) f(y)$ if and only if
$b_{x}^{i}\leq a_{y}^{i}$
holds for all $i<\lambda$.
\end{center}
Therefore,
\begin{center}
$x\overline{R}y\Leftrightarrow {\rm (}\forall i<\lambda{\rm )} \ [b_{x}^{i}\leq a_{y}^{i}] \Leftrightarrow f(x) <_{_{S(\widehat{\mathfrak{O}})}}) f(y)$.
\end{center}
\par\smallskip\par\noindent
{\bf Step 2} ($ispc(\mathfrak{P})\geq$ $idpc(\mathfrak{P})$).
To show this fact, it suffices to show that
the strict direct product
$(\widetilde{X},<_{_{S(\widehat{\mathfrak{O}})}})=\odot\{(X,$
$\prec_{_i})\vert i<\lambda\}$ can be box embedded into a direct
product 
of strong interval orders. Indeed,
let $\mathcal{I}=\{I_x^i \vert x\in X\}$, where for each $i<\lambda$,
$I_x^i=(a_{x}^{i},b_{x}^{i})$
be an interval representation of $\prec_{_i}$.
For each $i<\lambda$, define the ordering $\sqsubseteq_{_i}$
on $X$ by 
\begin{center}
$(x_{_j})_{_{j<\lambda}}\sqsubseteq_{_i}(y_{_j})_{_{j<\lambda}}$ if and only if either
(i) $b_{x}^{i}\leq a_{y}^{i}$ or (ii) $a_{x}^{i}=a_{y}^{i}, b_{x}^{i}=b_{y}^{i}$ and $b_{y}^{k}\leq a_{x}^{k}$\\
where 
$k=\min\{\mu\vert a_{x}^{\mu}\neq a_{y}^{\mu}\ {\rm or}\ b_{x}^{\mu}\neq b_{y}^{\mu})]$.
\end{center}
Clearly, for all $i$, $\sqsubseteq_{_i}$ is an extension of $<_{_{S(\widehat{\mathfrak{O}})}}$.
Since the reals satisfy the law of trichotomy we conclude that for each
$i<\lambda$, $\sqsubseteq_{_i}$ is a strong interval order.
We show that $(\widetilde{X},<_{_{S(\widehat{\mathfrak{O}})}})$ 
is box embedded in the direct product
$(\widetilde{\widetilde{X}},<_{_{Q(\mathfrak{O})}})=\otimes\{(\widetilde{X},\sqsubseteq_{_i})\vert i<\lambda\}$,
where $\mathfrak{O}=\{(\widetilde{X},\sqsubseteq_{_i})\vert i<\lambda\}$ and $\widetilde{\widetilde{X}}=
\displaystyle\prod_{i<\lambda}\widetilde{X}^{^i}$. Let $\widetilde{x}=(x_{_i})_{_{i<\lambda}}$,
where $x_{_i}=x$ for all $i<\lambda$.
By definition, the ordering $<_{_{Q(\mathfrak{O})}}$ is defined
\begin{center}
$(\widetilde{x}_{_j})_{_{j<\lambda}}<_{_{Q(\mathfrak{O})}}(\widetilde{y}_{_j})_{_{j<\lambda}}$
if and only if $(\widetilde{x}_{_j})_{_{j<\lambda}} \sqsubseteq_{_i}(\widetilde{y}_{_j})_{_{j<\lambda}}$
holds for all $i<\lambda$.
\end{center}
Let $f$ be the mapping
$f(\widetilde{x})=(\widetilde{x}_{_i})_{_{i<\lambda}}$,
where
$\widetilde{x}_{_i}=\widetilde{x}$
for all $i<\lambda$.
Clearly, there holds the following implication:
\begin{center}
$\widetilde{x} <_{_{S(\widehat{\mathfrak{O}})}}\widetilde{y}
\Leftrightarrow
{\rm (}\forall i<\lambda{\rm )}[b_{x}^{i}\leq a_{y}^{i}]
\Rightarrow {\rm (}\forall i<\lambda{\rm )}[\widetilde{x}\sqsubseteq_{_i}\widetilde{y}]
\Leftrightarrow f(\widetilde{x})<_{_{Q(\mathfrak{O})}} f(\widetilde{y})$.
\end{center}
Conversely, if $f(\widetilde{x})<_{_{Q(\mathfrak{O})}} f(\widetilde{y})$, then 
$\widetilde{x}\sqsubseteq_{_i}\widetilde{y}$ for all
$i<\lambda$. Therefore, for all $i<\lambda$,
\begin{center}
[$b_{x}^{i}\leq a_{y}^{i}] \vee [(a_{x}^{i}=a_{y}^{i}, b_{x}^{i}=b_{y}^{i})\wedge  (b_{y}^{k}\leq a_{x}^{k}$
where $k=\min\{\mu\vert a_{x}^{\mu}\neq a_{y}^{\mu}\ {\rm or}\ b_{x}^{\mu}\neq b_{y}^{\mu})]$
\end{center}
Suppose that $a_{x}^{i}=a_{y}^{i}$ and $b_{x}^{i}=b_{y}^{i}$ for some $i<\lambda$. Then,
there is some $k$ such that 
$b_{y}^{k}\leq a_{x}^{k}$. On the other hand, since $\widetilde{x}\sqsubseteq_{_k}\widetilde{y}$
and $a_{y}^{k}<b_{y}^{k}\leq a_{x}^{k}<b_{x}^{k}$ ($a_{x}^{k}\neq a_{y}^{k}$ and $b_{x}^{k}\neq b_{y}^{k}$),
we have that $b_{x}^{k}\leq a_{y}^{k}$. But then, $b_{y}^{k}\leq a_{x}^{k}<b_{x}^{k}\leq a_{y}^{k}$
implies $b_{y}^{k}< a_{y}^{k}$ which is impossible.
The last contradiction shows that,
for all $i<\lambda$ there holds
$b_{x}^{i}\leq a_{y}^{i}$, which implies that 
$\widetilde{x} <_{_{S(\widehat{\mathfrak{O}})}}\widetilde{y}$.
Therefore,
\begin{center}
$\widetilde{x} <_{_{S(\widehat{\mathfrak{O}})}}\widetilde{y}\Leftrightarrow 
f(\widetilde{x})<_{_{Q(\mathfrak{O})}} f(\widetilde{y})$.
\end{center}
\par\smallskip\par\noindent
{\bf Step 3} ($idpc(\mathfrak{P})\geq idim(\mathfrak{P})$).
Suppose that $idpc(\mathfrak{P})=\lambda$.
Then,
$\lambda$ is
the least cardinal such that there is a box embedding $f$ of $(X,\overline{R})$ into a direct product 
$(\widetilde{X},<_{_{Q(\mathfrak{T})}})=\otimes\{(X,\unlhd_{_i})\vert i<\lambda\}$
of strong interval orders $\{\unlhd_{_i}\vert i<\lambda\}=\mathfrak{T}$.
By definition, $\widetilde{X}=
\displaystyle\prod_{i<\lambda}X^{^i}$. On the other hand, if for all $i<\lambda$,
$\mathfrak{J}_x^i=(\alpha_{x}^{i},\beta_{x}^{i})$
is an interval representation of $\unlhd_{_i}$, then the ordering
$<_{_{Q(\mathfrak{T})}}$ is defined by 
\begin{center}
$(x_{_i})_{_{i<\lambda}}<_{_{Q(\mathfrak{T})}}(y_{_i})_{_{i<\lambda}}$
if and only if 
$\beta_{x}^{i}\leq \alpha_{y}^{i}$ holds for all $i<\lambda$.
\end{center}
Then, by definition we have 
\begin{center}
$x\overline{R}y\Leftrightarrow f(x)<_{_{Q(\mathfrak{T})}}f(y)$.
\end{center}
If $f(x)=(x_{_j})_{_{j<\lambda}}$, then
for each $i<\lambda$ we define the ordering ${\ll\!\!\!\!\!\!\!_{_{{-\!\!\!-\!\!\!-}}}}_{i}$
on $X$ by:
\begin{center}
 $x {\ll\!\!\!\!\!\!\!_{_{{-\!\!\!-\!\!\!-}}}}_{i}y$ if and only if either
$\beta_{x}^{i}\leq \alpha_{y}^{i}$ or $\alpha_{x}^{i}=\alpha_{y}^{i}, \beta_{x}^{i}=\beta_{y}^{i}$ and 
$\beta_{y}^{k}\leq \alpha_{x}^{k}$,
where 
$k=\min\{\mu\vert a_{x}^{\mu}\neq a_{y}^{\mu}\ {\rm or}\ b_{x}^{\mu}\neq b_{y}^{\mu})$.
\end{center}
Clearly, each ${\ll\!\!\!\!\!\!\!_{_{{-\!\!\!-\!\!\!-}}}}_{i}$ is a strong interval order extension of $\unlhd_{_i}$.
We prove that
\begin{center}
 $f(x)<_{_{Q(\mathfrak{T})}}f(y)\Leftrightarrow (\forall i<\lambda) (x {\ll\!\!\!\!\!\!\!_{_{{-\!\!\!-\!\!\!-}}}}_{i} y)$. 
 \end{center}
 Indeed, let $f(x)<_{_{Q(\mathfrak{T})}}f(y)$. Then, $x\neq y$ and for any $i<\lambda$ there holds 
$\beta_{x}^{i}\leq \alpha_{y}^{i}$ 
and so $x {\ll\!\!\!\!\!\!\!_{_{{-\!\!\!-\!\!\!-}}}}_{i} y$
for all $i<\lambda$. 
Converselly, let $x {\ll\!\!\!\!\!\!\!_{_{{-\!\!\!-\!\!\!-}}}}_{i} y$ for all $i<\lambda$. 
Then, either
\begin{center}
($\mathfrak{a}$) $[(a_{x}^{i}\neq a_{y}^{i})\vee (b_{x}^{i}\neq b_{y}^{i})]\wedge 
 (\beta_{x}^{i}\leq \alpha_{y}^{i})$
 \end{center}
or 
\begin{center}
($\mathfrak{b}$)
$[(a_{x}^{i}=a_{y}^{i})\wedge (b_{x}^{i}=b_{y}^{i})]\wedge [(\beta_{y}^{j}\leq \alpha_{x}^{j})$,
where
$j=min\{k<\lambda\vert a_{x}^{k}\neq a_{y}^{k}$\ or \  $b_{x}^{k}\neq b_{y}^{k}]$.
\end{center}
Suppose that 
$a_{x}^{i}=a_{y}^{i}$ and $b_{x}^{i}=b_{y}^{i}$
for some $i<\lambda$. Then,
$\beta_{y}^{j}\leq \alpha_{x}^{j}$ where $j$
has the meaning above mentioned. 
On the other hand, $\beta_{y}^{j}\leq \alpha_{x}^{j}$ implies that
$a_{x}^{j}\neq a_{y}^{j}$ or $b_{x}^{j}\neq b_{y}^{j}$.
Since
$x {\ll\!\!\!\!\!\!\!_{_{{-\!\!\!-\!\!\!-}}}}_{j} y$ we have that
$\beta_{x}^{j}\leq \alpha_{y}^{j}$.
It follows that $\beta_{y}^{j}< \alpha_{y}^{j}$ ($\alpha_{x}^{j}<\beta_{x}^{j}$)
which is impossible. The last conclusion shows that
for all $i<\lambda$ we have that the case ($\mathfrak{a}$) holds.
Therefore, 
\begin{center}
$(\forall i<\lambda) (x {\ll\!\!\!\!\!\!\!_{_{{-\!\!\!-\!\!\!-}}}}_{i} y)\Rightarrow (\forall i<\lambda)
(\beta_{x}^{i}\leq \alpha_{y}^{i}) \Rightarrow 
(x_{_i})_{_{i<\lambda}}<_{_{Q(\mathfrak{T})}}(y_{_i})_{_{i<\lambda}}\Rightarrow 
f(x)<_{_{Q(\mathfrak{T})}}f(y)$.
\end{center}
Therefore,
by combining the previous implications, we get 
\begin{center}
$(\forall i<\lambda) (x {\ll\!\!\!\!\!\!\!_{_{{-\!\!\!-\!\!\!-}}}}_{i} y)\Leftrightarrow 
f(x)<_{_{Q(\mathfrak{T})}}f(y)$.
\end{center}
Finally, by 
$x\overline{R}y\Leftrightarrow f(x)<_{_{Q(\mathfrak{T})}}f(y)$, we have that
\begin{center}
 $x\overline{R} y\Leftrightarrow f(x)<_{_{Q(\mathfrak{T})}}f(y)\Leftrightarrow (\forall i<\lambda) (x\
 {\ll\!\!\!\!\!\!\!_{_{{-\!\!\!-\!\!\!-}}}}_{i}y)$. 
 \end{center}
Since $R$ is acyclic, the last implication implies that
$\overline{R}=
\displaystyle\bigcap_{i<\lambda}({\ll\!\!\!\!\!\!\!_{_{{-\!\!\!-\!\!\!-}}}}_{i}\!\!\setminus \!\Delta)$,
where for all $i<\lambda$, ${\ll\!\!\!\!\!\!\!_{_{{-\!\!\!-\!\!\!-}}}}_{i}\!\!\setminus \!\Delta$ is an interval order.
As a consequence of the three steps above we conclude that
$dim(\mathfrak{R})=spc(\mathfrak{R})=dpc(\mathfrak{R})$, and the proof is complete.
\end{proof}

The following corollary is an immediate consequence of Theorem \ref{than}.

\begin{corollary}\label{a129} Let $\mathfrak{G}=(X,\prec)$ be a poset. Then the following statements are equivalent.
\par\noindent
(a) The interval order dimension of $\mathfrak{G}$ is the least cardinal $\lambda$ such that $\prec$ is the intersection of $\lambda$ interval orders.
\par\noindent
(b) The interval order dimension of $\mathfrak{G}$
is the least cardinal $\lambda$ such that there is a box embedding of $(X,\prec)$ into a strict  

\! \!\!direct product of $\lambda$ interval orders.
\par\noindent
(c) The interval order dimension of $\mathfrak{G}$
is the least cardinal $\lambda$ such that there is a box embedding of $(X,\prec)$ into a direct 

\! \!\!product 
of $\lambda$ strong interval orders.
\end{corollary}

Let $T$ be a triangle $ABC$. Denote $\kappa(T)=A$ and $\pi(T)=BC$. 
Let $L_{_1}$ and $L_{_2}$  be two distinct parallel lines.
A point-interval graph or PI graph is the intersection graph of a family of triangles $ABC$, such that $A$ is on $L_{_1}$ and $BC$ is on $L_{_2}$. Except for the definition we gave
in the introduction,
the linear-interval order can also be defined as follows: An acyclic binary relation $R$ is a linear-interval order if there is such a triangle $T_x$ for each element $x\in X$, and $(y,x)\in R$ if and only if $T_{_y}$ lies completely to the left of $T_{_x}$. 
In fact, 
the ordering of the apices $\kappa(T_{_x})=x$ of the triangles gives the linear order L, and the bases $\pi(T_{_x})=(a_{_x},b_{_x})$ of the triangles give an interval representation of the interval order $P$ for which $\overline{R}=L\cap P$.
As usual, the left and right extreme points of an interval $I_{_x}$ are denoted by 
$a_{_x}$ and $b_{_x}$, respectively. 
When $a_{_x}=b_{_x}=x$, we say that $I_{_x}$ is trivial. 
Let $\mathbb{R}^\lambda$ be the cartesian product of $\lambda$ many copies of $\mathbb{R}$. A {\it linear-interval point} $\gamma$ is the set 
$\displaystyle\prod_{i<\lambda}I_{_i}$ where $I_{_i}\subset \mathbb{R}$ for all $i<\lambda$
(notice that in this definition it is allowed that $I_{_i}$ be trivial).
With $x\in X$ we associate the box $\displaystyle\prod_{i<\lambda}(a_{_x}^{i},b_{_x}^{i})\subseteq \mathbb{R}^\lambda$. 
Such an assignment is called a {\it linear-interval box embedding} of $X$.
For linear-interval order dimension, the linear-interval box embedding play the role of the point embedding into $\mathbb{R}^\lambda$
introduced by Ore. 
The projections of  a linear-interval box embedding on each coordinate yields a 
linear order or an interval order.

In order to approach the linear-interval orders analogue of the 
Hiraguchi \cite{hir}, Ore \cite{ore} and Milner and Pouzet \cite{MP} results for posets, in a first step the concepts of direct product and strict direct product
have to be generalized from linear orders and interval orders to linear-interval orders
on $X$. 
The {\it direct product} of a family $\mathfrak{G}=\{(X_i,\preceq_i)\vert i<\lambda\}$
of
strong linear-interval orders is the Cartesian product 
$\displaystyle\prod_{i<\lambda}X^i$ equipped with 
the ordering $\preceq_{_{Q(\mathfrak{G})}}$ defined by 
\begin{center}
$x\preceq_{_{Q(\mathfrak{G})}} y$ if and only if 
either
$b_{x}^{i}\leq a_{y}^{i}$ or $a_{x}^{i}=a_{y}^{i}, b_{x}^{i}=b_{y}^{i}$
holds for all $i<\lambda$.
\end{center} 
The {\it strict direct product} of a family $\mathfrak{G}=\{(X_i,\prec_i)\vert i<\lambda\}$
of
linear-interval orders is the Cartesian product 
$\displaystyle\prod_{i<\lambda}X^i$
equipped with 
the ordering $\prec_{_{S(\mathfrak{G})}}$ defined by 
\begin{center}
$x\prec_{_{s(\mathfrak{G})}} y$ if and only if 
$b_{x}^{i}\leq a_{y}^{i}$
holds for all $i<\lambda$.
\end{center}

\begin{definition}\label{a127}{\rm Let $\mathfrak{P} = (X, R)$ be an abstract system. (i)We call $lidpc(\mathfrak{P})$, the least cardinal $\lambda$ such that there is a linear-interval box embedding of $(X,\overline{R})$ into a direct product of $\lambda$ strong linear-interval orders.
(ii) 
We call $lispc(\mathfrak{P})$), the least cardinal $\lambda$ such that there is a linear-interval box embedding of $(X,\overline{R})$ into a direct product of $\lambda$ linear-interval orders.
}
\end{definition}

The following theorem generalizes Theorem \ref{a126} and Theorem \ref{a1128}. The prove is omitted since it follows exactly the same scheme.

\begin{theorem}\label{aw128} Let $\mathfrak{P}=(X,R)$ be an abstract system where $R$ is acyclic. Then the following statements are equivalent.
\par\noindent
(a) The $(\lambda,\mu)$-linear-interval order dimension of $\mathfrak{P}$ is the least cardinal $\lambda$ such that $\overline{R}$ is the intersection of $\lambda$ linear-inter-
\par
\!\!\! \!\!
val orders which $\mu$ of them are not linear orders.
\par\noindent
(b) The $(\lambda,\mu)$-linear-interval order dimension of $\mathfrak{P}$
is the least cardinal $\lambda$ such that there is a linear-interval embedding 
\par
\!\!\! \!\!
of $(X,\overline{R})$ into a strict direct product of $\lambda$ linear-interval orders
which $\mu$ of them are not strict linear orders.
\par\noindent
(c) The $(\lambda,\mu)$-linear-interval order dimension of $\mathfrak{P}$
is the least cardinal $\lambda$ such that there is a strong linear-interval 
\par
\!\!\! \!\!
embedding of $(X,\overline{R})$ into a direct product of $\lambda$ strong linear-interval orders which $\mu$ of them are not linear
orders.
\end{theorem}

The following corollary is an immediate consequence of Theorem \ref{aw128}.

\begin{corollary}\label{wxa128} Let $\mathfrak{G}=(X,\prec)$ be a poset. Then the following statements are equivalent.
\par\noindent
(a) The $(\lambda,\mu)$-linear-interval order dimension of $\mathfrak{G}$ is the least cardinal $\lambda$ such that $\prec$ is the intersection of $\lambda$ linear-
\par
\!\!\! \!\!
interval orders which $\mu$ of them are not linear orders.
\par\noindent
(b) The $(\lambda,\mu)$-linear-interval order dimension of $\prec$
is the least cardinal $\lambda$ such that there is a linear-interval embedding 
\par
\!\!\! \!\!
of $(X,\prec)$ into a strict direct product of $\lambda$ linear-interval orders
which $\mu$ of them are not strict linear orders.
\par\noindent
(c) The $(\lambda,\mu)$-linear-interval order dimension of $\mathfrak{G}$
is the least cardinal $\lambda$ such that there is a strong linear-interval 
\par
\!\!\! \!\!
embedding of $(X,\prec)$ into a direct product of $\lambda$ strong linear-interval orders which $\mu$ of them are not linear
orders.
\end{corollary}

Using the previous approach for linear-interval orders,
we can define in a similar way the notion of 
(strong) linear-semiorder box embedding. The only difference is that
a semiorder is a poset whose elements correspond to unit length intervals.

\begin{definition}\label{au127}{\rm Let $\mathfrak{P} = (X, R)$ be an abstract system. (i)We call $sidpc(\mathfrak{P})$, the least cardinal $\lambda$ such that there is a linear-semiorder box embedding of $(X,\overline{R})$ into a direct product of $\lambda$ strong linear-semiorders.
(ii) 
We call $sispc(\mathfrak{P})$), the least cardinal $\lambda$ such that there is a linear-semiorder box embedding of $(X,\overline{R})$ into a direct product of $\lambda$ linear-semiorders.
}
\end{definition}

The following two theorems are proved in a similar way to the proof of Theorems
4.11 and 4.13, by using Theorem \ref{awerr}, Theorem \ref{awer1} and definition \ref{au127}.

\begin{theorem}\label{aw128} Let $\mathfrak{P}=(X,R)$ be an abstract system where $R$ is acyclic. Then the following statements are equivalent.
\par\noindent
(a) The $(\lambda,\mu)$-linear-semiorder dimension of $\mathfrak{P}$ is the least cardinal $\lambda$ such that $\overline{R}$ is the intersection of $\lambda$ linear-
\par
\!\!\! \!\!
semiorders which $\mu$ of them are not linear orders.
\par\noindent
(b) The $(\lambda,\mu)$-linear-semiorder dimension of $\mathfrak{P}$
is the least cardinal $\lambda$ such that there is a linear-semiorder embedding 
\par
\!\!\! \!\!
of $(X,\overline{R})$ into a strict direct product of $\lambda$ linear-semiorders
which $\mu$ of them are not strict linear orders.
\par\noindent
(c) The $(\lambda,\mu)$-linear-semiorder dimension of $\mathfrak{P}$
is the least cardinal $\lambda$ such that there is a strong linear-semiorder 
\par
\!\!\! \!\!
embedding of $(X,\overline{R})$ into a direct product of $\lambda$ strong linear-semiorders which $\mu$ of them are not linear
orders.
\end{theorem}

The following corollary is an immediate consequence of Theorem \ref{aw128}.

\begin{corollary}\label{awxa128} Let $\mathfrak{G}=(X,\prec)$ be a poset. Then the following statements are equivalent.
\par\noindent
(a) The $(\lambda,\mu)$-linear-semiorder dimension of $\mathfrak{G}$ is the least cardinal $\lambda$ such that $\prec$ is the intersection of $\lambda$ linear-
\par
\!\!\! \!\!
semiorders which $\mu$ of them are not linear orders.
\par\noindent
(b) The $(\lambda,\mu)$-linear-semiorder dimension of $\prec$
is the least cardinal $\lambda$ such that there is a linear-semiorder embedding 
\par
\!\!\! \!\!
of $(X,\prec)$ into a strict direct product of $\lambda$ linear-semiorders
which $\mu$ of them are not strict linear orders.
\par\noindent
(c) The $(\lambda,\mu)$-linear-interval order dimension of $\mathfrak{G}$
is the least cardinal $\lambda$ such that there is a strong linear-semiorder 
\par
\!\!\! \!\!
embedding of $(X,\prec)$ into a direct product of $\lambda$ strong linear-interval orders which $\mu$ of them are not linear
orders.
\end{corollary}

\par\bigskip\smallskip\par\noindent

\par\noindent
{\it Address}: {\tt {Athanasios Andrikopoulos} \\ {Department of Computer Engineering \& Informatics\\ University of Patras\\ Greece}
\par\noindent
{\it E-mail address}:{\tt aandriko@ceid.upatras.gr}

\end{document}

give an analogue of the Szpilran and dushnik-Miller theories in the interval order case.

analogue of the

Cerioli M., Oliveira F., Szwarcfiter S., (2008), Linear-Interval Dimension and PI Orders,
{\it Electronic Notes in Discrete Mathematics}, 30, 111-116.

\end{document}

Let $X$ be a non-empty decision space and
$\precsim\subseteq X\times X$ be a binary relation on $X$. As usual, $\prec$ denotes the {\it strict part} of $\precsim$.
We sometimes
abbreviate $(x,y)\in \precsim$ (resp. $(x,y)\in \prec$) as $x\precsim y$ (resp. $x\prec y$).
We say that $\precsim$ on $X$ is (i)
{\it reflexive} if for each $x\in X$, $x\precsim x$; (ii)
{\it asymmetric} if for all $x, y\in X$, $x\precsim y$
 and $y \not\precsim x$;
(iii) {\it transitive} if for all $x,y,z\in X$, [$x\precsim z$ 
and
$z\precsim y$ $\Rightarrow x\precsim y$; 
(iv)
{\it antisymmetric} if for each $x,y\in X$,
[$x\precsim y$ and
$y\precsim x] \Rightarrow x=y$; 
(vi) {\it complete} if
for all $x$ and $y$, $x\precsim y$ or $y\precsim x$; 
(v) {\it total} if for each $x,y\in X$,
$x\neq y$ we have $x\precsim y$ or $y\precsim x$.
The following combination of properties are considered in the next
theorems. A binary relation $\precsim$ on $X$ is:  (1) {\it A preorder} if
$\precsim$ is reflexive and transitive;
(2) {\it A partial order} if
$\precsim $ is reflexive, transitive and antisymmetric;
(3) {\it A linear order} if
$\precsim$ is a total partial order.
A {\it preordered set} is a pair $(X,\precsim)$ consisting of a set $X$ and a preorder $\precsim$ on $X$.
Let $\mathcal{L}=(X,\precsim)$ be a preordered set.
If $\precsim$ is an order, then $\mathcal{L}$ is called {\it partially ordered set} or {\it poset}.
A {\it lattice} is a mathematical structure studied in the mathematical subdisciplines of order theory and abstract algebra, among others.
It consists of a poset $(X,\precsim)$ in which every pair of elements $x$ and $y$
has a unique supremum $x\vee y$ (also called a least upper bound or join) and a unique 
infimum  $x\wedge y$ (also called a greatest lower bound or meet).
A subspace $(A,\precsim)$ of a lattice $(X,\precsim)$ is called a {\it sublattice} of $(X,\precsim)$ if it becomes a 
lattice with respect to $\precsim$.
If $\precsim$ is a preorder on $X$, then we denote the associated asymmetric relation $\prec$
and the associated equivalence relation $\sim$, respectively, by 
[$x\prec y\Leftrightarrow (x\precsim y)\wedge (y\not\precsim x)]$ and 
[$x\sim y\Leftrightarrow (x\precsim y)\wedge (y\precsim x)]$. 
We recall that $f: (X,\precsim)\to (\mathbb{R},\leq)$
is: (i) {\it increasing} if, for each $x,y\in X$, $x\precsim y$ implies $f(x)\leq f(y)$;
(ii) {\it order preserving} if for all $x ,y\in X$,
 $x\prec y\Rightarrow f(x)<f(y)$].
A {\it Richter-Peleg utility} is an increasing function $f: (X,\precsim)\to (\mathbb{R},\leq)$ 
that is also an order-preserving function
(see e.g. Peleg \cite{pel} and Richter \cite{ric0}). 
Equivalently, if $\sim$ implies $f(x)=f(y)$ and it is order preserving.
We say that a preorder $\precsim$ admits
a {\it Richter-Peleg multi-utility representation} by a family of functions $\mathcal{V}$ if
$x\precsim y \Leftrightarrow  f(x)\leq f(y)$ for all $f\in\mathcal{V}$. 
It is obvious that a Richter-Peleg multi-utility representation $\mathcal{V}$ of a preordered set
$(X,\precsim)$
characterizes the strict part $\prec$ of $\precsim$, in the sense that for each $x,y \in X$,
$x\prec y \Leftrightarrow  f(x)<f(y)$ for all $f\in\mathcal{V}$.

\par
Let $(X,\tau)$ be a topological space. 
We say that $(X,\tau)$
is a: ($\mathfrak{i}$) $T_{_0}$ topological space if given two distinct points in it, there exists an open neighborbood of it that contains exactly one of them;
($\mathfrak{ii}$) $T_{_2}$ or {\it Hausdorff} topological space
if its distinct points are contained in disjoint neighborhoods and ($\mathfrak{iii}$)
{\it quasi-compact} if for each collection of open sets which covers $X$ there exists a finite subcollection that also covers $X$.
If it is quasi-compact and Hausdorff, it is called {\it compact}. 
A function $f$ in $(X,\tau)$ is {\it upper} (resp., {\it lower})  
 {\it semicontinuous}
at $x\in X$ if for each $\varepsilon>0$,
there exists a neighborhood $U_{_x}$ of $x$ 
such that for all $y\in U_{_x}$ we have $f(y)<f(x)+\varepsilon$ (resp. $f(y)>f(x)-\varepsilon$).
An {\it isomorphism} is a structure-preserving mapping $f$ that can be reversed by $f^{-1}$ between two structures of the same type.
If there is an isomorphism between two mathematical structures, they are said to be {\it isomorphic}.
If $A\subseteq X$, then is the {\it interior} (resp. {\it closure}) of $A$ with respect to $\tau$ topology is denoted by 
$int_{_{\tau}}A$ (resp. $cl_{_{\tau}}A$).
Every topology $\tau$ on a set $X$ induces a preorder $\precsim_{_\tau}$ on this set, called {\it specialization preorder}, as follows:
$x\precsim_{_\tau} y$ if and only if $x\in cl_{_\tau}y$, $x, y\in X$.
This preorder is an order if and only if $(X,\tau)$ is $T_{_0}$. The {\it compatible topologies}
on a preordered set are those which induce the given preorder.
Nachbin \cite{nac} defines a {\it topological preordered space} $(X,\tau,\precsim)$ as a topological space $(X,\tau)$ with a preorder $\precsim$ that is closed as a subset of $X\times X$.
If $(X,\tau)$ is Hausdorff, then $\precsim$ is a partial order and $(X,\tau,\precsim)$ 
is called a {\it topological ordered space}.
For any $A\subseteq X$, define $\displaystyle\uparrow A=\{y\in X\vert x\precsim y\ {\rm for\ some}\ x\in A\}$.
We also write $\displaystyle\uparrow \{x\}$ as $\displaystyle\uparrow x$.
The sets $\displaystyle\downarrow A$ as $\displaystyle\downarrow x$
are defined dually.
A subset $A$ of $X$ is said to be an increasing (resp. decreasing) set, or to be an upper (lower) set,
if
$A=\displaystyle\uparrow A$ (resp. $A=\displaystyle\downarrow A$).
In what follows, we use two names for the sets in the preceding definition: "increasing (resp. decreasing) set" as in general topology, or "upper (lower) set" as in order theory.
A subset $A$ of $X$ is increasing if and only if $X\setminus A$ is decreasing. For each subset $A$ of $X$ 
there is a smallest increasing set $i(A)$ (decreasing set $d(A)$) that contains $A$.
If $A=\{a\}$ for some $a\in X$, then the notation is $i(a)$ (resp. $d(a)$).
A set $X$ equipped with two topologies $\tau_{_1}$ and $\tau_{_2}$
is called a {\it bitopological space}.
The concept of bitopological spaces was first introduced by Kelly \cite{kel}
in order to generalize the notion of topological space.
Every bitopological space $\mathfrak{D}=(X,\tau_{_1},\tau_{_2})$
can be regarded as a topological space 
$(X,\tau)$ with $\tau_{_1}=\tau_{_2}=\tau$. 
The dual of $\mathfrak{D}$ is the bitopological space
$\mathfrak{D}^{\ast}=(X,\tau_{_2},\tau_{_1})$.
A bitopological space $\mathfrak{D}=(X,\tau_{_1},\tau_{_2})$ is called:
(i) {\it  Pairwise Hausdorff} if for each two points $x$ and $y$ in $X$, there exists$\tau_{_i}$-open neighborhood of $x$
and $\tau_{_j}$-open neighborhood $V_j$ of $y$
such that $U_i\cap V_j=\emptyset$, \  $i, j\in \{1,2\}, i\neq j$.
(ii) {\it Pairwise normal} if for every pair of a $\tau_{_1}$-closed set $F$ and a $\tau_{_2}$-closed set $G$ with 
$F\cap G=\emptyset$, there exist a $\tau_{_1}$-open set $U$ and a $\tau_{_2}$-open set $V$
such that $F\subset V$, $G\subset U$ and $V\cap U=\emptyset$.
A bitopological ordered space $\mathfrak{X}=(X,\tau_{_1},\tau_{_2},\precsim)$ is called
{\it Pairwise normally ordered} if for every pair of a decreasing $\tau_{_1}$-closed set $F$ and an increasing $\tau_{_2}$-closed set $G$ with 
$F\cap G=\emptyset$, there exist an increasing $\tau_{_1}$-open set $U$ and a decreasing $\tau_{_2}$-open set $V$
such that $F\subset V$, $G\subset U$ and $V\cap U=\emptyset$.

\section{Bitopological preordered spaces}

A {\it bitopological preordered space} 
$\mathfrak{X}=(X,\tau_{_1},\tau_{_2},\precsim)$ (see \cite[Definition 1.4]{and}) is a bitopological space 
$\mathfrak{D}=(X,\tau_{_1},\tau_{_2})$
equipped with a $\tau_{_1}\times \tau_{_2}$-closed preorder $\precsim$ in $X\times X$. 
The definition of $\mathfrak{X}$ allows for the definition of the dual bitopological preordered space $\mathfrak{X}^{\ast}=(X,\tau_{_2},\tau_{_1},\precsim^{\ast})$
of $\mathfrak{X}$, where $\precsim^{\ast}=(\precsim)^{-1}$.
If $\mathfrak{X}$ is pairwise Hausdorff, then
$(X,\tau_{_1}\vee\tau_{_2})$
is Hausdorff and $\precsim$ is a partial order ($\precsim$ is antisymmetric). In this case, $\mathfrak{X}$ is called a {\it bitopological ordered space}.
This definition extends the notion of topological ordered space\footnote{Acording to Nachbin, a {\it topological preordered space} $(X,\tau,\precsim)$ is a topological space $(X,\tau)$ equipped with a 
preorder $\precsim$ which is $\tau\times \tau$-closed subset of $X\times X$.
If $(X,\tau)$ is Hausdorff, then $\precsim$ is a partial order and $(X,\tau,\precsim)$ is called a {\it topological ordered space}.}
of Nachbin to the bitopological spaces.
Every bitopological ordered space $\mathfrak{X}$ can be thought of as a
bitopological space $\mathfrak{D}=(X,\tau_{_1},\tau_{_2})$, where $\precsim$ is the equality relation 
$\Delta=\{(x,x)\vert x\in X\}$.
Each member of $\tau_{_1}$ (resp. $\tau_{_2}$) in 
$\mathfrak{X}$
is called an {\it open}
set, or a $\tau_{_1}$-{\it open} (resp. $\tau_{_2}$-{\it open}) set if one wishes to emphasize the topologies $\tau_{_1}$ (resp. $\tau_{_2}$) on $X$. 
This notion comes in handy when addressing many topologies on a given $X$ at the same time.

\begin{definition}\label{osr}{\rm Let $\mathfrak{X}=(X,\tau_{_1},\tau_{_2},\precsim)$ 
be a bitopological preordered space. A subset $N$ of $X$ is a 
$\tau_{_1}$-neighborhood (resp. $\tau_{_2}$-neighborhood) of an $x\in X$
if and only if $N$ includes a $\tau_{_1}$-open (resp. $\tau_{_2}$-open) set
containing $x$. If $N$ is a 
$\tau_{_1}$-open (resp. $\tau_{_2}$-open) set then it is called
$\tau_{_1}$-open neighborhood (resp. $\tau_{_2}$-open neighborhood) of $x$.
}
\end{definition}

\begin{example}\label{ex0}{\rm Let $\mathfrak{I}=(\mathbb{R},\mathfrak{U},\mathfrak{L},\precsim)$,
where $\mathfrak{U}$, $\mathfrak{L}$ are the upper and lower topologies:
$\{(a,\infty)\vert a\in \mathbb{R}\}$, $\{ (-\infty,a)\vert \ a\in \mathbb{R}\}$ and 
$\precsim$ is the usual order in $\mathbb{R}$. Then it is easy to check that $\mathfrak{I}$ 
is a bitopological preordered space.
}
\end{example}

In what follows,
$\mathfrak{I}$ will denote the bitopological preordered space 
$(\mathbb{R},\mathfrak{U},\mathfrak{L},\precsim)$.

\begin{example}\label{ex1}{\rm Let $(X,\tau)$ be a topological space, $\tau^{\ast}$ be the discrete topology on $X$ and $\precsim_{_{\tau}}$ be the specialization order of $\tau$. Then,
$(X,\tau,\tau^{\ast},\precsim_{_{\tau}})$ is a bitopological preordered space, that is,  $\precsim_{_{\tau}}$ is 
$\tau\times \tau^{\ast}$-closed subset of $X\times X$. Indeed,
if 
$(a,b)\in X\times X\setminus \{(x,y)\in X\times X\vert \ x\precsim_{_{\tau}}y\}$, then there exists a $\tau$-open
neighborhood $V_{a}$ of $a$ such that 
$V_{a}\cap \{b\}=\emptyset$ or equivalently
$b\notin V_{_a}$. Then,
$V_{_a}\times \{b\}\subset X\times X\setminus \{(x,y)\in X\times X\vert \ x\precsim_{_{\tau}}y\}$.
Indeed, suppose to the contrary that there exists $z\in X$ such that
$(z,b)\in V_{_a}\times \{b\}$ and $(z,b)\notin X\times X\setminus \{(x,y)\in X\times X\vert \ x\precsim_{_{\tau}}y\}$. Then, $z\precsim_{_{\tau}}b$ which implies that $b\in \bigcap V_{z}$ ($V_{z}$ denotes
an arbitrary $\tau$-open neighborhood of $z$). Since $z\in V_a$ and $V_a$ is 
$\tau$-open, there exists a $\tau$-open neighborhood $V^{\ast}_{z}$ of $z$ such that
$V^{\ast}_{_z}\subset V_a$. But then, $b\in \bigcap V_{z}\subseteq V^{\ast}_{_z}\subset V_a$, a contradiction.
In light of this last contradiction, it can be concluded that
$\precsim_{_{\tau}}$ is $\tau\times \tau^{\ast}$-closed.
}
\end{example}

\begin{example}\label{ex2}{\rm Let $(X,d)$ be a quasi-pseudometric space\footnote{A quasi-pseudometric space $(X,d)$ is a set $X$ together with a non-negative real-valued function $d: X\times X\longrightarrow \mathbb{R}$ (called a quasi-pseudometric) such that, for every $x, y, z\in X$: (i) $d(x,x)=0$;
(ii) $d(x,y)\leq d(x,z)+d(z,y)$. A quasi-pseudometric $d$ on $X$ induces a topology $\tau_{d}$ on $X$ 
which has as a base the family of $d$-balls $\{ B_{_d}(x,r): x\in X, r>0\}$ where $B_{_d}(x,r)=\{y\in X: d(x,y)<r\}$.}.
Then, $(X,\tau_{_d},\tau_{_{d^{-1}}},\precsim_{_{\tau_{_d}}})$ is a bitopological preordered space.
To demonstrate this, we must show that $cl_{_{\tau_{_d}\times \tau_{_{d^{-1}}}}}G_{_{\tau_{_d}}}=G_{_{\tau_{_d}}}$,
where $G_{_{\tau_{_d}}}$ is the graph of $\precsim_{_{\tau_{_d}}}$ and $cl_{_{\tau_{_d}\times \tau_{_{d^{-1}}}}}G_{_{\tau_{_d}}}$
is the closure of $G_{_{\tau_{_d}}}$ in $(X\times X, \tau_{_d}\times \tau_{_{d^{-1}}})$. 
Indeed, let $(a,b)\in cl_{_{\tau_{_d}\times \tau_{_{d^{-1}}}}}G_{_{\tau_{_d}}}$. Then, for each $\varepsilon>0$, we have 
$B_{_d}(a,\frac{\varepsilon}{3})\times B_{_{d^{-1}}}(b,\frac{\varepsilon}{3})\cap G_{_{\tau_{_d}}}\neq \emptyset$. Therefore, there exists 
$(z_{_1},z_{_2})\in B_{_d}(a,\frac{\varepsilon}{3})\times B_{_{d^{-1}}}(b,\frac{\varepsilon}{3})$ and 
$(z_{_1},z_{_2})\in G_{_{\tau_{_d}}}$ or equivalently
$d(a,z_{_1})<\frac{\varepsilon}{3}$, $d(z_{_2},b)<\frac{\varepsilon}{3}$ and $d(z_{_1},z_{_2})=0<\frac{\varepsilon}{3}$. We have
$d(a,b)\leq d(a,z_{_1})+d(z_{_1},z_{_2})+d(z_{_2},b)<\frac{\varepsilon}{3}+\frac{\varepsilon}{3}+\frac{\varepsilon}{3}=\varepsilon$.
As a result, we can conclude that $d(a,b)<\varepsilon$ is true for every $\varepsilon>0$.
It follows that
$d(a,b)=0$ which implies that $(a,b)\in G_{_{\tau_{_d}}}$.

}
\end{example}

\begin{theorem}\label{a3}  {\rm Let $\mathfrak{D}=(X,\tau_{_1},\tau_{_2})$ be a bitopological space
and $\precsim$ a preorder on $X$. Then,
\par
($\mathfrak{a}$) The preorder is $\tau_{_1}\times \tau_{_2}$-closed if and only if for every two points $a, b\in X$ such that 
$a\not\precsim b$ there exist a 
increasing $\tau_{_1}$-neighborhood $V_{_a}$ of $a$ and a 
decreasing $\tau_{_2}$-neighborhood $V_{_b}$ of $b$ which are disjoint.
\par
($\mathfrak{b}$) If the preorder is $\tau_{_1}\times \tau_{_2}$-closed, then
for every $a\in X$ the set $d(a)$ is $\tau_{_1}$-closed and 
the set $i(a)$ is $\tau_{_2}$-closed.
\par
($\mathfrak{c}$)} If the preorder is $\tau_{_1}\times \tau_{_2}$-closed, then
$ a\in cl_{_{\tau_{_1}}}\{b\}\Rightarrow
a\precsim b$ and $ b\in cl_{_{\tau_{_2}}}\{a\}\Rightarrow
a\precsim b$.
\end{theorem}
\begin{proof}
$\mathfrak{a})$ Let $\precsim$ is a $\tau_{_1}\times \tau_{_2}$-closed
and
$a, b\in X$ such that $a\not\precsim b$. 
Since $(a,b)$ does not belong to the graph $G_{_{\precsim}}$ of $\precsim$ and $G_{_{\precsim}}$ is $\tau_{_1}\times \tau_{_2}$-closed, there exist a  $\tau_{_1}$-open neighborhood $V_a$ of $a$ and
a $\tau_{_2}$-open neighborhood $V_b$ of
$b$ such that
$V_a\times V_b\subset X\times X\setminus G_{_{\precsim}}$.
Let $\widetilde{V}_a=\{\lambda\in X\vert \exists \mu\in V_a\ \ {\rm such\ that}\ \ \mu\precsim \lambda\}$ and $\widetilde{V}_b=\{\kappa\in X\vert \exists \nu\in V_b\ \ {\rm such\ that}\ \ \kappa\precsim \nu\}$.
Clearly, $\widetilde{V}_a$ is an increasing $\tau_{_1}$-neighborhood of $a$
and
$\widetilde{V}_b$ is an decreasing $\tau_{_2}$-neighborhood of $b$ 
such that
$ \widetilde{V}_a\cap \widetilde{V}_b=\emptyset$. Indeed,
suppose to the contrary that it is not the case and let
$z\in   \widetilde{V}_a\cap \widetilde{V}_b$. 
Then, $\mu\precsim z$ and $z\precsim \nu$ implies that $\mu\precsim \nu$. Hence,
$(\mu,\nu)\in G_{_{\precsim}}$ and 
$(\mu,\nu)\in  V_{_a}\times V_{_{b}}\subset X\times X\setminus G_{_{\precsim}}$ which is impossible.

Conversely, suppose to the contrary that 
for each $a, b\in X$, $a\not\precsim b$ implies that
there exists 
an increasing $\tau_{_1}$-neighborhood $V_a$ of $a$
and a decreasing $\tau_{_2}$-neighborhood $V_b$ of $b$ 
such that
$ V_a\cap V_b=\emptyset$,
while $\precsim$ is not $\tau_{_1}\times \tau_{_2}$-closed.
Thus, there exists
$(a,b)\in X\times X$ such that $(a,b)\in X\times X\setminus G_{_{\precsim}}$ and
$(a,b)\in cl_{_{\tau_{_1}\times\tau_{_2}}}G_{_{\precsim}}$, or equivalently,
$a\not\precsim b$ and for each
$\tau_{_1}$-neighborhood $V_a$ of $a$
and each $\tau_{_2}$-neighborhood $V_b$ of $b$ there holds $V_a\times V_b\cap G_{_{\precsim}}
\neq\emptyset$. If 
$V_a$ is increasing and $V_b$ is decreasing, then
from $V_a\times V_b\cap G_{_{\precsim}}\neq \emptyset$, we conclude that there exist $x, y\in X$ such that 
$(x,y)\in V_a\times V_b$ and $x\precsim y$.
It follows that $x, y\in V_a\cap V_b$, a contradiction
to our assumption that $V_a\cap V_b$ must be empty for at least one pair of these neighborhoods. 
Therefore, $\precsim$ is $\tau_{_1}\times \tau_{_2}$-closed.
\par\noindent
$\mathfrak{b})$ Let $a\in X$. To show that $d(a)$ is $\tau_{_1}$-closed we prove that
$X\setminus d(a)$ is $\tau_{_1}$-open. Indeed, let $b\in X\setminus d(a)$.
It follows that $b\not\precsim a$. By the first part of proposition
there exists 
an increasing $\tau_{_1}$-neighborhood $\widetilde{V}_b$ of $b$
and a decreasing $\tau_{_2}$-neighborhood $\widetilde{V}_a$ of $a$ 
such that
$ \widetilde{V}_b\cap \widetilde{V}_a=\emptyset$. Hence, there exists a $\tau_{_1}$-open neighborhood $V_b$ of $b$ such that $b\in V_b\subseteq \widetilde{V}_b$ and 
$d(a)\subseteq \widetilde{V}_a$. It follows that $d(a)\cap V_b=\emptyset$ which implies
that $X\setminus d(a)$ is $\tau_{_1}$-open ($b\in V_b\subset X\setminus d(a)$). 
\par\noindent
$\mathfrak{c})$ Let
$a\in cl_{_{\tau_{_1}}}\{b\}$. Then, $b$ belong to all the $\tau_{_1}$-open neighborhoods
$V_a$ of $a$. Suppose to the contrary that $a\not\precsim b$. Then, by the part 
($\mathfrak{a})$
of proposition,
there exist a $\tau_{_1}$-neighborhood $ \widetilde{V}_a$ of $a$
and a
$\tau_{_2}$-neighborhood $ \widetilde{V}_b$ of $b$ 
such that
$ \widetilde{V}_a\cap \widetilde{V}_b=\emptyset$.
By definition, there exists a $\tau_{_1}$-open neighborhood $V_a$ of $a$ such that
$V_a\subseteq \widetilde{V}_a$.
It follows that $b\notin V_a$, a contradiction. Therefore, $a\precsim b$. Similarly we prove that 
$ b\in cl_{_{\tau_{_2}}}\{a\}\Rightarrow
a\precsim b$.
\end{proof}

\begin{proposition}\label{pan}{\rm Each bitopological space $\mathfrak{D}=(X,\tau_{_1},\tau_{_2})$
equipped with a closed preorder is 
a pairwise Hausdorff space; that is $\tau_1\bigvee \tau_2$ is a Hausdoff topology. }
\end{proposition}
\begin{proof}
Take two distinct points into consideration $x, y\in X$. Due to the fact that we are dealing with a preorder $\precsim$, one of the two relationships $x\precsim y$,  $y\precsim x$
is false. Assume the first is false (the case of the second is analogous).
By Theorem \ref{a3},
there exist a 
increasing $\tau_{_1}$-neighborhood $V_{_x}$ of $x$ and a 
decreasing $\tau_{_2}$-neighborhood $V_{_y}$ of $y$ which are disjoint. Since the dual space
$\mathfrak{X}^{\ast}=(X,\tau_{_2},\tau_{_1},(\precsim)^{-1})$,  is also a bitopological preordered space, 
there exist a 
increasing $\tau_{_2}$-neighborhood $V_{_y}$ of $y$ and a 
decreasing $\tau_{_1}$-neighborhood $V_{_x}$ of $x$ which are disjoint. Therefore, 
$\mathfrak{D}$ is a pairwise Hausdorff space, $\tau_1\bigvee \tau_2$ is a Hausdoff topology.
\end{proof}

In \cite{and} Andrikopoulos presents a version of Nachbin's theory based on bitopological ordered spaces. The following theorem is taken from Andrikopoulos \cite[Theorem 1.8]{and}.

\begin{theorem}\label{a231}{\rm A bitopological ordered space $\mathfrak{X}=(X,\tau_{_1},\tau_{_2},\precsim)$
is pairwise normally ordered if and only if given a decreasing $\tau_{_1}$-closed
set $A$ and an increasing $\tau_{_2}$-closed
set $B$ with $A\cap B=\emptyset$, there exists an increasing real-valued function $f$ on $X$ such that
\par
(i) $f(A)=0$, $f(B)=1$ and $0\leq f(x)\leq 1$,
\par
(ii) $f$ is 
$\tau_{_1}$-lower semicontinuous and 
$\tau_{_2}$-upper semicontinuous.}
\end{theorem}

\begin{proof} To prove sufficiency, we follow the line of the proof of the classical theorems of Nachbin \cite[Theorem 2]{nac} and Kelly \cite[Theorem 2.7]{kel} (see also \cite[Theorem 1.8]{and}. 
Let $A$ be
a decreasing $\tau_{_1}$-closed
set in $X$ and $B$ be an increasing $\tau_{_2}$-closed
set $B$ in $X$ with $A\cap B=\emptyset$. 
Since $\mathfrak{X}$ is pairwise normal there exist an 
increasing 
$\tau_{_1}$-open set $O_{_1}$ and a
decreasing $\tau_{_2}$-open set $O_{_2}$ such that $A\subset O_{_2}$, $B\subset O_{_1}$and $O_{_1}\cap O_{_2}=\emptyset$. We put 
$A_{_0}=A$, $A_{_{1\over 2}}=X\setminus O_{_1}$, $B_{_{1\over 2}}=O_{_2}$ and 
$B_{_1}=X\setminus B$.
Then, we have 
$A_{_0}\subseteq B_{_{1\over 2}}\subseteq A_{_{1\over 2}}\subseteq B_{_1}$, $cl_{_{\tau_{_1}}}B_{_{1\over 2}}
\subseteq B_{_1}$. When we apply our hypothesis on $\mathfrak{X}$ to each pair of sets $A_{_0}=A, B_{_{1\over 2}}$
and $A_{_{1\over 2}}, B_{_1}$,
we get decreasing $\tau_{_1}$-closed sets $A_{_{1\over 4}}, A_{_{3\over 4}}$ 
and increasing $\tau_{_2}$-oped sets $B_{_{1\over 4}}, B_{_{3\over 4}}$ 
such that
\begin{center}
$A_{_0}\subseteq B_{_{1\over 4}}\subseteq A_{_{1\over 4}}\subseteq B_{_{1\over 2}}\subseteq A_{_{1\over 2}}\subseteq 
B_{_{3\over 4}}\subseteq A_{_{3\over 4}}\subseteq B_{_1}$.
\end{center}
Continuing this process, keeping in mind that dyadic rationals are dense in R,
we obtain by induction families
$(A_{_s})_{_{s\in S}}$, $(B_{_s})_{_{s\in S}}$, where 
$S=\{\frac{p}{2^{^q}}\vert \ p=1,2,..., 2^q-1,\ q\in \mathbb{N}\setminus \{0\}\}$. 
We put $A_{_s}=\emptyset$ if $s<0$, $A_{_s}=X$ if $s\geq 1$ and 
$B_{_s}=\emptyset$ if $s\leq 0$, $B_{_s}=X$ if $s>1$.
Then,
\begin{center}
$B_{_r}\subseteq B_{_s}\subseteq A_{_s}\subseteq B_{_t}$ if $r\leq s\leq t$, and $A_{_s}\subseteq B_{_t}$ if $s<t$.
\end{center}
Define 
\begin{center}
$f(x)=inf\{t\in S\vert x\in B_{_t}\}$.
\end{center}
Then,
\begin{center}
$f(x)=inf\{t\in S\vert x\in A_{_t}\}$.
\end{center}

We have $f(A_{_0})=f(A)=0$, $f(B_{_1})=f(X\setminus B)=1$ and $0\leq f(x)\leq 1$ for all
$x\in X$.
To prove that $f$ is increasing, let $x\precsim y$ for some $x, y\in X$. 
As a result, for some $t^{\ast}$ in $S$, $y\in B_{_{t^{\ast}}}$ holds true.
 Since $B_{_{t^{\ast}}}$ is decreasing we have that 
$x\in B_{_{t^{\ast}}}$.
By the definition of $f$ we have $f(x)=inf\{t\in S\vert x\in B_{_{t}}\}\leq t^{\ast}$.
Due to the arbitrariness of $t^{\ast}$, we conclude that
$f(x)\leq inf\{t^{\ast}\in S\vert y\in B_{_{t^{\ast}}}\}=f(y)$.

It remains to prove that $f$ is a 
$\tau_{_1}$-lower semicontinuous 
$\tau_{_2}$-upper semicontinuous function on $\mathfrak{X}$.
By the above construction, we have
\begin{center}
$A_{_t}\subseteq int_{_{\tau_{_2}}}A_{_s}$ and $cl_{_{\tau_{_1}}}B_{_t}\subseteq B_{_s}$ for $t<s$. 
\end{center}
To prove that $f$ is a 
$\tau_{_1}$-lower semicontinuous, retaining the previous proof process's symbolism, 
let $t, r, s\in S$ and $\varepsilon>0$ such that $f(x)-\varepsilon<t<r<s<f(x)$.
We know that $x\notin B_{_s}$, so $x\in X\setminus cl_{_{\tau_{_1}}}B_{_r}=G$.
For each $t<r$ and $y\in G$ we have that $y\notin B_{_t}$ and thus $f(x)-\varepsilon<r<f(y)$.
It follows that  $f$ is a 
$\tau_{_1}$-lower semicontinuous function.
To prove that $f$ is a 
$\tau_{_2}$-upper semicontinuous,
let $t, s\in S$ and $\varepsilon>0$ such that $f(x)<t<s<f(x)+\varepsilon$.
We have that $x\in A_{_t}\subseteq int_{_{\tau_{_2}}}A_{_s}$. If $y \in int_{_{\tau_{_2}}}A_{_s}$,
then $f(y)<s<f(x)+\varepsilon$ which completes the proof.

To prove necessity, let
$A, B$ be two disjoint subsets of $X$ such that $A$ is $\tau_{_1}$-closed and $B$ is $\tau_{_2}$-closed.  By hypothesis, there exists an increasing  
$\tau_{_1}$-lower semicontinuous and $\tau_{_2}$-upper semicontinuous function $f$ such that
$f(A)=0, f(B)=1$ and $0\leq f(x)\leq 1$ for each $x\in X$. By the $\tau_{_1}$-lower semicontinuouity of $f$, the set 
$U=f^{-1}((\frac{1}{2},1])$ is $\tau_{_1}$-open such that $B\subset U$ and 
by the $\tau_{_2}$-upper semicontinuouity of $f$, the set 
$V=f^{-1}([0,\frac{1}{2}))$ is $\tau_{_2}$-open set
such that $A\subset V$. 
Clearly, $U\cap V=\emptyset$. On the other hand,
since $f$ is increasing we have that $V$ is decreasing and $U$ is increasing.
\end{proof}

In case where the order considered is
the equality relation 
$\Delta=\{(x,x)\vert x\in X\}$, Theorem \ref{a231} reduces to Kelly's theorem \cite[Theorem 2.7]{kel}.

\begin{proposition}\label{a4}  {\rm Let $\mathfrak{X}=(X,\tau_{_1},\tau_{_2},\precsim)$ be a bitopological 
ordered space. 
\par
($\mathfrak{a}$) If $A$ is a $\tau_{_1}$-compact subset of $X$ and  
$B$ is a $\tau_2$-compact subset of $X$, then the increasing subset $i(A)$ generated 
by $A$ is $\tau_{_2}$-closed
and the decreasing subset $d(B)$ generated by $B$ is $\tau_{_1}$-closed.
\par
($\mathfrak{b}$) If $\tau_{_1}\vee \tau_{_2}$ is quasicompact, then each $\tau_{_1}$-closed subset of $X$ is $\tau_{_2}$-quasicompact and each $\tau_{_2}$-closed subset of $X$ is $\tau_{_1}$-quasicompact.
}
\end{proposition}
\begin{proof}($\mathfrak{a}$) To show that $i(A)$ is $\tau_{_2}$-closed, it is enough to prove that 
$X\setminus i(A)$ is $\tau_{_2}$-open. Indeed, let $b\in  X\setminus i(A)$. Then, 
for each $a\in A$, we have $a\not\precsim b$.
By Theorem \ref{a3}($\mathfrak{a}$), there is an increasing $\tau_{_1}$-neighborhood 
$\widetilde{V}_{_a}$ of $a$ and a decreasing $\tau_{_2}$-neighborhood 
$\widetilde{V}_{_{b(a)}}$ of $b$ such that
$ \widetilde{V}_{_a}\cap \widetilde{V}_{_{b(a)}}=\emptyset$.
On the other hand, for each $a\in A$ there exists an $\tau_{_1}$-open neighborhood $V_a$ of
$a$ such that $V_a\subseteq \widetilde{V}_{_a}$. 
Let $\mathcal{C}=\{V_{_a}\vert a\in A\}$. Then, $\mathcal{C}$ is a $\tau_1$-open cover of $A$.
Therefore, by $\tau_{_1}$-compactness of $A$, there exist $a_{_1}, a_{_2},..., a_{_n}\in A$ such that 
$A\subset V_{_{a_{_1}}}\cup V_{_{a_{_2}}}...
\cup V_{_{a_{_n}}}
\subset \widetilde{V}_{_{a_{_1}}}\cup \widetilde{V}_{_{a_{_2}}}...
\cup \widetilde{V}_{_{a_{_n}}}$.
Let
$ \widetilde{V}_{_{b}}= \widetilde{V}_{_{b(a_{_1})}}\cap \widetilde{V}_{_{b(a_{_2})}}...
\cap \widetilde{V}_{_{b(a_{_n})}}$. Then, $b\in \widetilde{V}_{_{b}}$
and 
\begin{center}
$\widetilde{V}_{_{b}}\cap A\subseteq \widetilde{V}_{_{b}}\cap
(\displaystyle\bigcup_{i\in \{1,2,...,n\}} \widetilde{V}_{_{a_{_i}}})=
\displaystyle\bigcup_{i\in \{1,2,...,n\}} \widetilde{V}_{_{b}}\cap\widetilde{V}_{_{a_{_i}}}
\subseteq \displaystyle\bigcup_{i\in \{1,2,...,n\}} \widetilde{V}_{_{b(a_{_i})}}\cap\widetilde{V}_{_{a_{_i}}}=\emptyset.
$
\end{center}
Therefore, $A\subset X\setminus \widetilde{V}_{_{b}}$. Since 
$X\setminus \widetilde{V}_{_{b}}$ is increasing, we conclude that 
$i(A)\subset X\setminus \widetilde{V}_{_{b}}$ or equivalently 
$i(A)\cap \widetilde{V}_{_{b}}=\emptyset$. Since $\widetilde{V}_{_{b}}$ is a 
$\tau_{_2}$-neighborhood of $b$, there exists a 
$\tau_{_2}$-open neighborhood $W_{_b}$ of $b$ such that $b\in W_{_b}\subseteq \widetilde{V}_{_{b}}$. It follows that 
$i(A)\cap W_{_{b}}=\emptyset$ which implies that  $i(A)$ is 
$\tau_{_2}$-closed. In a similar way we can prove that the subset
$d(A)$ generated by $A$ is $\tau_{_1}$-closed.
\par\smallskip\par\noindent
($\mathfrak{b}$) Let $A$ be a $\tau_{_1}$-closed subset of $X$, then it is $\tau_{_1}\vee \tau_{_2}$-closed.
Since $\tau_{_1}\vee \tau_{_2}$ is quasicompact, we have that 
 $A$ be a $\tau_{_1}\vee \tau_{_2}$-quasicompact subset of $X$. Then,  $A$ is quasicompact
in the weaker topology $\tau_{_2}$. 
\end{proof}

\section{Joincompact bitopological ordered spaces}

We can now proceed to the definition of joincompact bitopological spaces.

\begin{definition}\label{plk}{\rm A bitopological space $\mathfrak{D}=(X,\tau_{_1},\tau_{_2})$
is
{\it joincompact} if it is pairwise Hausdorff and the topology $\tau_{_1}\bigvee \tau_{_2}$ is quasi-compact.} 
\end{definition}

Joincompact spaces and
Lawson-closed subsets of continuous lattices are the same objects with different names. 
In addition, the order on a continuous lattice corresponds to the specialization order of the Scott 
topology of the lattice (see paragraph 6 below).

The following proposition is an immediate consequence of Definition \ref{plk} and Proposition \ref{pan}.

\begin{proposition} {\rm A quasi-compact bitopological preordered space $\mathfrak{X}=(X,\tau_{_1},\tau_{_2},$
$\precsim)$
is a joincompact bitopological ordered space.}
\end{proposition}

\begin{corollary}\label{a5}{\rm  Let $\mathfrak{X}=(X,\tau_{_1},\tau_{_2},\precsim)$ be a joincompact bitopological ordered space. Then,
the $\tau_{_1}$-closed sets are precisely the 
$\tau_{_2}$-compact sets and the $\tau_{_2}$-closed sets are precisely the 
$\tau_{_1}$-compact sets.
}
\end{corollary}

\begin{proposition}\label{a117} {Let $\mathfrak{X}=(X,\tau_1,\tau_2,\precsim)$ be a joincompact bitopological ordered space. 
If $A\subset X$ is a decreasing set and $O_{_1}$ is a $\tau_{_2}$-neighborhood 
of $A$, then there exists a $\tau_{_2}$-open decreasing neighborhood $O_{_2}$
of $A$ such that $A\subset O_{_2}\subset O_{_1}$.}
\end{proposition}
\begin{proof}Let $A$ and $O_{_1}$
be as in the supposition of the proposition.
Put 
$O_{_2}=X\setminus i(cl_{_{\tau_{_2}}}(X\setminus O_{_1}))$. Then, 
$O_{_2}$ is a decreasing $\tau_{_2}$-open subset of $X$. 
By Corollary \ref{a5},  $cl_{_{\tau_{_2}}}(X\setminus O_{_1})$ is $\tau_{_1}$-compact and
by Theorem \ref{a4}($\mathfrak{a}$) $ i(cl_{_{\tau_{_2}}}(X\setminus O_{_1}))$ is $\tau_{_2}$-closed. It follows that $O_{_2}$ is $\tau_{_2}$-open. Now, we have
\begin{center}
$X\setminus O_{_1}\subseteq cl_{_{\tau_{_2}}}(X\setminus O_{_1})\subseteq i(cl_{_{\tau_{_2}}}(X\setminus O_{_1}))$.
\end{center}
We show that $A\subset O_{_2}\subset O_{_1}$. Clearly, $O_{_2}\subset O_{_1}$. To prove that $A\subset O_{_2}$ it suffices to show that $A\cap i(cl_{_{\tau_{_2}}}(X\setminus O_{_1}))=\emptyset$. Indeed, suppose to the contrary that there exists $z\in X$ such that
$z\in A$ and $z\in i(cl_{_{\tau_{_2}}}(X\setminus O_{_1}))$. Then, there exists $w\in 
i(cl_{_{\tau_{_2}}}(X\setminus O_{_1}))$ such that $w\precsim z$. From $z\in A$ we get $w\in A$, a contradiction because $O_{_1}$ is a $\tau_{_2}$-neighborhood of $A$.
\end{proof}

The following proposition 
is the dual of Proposition \ref{a117} and the proof is similar to the one given in
Proposition \ref{a117}.

\begin{proposition}\label{a217} {Let $\mathfrak{X}=(X,\tau_1,\tau_2,\precsim)$ be a joincompact bitopological ordered space. 
If $B\subset X$ is an increasing set and $F_{_1}$ is a $\tau_{_1}$-neighborhood 
of $B$, then there exists a $\tau_{_1}$-open increasing neighborhood $F_{_2}$
of $B$ such that $B\subset F_{_2}\subset F_{_1}$.}
\end{proposition}

\begin{theorem}\label{a7}{\rm Every joincompact bitopological ordered space is
pairwise normally ordered space.}
\end{theorem}
\begin{proof} Consider a decreasing $\tau_{_1}$-closed set $A\subset X$ and an increasing $\tau_{_2}$-closed set $B\subset X$ such that $A\cap B=\emptyset$. Let $b\in B$.
Then, for each $x\in A$ we have $b\not\precsim x$. By Theorem \ref{a3}($\mathfrak{a}$),
there exists a decreasing $\tau_{_2}$-neighborhood $V_{_x}$ of $x$ and an increasing 
$\tau_{_1}$-neighborhood $V_{_b}$ of $b$ such that $V_{_x}\cap V_{_b}=\emptyset$.
By theorem \ref{a3}($\mathfrak{b}$), $d(x)$ is $\tau_{_1}$-closed and 
$i(b)$ is $\tau_{_2}$-closed. Since $\mathfrak{X}$ is joincompact, by \cite[Theorem 3.6(c)]{kop}
we have that
$\mathfrak{X}$ is pairwise normal. Thus, there exist
a $\tau_{_1}$-open set $W$ and 
a $\tau_{_2}$-open set $O$ such that $d(x)\subset O$, $i(b)\subset W$ and $O\cap W=\emptyset$. Using the fact that $d(x)$ is decreasing and applying Proposition \ref{a117},
we obtain a decreasing $\tau_{_2}$-open set $\widetilde{O}$ such that 
$d(x)\subset \widetilde{O}\subset O$. By a dual argument from Proposition \ref{a217},
we can get an increasing $\tau_{_1}$-open set $\widetilde{W}$ such that 
$i(b)\subset \widetilde{W}\subset W$. Since $b\not\precsim x$ for any $x\in A$,
we can find a decreasing $\tau_{_2}$-open set $O_{_x}$ containing $x$ and 
an increasing $\tau_{_1}$-open set $W_{_b}$ containing $b$ such that 
$O_{_x}\cap W_{_b}=\emptyset$. Since $\mathfrak{X}$ is joincompact and $A$ is
$\tau_{_1}$-closed, by Proposition \ref{a4}($\mathfrak{d}$), we have that $A$ is $\tau_{_2}$-quasicompact. Therefore, there exists a finite number of $\tau_{_2}$-open sets $O_{_{x_{_i}}}$, $i=\{1,2,...,n\}$ such that $x_{_i}\in A$ and $A\subseteq \displaystyle\bigcup_{_{i\in \{1,2,...,n\}}}O_{_{x_{_i}}}=\widehat{O}$.
Let $\widehat{W}=\{\displaystyle\bigcap_{_{i\in\{1,2,...,n\}}}W^{x_{_i}}_{_b}\vert O_{_{x_{_i}}}\cap W^{x_{_i}}_{_b}=\emptyset\}.$ It is then clear that $\widehat{W}$ is an increasing $\tau_{_1}$-open set containing $b$ and $\widehat{O}$ is a decreasing $\tau_{_2}$-open set containing $A$  such that
$\widehat{O}\cap \widehat{W}=\emptyset$.
Now, for each $y\in B$, this point $y$ and the set $A$ are in exactly
the same relation as $b$ and $x$ in the preceding procedure. Similarly, 
we find a $\tau_{_1}$-open set $O_{_1}$ and 
a $\tau_{_2}$-open set $O_{_2}$ such that $A\subset O_{_1}$, $B\subset O_{_2}$
and $O_{_1}\cap O_{_1}=\emptyset$.
\end{proof}

\begin{proposition}\label{a112} {Let $\mathfrak{X}=(X,\tau_1,\tau_2,\precsim)$ be a joincompact bitopological ordered space and let $a, b\in X$ such that
$a\not\precsim b$. Then, there exists an increasing $\tau_1$-lower and $\tau_2$-upper semicontinous function $f$ on $X$ such that $f(a)>f(b)$.}
\end{proposition}
\begin{proof} By Theorem \ref{a3}($\mathfrak{b}$) we have that $d(b)$ is a decreasing
$\tau_{_1}$-closed set and 
$i(a)$ is an increasing
$\tau_{_2}$-closed set. But then, Theorems \ref{a231} and \ref{a7} imply that
there exists a increasing
$\tau_{_1}$-lower semicontinuous and 
$\tau_{_2}$-upper semicontinuous function such that $f(b)=0<1=f(a)$.
\end{proof}

\begin{proposition}\label{a2}{\rm Let $\mathfrak{X}=(X,\tau_{_1},\tau_{_2},\precsim)$
be a bitopological ordered space and let $\mathfrak{I}=(\mathbb{R},\mathfrak{U},\mathfrak{L},\leq)$,
where $\mathfrak{U}$, $\mathfrak{L}$ are the upper and lower topologies in $\mathbb{R}$, and $\leq$ 
is the usual order in $\mathbb{R}$. 
If $\mathfrak{X}$ is joincompact and 
$f:\mathfrak{X}\to \mathfrak{I}$
is $\tau_{_1}$-lower and $\tau_{_2}$-upper semicontinuous
function, then there exist points at which the function takes its maximum and minimum value.
}
\end{proposition}
\begin{proof} Let's assume that $f$ is $\tau_{_1}$-lower and $\tau_{_2}$-upper semicontinuous
function on $X$.
Since $\mathfrak{X}$ is joincompact,
$f(X)$ is joincompact and hence it is compact in the weaker 
topologies $\mathfrak{U}$ and $\mathfrak{L}$. Therefore, 
$f(X)$ is $\mathfrak{U}$-closed and $\mathfrak{L}$-closed. Since the topologies $\mathfrak{U}$ and $\mathfrak{L}$
are compact $f$ is bounded, that is, $m\leq f(x)\leq M$ for all $x\in X$ and $m, M\in R$.
Let $\lambda$ be the infimum of set of values of $f(x)$ on $X$. Suppose to the contrary that $f(x)\neq \lambda$ for all $x\in X$.
Let $h(x)=\displaystyle\frac{1}{f(x)-\lambda}$ ($h(x)>0$). 
Then,  $h$ is 
$\tau_{_1}$-lower and $\tau_{_2}$-upper semicontinuous
function on $X$. As in the case of $f$, we conclude that $h$ is bounded. 
Because $\lambda$ is the infimum of $f(x)$ values, given $\epsilon>0$, we can find a value $f(x)<\lambda+\epsilon$.
Hence, $h(x)=\displaystyle\frac{1}{f(x)-\lambda}>\displaystyle\frac{1}{\epsilon}$.
Since $\epsilon$ is arbitrary the inequality 
$h(x)>\displaystyle\frac{1}{\epsilon}$ shows that $h(x)$ is unbounded in $X$, a contradiction.
Therefore,
$\lambda=f(x)$ for some $x\in X$.
Similarly, we conclude that if $\mu=\sup f(X)$ then $\mu=f(x)$ for some $x\in X$.
\end{proof}

Based on the previous proposition, 
in the following where the space $\mathfrak{X}=(X,\tau_{_1},\tau_{_2},\precsim)$ is joincompact, it will also be valid that the 
$\tau_{_1}$-lower and $\tau_{_2}$-upper semicontinuous
functions defined in it are bounded.

\section{Main result}

\begin{definition}\label{a222}{\rm
A function $u$ on $\mathfrak{X}=(X,\tau_{_1},\tau_{_2},\precsim)$ is said to be a
{\it Richter-Peleg utility representation function} for a preorder $\precsim$ on $\mathfrak{X}$ if 
$x \precsim y$ implies $u(x)\leq u(y)$ and 
$x \prec y$ implies $u(x)<u(y)$, where $\prec$ stands for the 
strict part of $\precsim$.
A {\it Richter-Peleg multi-utility representation} $\mathcal{V}$ for a preorder $\precsim$ on 
$\mathfrak{X}$ is a 
multi-utility representation for $\precsim$ such that every function $u\in V$ is a 
Richter-Peleg utility for $\precsim$.
A Richter-Peleg multi-utility representation $\mathcal{V}$
is {\it maximal} if none of its proper supersets has this property.
}
\end{definition}

\begin{proposition}\label{a111} {\rm Let $\mathfrak{X}=(X,\tau_1,\tau_2,\precsim)$ be a joincompact bitopological ordered space. 
Then, $\precsim$ has a
$\tau_1$-lower and $\tau_2$-upper semicontinous Richter-Peleg utility representation.}
\end{proposition}
\begin{proof} By Theorems \ref{a231}, \ref{a7} and Proposition \ref{a2}, $\mathfrak{X}$ has an increasing 
real-valued function $f$ on $X$ such that
$f$ is 
$\tau_{_1}$-lower semicontinuous and 
$\tau_{_2}$-upper semicontinuous. Therefore, if $x\precsim y$ the we have $f(x)\leq f(y)$.
On the other hand, if $x\prec y$, then $y\not\precsim x$ and by Proposition
\ref{a112} we have that $f(x)<f(y)$. It follows that 
$f$ is a $\tau_1$-lower and $\tau_1$-upper semicontinous Richter-Peleg utility representation of $\precsim$.
\end{proof}

According to Proposition \ref{a2}, 
a Richter-Peleg multi-utility representation $\mathcal{V}$ of $\precsim$ in a joincompact bitopological 
ordered space
$\mathfrak{X}=(X,\tau_1,\tau_2,\precsim)$
consists of bounded functions. In what follows, $\mathcal{L}_{_1}\mathcal{U}_{_2}(\mathfrak{X})$
(resp. $\mathcal{B}\mathcal{L}_{_1}\mathcal{U}_{_2}(\mathfrak{X})$)
denotes the family of (resp. bounded) $\tau_{_1}$-lower and 
$\tau_{_2}$-upper semicontinuous
functions of $\mathfrak{X}$ to 
$\mathfrak{I}=(\mathbb{R},\mathfrak{U},\mathfrak{L},\leq)$
and 
$\mathcal{R}\mathcal{P}\mathcal{B}\mathcal{L}_{_1}\mathcal{U}_{_2}(\mathfrak{X})$
denotes the members of $\mathcal{B}\mathcal{L}_{_1}\mathcal{U}_{_2}(\mathfrak{X})$
which are Richter-Peleg utility representations of $\precsim$.
Define the {\it supremum norm} on 
$\mathcal{B}\mathcal{L}_{_1}\mathcal{U}_{_2}(\mathfrak{X})$ as
\begin{center}
$\norm{f}_{_\infty}=\sup \{|f(x)|\ \vert\ x\in X\}$.
\end{center}
For each subfamily $\mathcal{V}$ of $\mathcal{L}_{_1}\mathcal{U}_{_2}(\mathfrak{X})$,
$cl_{_{\mathcal{U}n}}\mathcal{V}$ denotes the closure of 
$\mathcal{V}$
with respect to the usual (uniform) norm topology.

The next proposition is a dual version of Lemma 3.4 in \cite{kel}.

\begin{proposition}\label{lase}{\rm The class $\mathcal{L}_{_1}\mathcal{U}_{_2}(X)$
of real-valued functions which are $\tau_{_1}$-lower semicontinuous and 
$\tau_{_2}$-upper semicontinuous on a space $\mathfrak{D}=(X,\tau_{_1},\tau_{_2})$ is complete
with respect to the uniform norm on $X$.}
\end{proposition}

\begin{definition}\label{da1}{\rm Let $X$ be a set and let $\mathfrak{F}$ be a family of functions each of which has domain $X$. Then the family $\mathfrak{F}$ is {\it separating}
if, for each pair $x, y\in X$ of distinct elements of $X$, there exists an $f_{_{xy}}$ in
$\mathfrak{F}$ such that $f_{_{xy}}(x)\neq f_{_{xy}}(y)$.
}
\end{definition}

\begin{lemma}\label{a098}{\rm  Let $\mathfrak{X}=(X,\tau_1,\tau_2,\precsim)$ be a bitopological 
joincompact ordered space and
let $\mathcal{V}$ be a $\tau_{_1}$-lower and $\tau_{_2}$-upper semicontinuous 
Richter-Peleg multi-utility representation of $\precsim$.
Then, $cl_{_{\mathcal{U}n}}\mathcal{V}$ is also a 
a Richter-Peleg multi-utility representation of $\precsim$.}
\end{lemma}
\begin{proof} Let $f$ belongs to the closure of 
$\mathcal{V}$
with respect to the uniform norm topology on $\mathfrak{X}$. 
Then,
there exists a sequence $(f_{_n})_{_{n\in\mathbb{N}}}\in \mathcal{V}$ such that
$f_{_n}$ converges to $f$ in the uniform norm topology.
We prove that $f$ is a bounded, increasing, order preserving, $\tau_{_1}$-lower and $\tau_{_2}$-upper 
semicontinuous function.
Uniform convergence implies that for any $\varepsilon>0$ there is an $N_{_{\varepsilon}}\in\mathbb{N}$ such that 
$|f_{_n}(x)-f(x)|<\frac{\varepsilon}{3}$
for all $n\geq N_{_{\varepsilon}}$ and all $x\in X$. 
\par\smallskip\noindent
($\mathfrak{a}$) To prove that $f$ is bounded, let $\varepsilon=3$. 
Since each $f_{_n}$ is bounded on $X$, there exists $M_{_n}$ such that $|f_{_n}(x)|\leq M_{_n}$ for all $x\in X$. By uniform convergence, there exists an $N_{_3}\in\mathbb{N}$ such that $|f_{_n}(x)-f(x)|\leq 1$ for all $n\geq N_{_3}$. Then, 
$|f(x)|-|f_{_{N_{_3}}}(x)|\leq |f(x)-f_{_{N_{_3}}}(x)|\leq  1$ which implies that 
$|f(x)|\leq |f_{_{N_{_3}}}(x)|+1\leq M_{_{N_{_3}}}+1$.
\par\smallskip\noindent
($\mathfrak{b}$) To demonstrate that $f$ is a Richter-Peleg utility function, 
we begin by proving that if $a\sim b$ for some $a, b$ in X, then $f(a)=f(b)$.
Suppose to the contrary that 
$f(a)\neq f(b)$.
Since the sequence $(f_{_n})_{_{n\in\mathbb{N}}}$ converges uniformly to $f$, we have that
$f_{_n}(a)\to f(a)$, $f_{_n}(b)\to f(b)$ and for each $n\in \mathbb{N}$, $f_{_n}(a)=f_{_n}(b)$.
Therefore, $f_{_n}(a)$ converges 
to both $f(a)$ and $f(b)$, which contradicts the
Hausdorffness of $(X,\tau)$, where $\tau$ is the usual order topology on $\mathbb{R}$.
To prove that $f$ is order preserving, suppose to the contrary that $f(a)\geq f(b)$
for some $a$ and $b$, where $a\prec b$. Let $\varepsilon=f(a)-f(b)$.
Since the sequence $(f_{_n})_{_{n\in\mathbb{N}}}$ converges uniformly to $f$, 
we have that for each $n\geq N_{_\varepsilon}$ there holds
$\sup \{|f(x)-f_{_n}(x)|\vert x\in X\}<\frac{\varepsilon}{3}$.
It follows that $f(a)-f_{_n}(a)<\frac{\varepsilon}{3}$  and 
$f_{_n}(b)-f(b)<\frac{\varepsilon}{3}$ for $n>n_{_0}$. Therefore,
$f_{_n}(a)-f_{_n}(b)>-\frac{2\varepsilon}{3}+(f(a)-f(b))=\frac{\varepsilon}{3}>0$, 
contradicting the fact that $f_{_n}$ is order preserving.
Hence, $f$ is order preserving.

\par\smallskip\noindent
($\mathfrak{c}$) It remains to prove that
$f$ is $\tau_{_1}$-lower and $\tau_{_2}$-upper semicontinuous function.
Fix $x_{_0}$ in $X$ and let $x\in X$.
For any $n\in \mathbb{N}$ we have
\begin{center}
$f(x)-f(x_{_0})=(f(x)-f_{_n}(x))+(f_{_n}(x)-f_{_n}(x_{_0}))+(f_{_n}(x_{_0})-f(x_{_0})\leq$\\
$2\sup \{|f_{_n}(x)-f(x)|\vert x\in X\}+(f_{_n}(x)-f_{_n}(x_{_0}))$.
\end{center}
Therefore, for any $n$ we have,
\begin{center}
$f(x)-f(x_{_0})\leq 2\sup \{|f_{_n}(x)-f(x)|\vert x\in X\}+(f_{_n}(x)-f_{_n}(x_{_0}))$.
\end{center}
Since the previous inequality is true for an arbitrary value of $n$, we can choose 
a positive number $\varepsilon>0$, as small as we like. 
Therefore, there must be an $N>0$ for which we have 
$\sup \{|f_{_n}(x)-f(x)|\vert x\in X\}<\frac{\varepsilon}{3}$ for all $n>N$.
Because
$f_{_N}(x)$ is $\tau_{_2}$-upper semicontinuous, so for any choice of the number
$\varepsilon>0$ there is $\tau_{_2}$-open neighborhood $\mathcal{M}$ 
such that for all $x\in \mathcal{M}$,
$f_{_N}(x)-f_{_N}(x_{_0})<\frac{\varepsilon}{3}$.
Using this, we now have
\begin{center}
$f(x)-f(x_{_0})\leq 2\sup \{|f_{_n}(x)-f(x)|\vert x\in X\}+(f_{_n}(x)-f_{_n}(x_{_0}))<
2\frac{\varepsilon}{3}+\frac{\varepsilon}{3}=\varepsilon$
\end{center}
for all $x\in \mathcal{M}$. Therefore there exists a $\tau_{_2}$-open neighborhood 
$\mathcal{M}\subseteq X$ of $x_{_0}$ such that
for each $x\in \mathcal{M}$ we have $f(x)<f(x_{_0}+\varepsilon$.
It follows that $f$ is $\tau_{_2}$-upper semicontinous.
Similarly, interchanging the order of the difference of $f(x)$ and $f(x_{_0})$
we can prove that
$f$ is $\tau_{_1}$-lower semicontinous. 
Hence,
$cl_{_{\mathcal{U}n}}\mathcal{V}$ is also
a $\tau_{_1}$-lower and $\tau_{_2}$-upper semicontinous
Richter-Peleg multi-utility representation of $\precsim$ and this complete the proof.

Obviously, if members of $\mathcal{V}$ are bounded, so will members of $cl_{_{\mathcal{U}n}}\mathcal{V}$.
\end{proof}

\begin{proposition}\label{a812} {Let $\mathfrak{X}=(X,\tau_1,\tau_2,\precsim)$ be a 
joincompact bitopological ordered space and $\mathcal{L}$ be
a closed, with respect to uniform norm topology, sublattice of 
$\mathcal{B}\mathcal{L}_{_1}\mathcal{U}_{_2}(\mathfrak{X})$.
Then, a $\tau_1$-lower and $\tau_2$-upper semicontinous function $\phi$
belongs to
 $\mathcal{L}$
if and only if for each
$x, y\in X$ there exists $f_{_{xy}}\in \mathcal{L}$ such that 
$f_{_{xy}}(x)=\phi(x)$ and $f_{_{xy}}(y)=\phi(y)$.}
\end{proposition}
\begin{proof}{\it Necessity}:
If $\phi\in\mathcal{L}$, then $\phi$ itself satisfies the requirements of the proposition.
\par\noindent
{\it Sufficiency}: Suppose that $\phi$ is a $\tau_1$-lower and $\tau_2$-upper semicontinous function
on $\mathfrak{X}$.
Let $x, y \in X$ with $x\neq y$ and let $f_{_{xy}}\in \mathcal{L}$ such that
$f_{_{xy}}(x)=\phi(x)$, $f_{_{xy}}(y)=\phi(y)$. If $\phi \in \mathcal{L}$, then we have nothing to prove.
Otherwise
assume that a
$\phi \notin \mathcal{L}$. So, $f_{_{xy}}$ approximates $\phi$
in neighborhoods around $x$ and $y$. If $\phi(x)=\phi(y)$ then we can take a constant.
If not, since $\mathcal{L}$ is closed,
it suffices to show that for each $\varepsilon>0$ there exists $f\in \mathcal{L}$ such that
for all $z\in X$ we have
\begin{center}
$\varphi(z)-\varepsilon<f(z)<\varphi(z)+\varepsilon$,
\end{center}
or equivalently
\begin{center}
$\sup\{|\varphi(z)-f(z)|\vert z\in X\}<\varepsilon$,
\end{center}
for it will follow from this that $\norm{f-\phi}_{_{\infty}}<\varepsilon$.
Now, fix an $x\in X$, and let $y\in X$ vary. Put
\begin{center}
$U_{_y}=\{z\in X\vert \ f_{_{xy}}(z)<\phi(z)+\varepsilon\}=\{z\in X\vert \ f_{_{xy}}(z)-\phi(z)<\varepsilon\}$.
\end{center}
Since $f_{_{xy}}-\varphi$ is $\tau_{_1}\vee \tau_{_2}$-lower semicontinuous and 
$\tau_{_1}\vee \tau_{_2}$-upper semicontinuous we conclude that
 $f_{_{xy}}-\varphi$ is $\tau_{_1}\vee \tau_{_2}$-continuous.
Then, $U_{_y}$ is $\tau_{_1}\vee \tau_{_2}$-open set because $f_{_{xy}}-\varphi$ is $\tau_{_1}\vee \tau_{_2}$-continuous.
Also, $y\in U_{_y}$ and thus $\{U_{_y}\vert y\in X\}$ is a 
$\tau_{_1}\vee \tau_{_2}$-open cover of $X$. 
Since $\mathfrak{X}$ is joincompact, we have that $X$ is $\tau_{_1}\vee \tau_{_2}$-compact. Hence, there exists finitely many $y_{_1},...,y_{_n}\in X$ such that
\begin{center}
$X=\displaystyle\bigcup_{i\in\{1,...,n\}}U_{_{y_{_i}}}$.
\end{center}
Let $f_{_{xy_{_1}}}, f_{_{xy_{_2}}},...,f_{_{xy_{_n}}}$ are the functions of $\mathcal{L}$ which correspond to the sets $U_{_{y_{_1}}},U_{_{y_{_2}}},$
$...,U_{_{y_{_n}}}$ respectively.
Put
\begin{center}
$g_{_x}=f_{_{xy_{_1}}}\wedge f_{_{xy_{_2}}}\wedge...\wedge f_{_{xy_{_n}}}$.
\end{center}
Since  $\mathcal{L}$ is semi-vector lattice we have that $g_{_x}\in \mathcal{L}$.
On the other hand, $g_{_x}(x)=\phi(x)$ and $g_{_x}(z)<\phi(x)+\varepsilon$ for all
$z\in X$.
We next consider the open set 
\begin{center}
$V_{_x}=\{z\in X\vert \ g_{_x}(z)>\phi(z)-\varepsilon\}$.
\end{center}
Then, $V_{_x}$ is $\tau_{_1}\vee \tau_{_2}$-open set because $g_{_x}$ is 
$\tau_{_1}\vee \tau_{_2}$-continuous. Also, $x\in V_{_x}$ and thus 
$\{V_{_x}\vert x\in X\}$ is a $\tau_{_1}\vee \tau_{_2}$-open cover of $X$. Since $\mathfrak{X}$ is joincompact, we have that $X$ is $\tau_{_1}\vee \tau_{_2}$-compact.
Therefore, this open cover has a finite subcover $\{V_{_{x_{_1}}}, V_{_{x_{_1}}}, ..., V_{_{x_{_n}}}\}$.
We denote the corresponding
functions in 
$\mathcal{L}$ by $g_{_{x_{_1}}}, g_{_{x_{_2}}},...,g_{_{x_{_n}}}$,
and we define $f$ by $f=g_{_{x_{_1}}}\vee g_{_{x_{_2}}}\vee,...,\vee g_{_{x_{_n}}}$.
It is clear that $f\in\mathcal{L}$ with the
property that
\begin{center}
$\varphi(z)-\varepsilon<f(z)<\varphi(z)+\varepsilon$,
\end{center}
for all $z\in X$. Therefore, for each $\varepsilon>0$ we have 
$\norm{f-\phi}_{_{\infty}}<\varepsilon$ and
so our
proof is complete.
\end{proof}

\begin{definition}{\rm Let $\mathcal{L}=(X,\precsim)$ be a lattice. We say that 
$\mathcal{L}$
is a {\it generalized cone lattice} if it is closed under multiplication by scalars and addition
of constants (in short, $f\in \mathcal{L}$ and $a, b\in \mathbb{R}$ imply 
$a f+b\in\mathcal{L}$).}
\end{definition}

\begin{theorem}\label{a119} {Let $\mathfrak{X}=(X,\tau_1,\tau_2,\precsim)$ be a joincompact bitopological ordered space. 
Then, the set of all 
$\tau_{_1}$-lower and $\tau_{_2}$-upper semicontinous
Richter-Peleg utility representations of $\precsim$
is a separating and closed (with respect to uniform norm topology) 
generalized cone lattice of the set of $\tau_{_1}$-lower and $\tau_{_2}$-upper semicontinous
functions on $\mathfrak{X}$.
Conversely, given a joincompact  bitopological space
$\mathfrak{D}=(X,\tau_1,\tau_2)$ and a
separating and closed (with respect to uniform norm topology) 
generalized cone lattice $\mathcal{U}$ of 
$\tau_{_1}$-lower and $\tau_{_2}$-upper semicontinous functions on $\mathfrak{X}$, there is one and only one way
to turn $\mathfrak{D}$ into a joincompact bitopological ordered space $\mathfrak{X}=(X,\tau_1,\tau_2,\precsim)$ where $\mathcal{U}$ is the set of all increasing
$\tau_{_1}$-lower and $\tau_{_2}$-upper semicontinous
Richter-Peleg utility representations of $\precsim$.
}
\end{theorem}
\begin{proof}
Let $\mathfrak{X}=(X,\tau_1,\tau_2,\precsim)$ be a bitopological ordered space  and let 
$\mathcal{V}$ be the set of all 
Richter-Peleg 
utility representation
of
$\precsim$. 
We
will prove that $\mathcal{V}$ satisfies the properties given in the statement of the present proposition. 
\par\smallskip\par\noindent
(i) $\mathcal{V}$ {\it separates points on $X$}.  Indeed, let $x, y\in X$ with $x\not\sim y$. 
Then, $x\not\precsim y$ or $y\not\precsim x$.
By Proposition \ref{a112} we have that $f(y)<f(x)$ or $f(x)<f(y)$ for some
$f\in \mathcal{V}$. Therefore, in any case there exists $f\in \mathcal{V}$
such that $f(x)\neq f(y)$, so $\mathcal{V}$
separates points. 
\par\smallskip\par\noindent
(ii) $\mathcal{V}$ {\it is closed under the uniform norm topology.} 
Let $f\in cl_{_{\mathcal{U}n}}\mathcal{V}$.
By Lemma \ref{a098}
$cl_{_{\mathcal{U}n}}\mathcal{V}\supset \mathcal{V}$
is a Richter-Peleg 
multi-utility representation of $\precsim$, a contradiction to the maximal character of $\mathcal{V}$. Hence, $cl_{_{\mathcal{U}n}}\mathcal{V}=\mathcal{V}$
which implies that $\mathcal{V}$ is closed with respect to uniform norm topology.
\par\smallskip\par\noindent
(iii) $\mathcal{V}$ {\it is a lattice}. Let $f, g \in \mathcal{V}$ and fix 
an $\varepsilon>0$ and a $x_{_0}\in X$. Since $f, g$ are $\tau_{_1}$-lower and $\tau_{_2}$-upper semicontinuous, there are
$\tau_{_1}$-open neighborhods $O, O^{\prime}$
and $\tau_{_2}$-open neighborhods 
$P, P^{\prime}$ such that for all $y\in O, y\in O^{\prime}, y\in P$ and $y\in P^{\prime}$
we have, respectively, that $f(y)<f(x_{_0})+\varepsilon$, 
$f(x_{_0})-\varepsilon <f(y)$,
$g(y)<g(x_{_0})+\varepsilon$ and 
$g(x_{_0})-\varepsilon <g(y)$. 
If $\widetilde{K}=O\cap P$ and 
$\widetilde{K^{\prime}}=O^{\prime}\cap P^{\prime}$, then for each $y\in \widetilde{K}$ we have
that $\sup \{f(y),g(y)\}<\sup\{f(x_{_0}),g(x_{_0})\}+\varepsilon$ and 
for each 
$y\in \widetilde{K^{\prime}}$ we have
that $\inf\{f(x_{_0}),g(x_{_0})\}-\varepsilon<\inf \{f(y),g(y)\}$. 
So, $f\vee g$ is a
$\tau_{_1}$-lower and $\tau_{_2}$-upper semicontinuous function.
On the other hand, since $f, g$ are increasing and order preserving, for each $x_{_1}, x_{_2}\in X$ with
$x_{_1}\precsim x_{_2}$ we have that $\sup \{f(x_{_1}),g(x_{_1})\}\precsim \sup \{f(x_{_2}),g(x_{_2})\}$ and for each $x_{_1}, x_{_2}\in X$ with
$x_{_1}\prec x_{_2}$ we have that $\sup \{f(x_{_1}),g(x_{_1})\}\prec \sup \{f(x_{_2}),g(x_{_2})\}$. 
Thus $f\vee g$ is increasing and order preserving. Since 
$\mathcal{V}$ is maximal with respect to set inclusion, we conclude that 
$f\vee g\in \mathcal{V}$.
Similarly we prove that 
$f\wedge g\in \mathcal{V}$ which implies that $\mathcal{V}$ is a lattice.
\par\smallskip\par\noindent
(iv) $\mathcal{V}$ {\it is a generalized cone lattice}. 
Let
$f$ be an increasing $\tau_{_1}$-lower and $\tau_{_2}$-upper semicontinuous
function on $\mathcal{V}$ and let $\lambda, \kappa\in\mathbb{R}$ with $\lambda\geq 0$.
Clearly, $\lambda f+\kappa$ is increasing. It rermains to prove that
$\lambda f+\kappa\in \mathcal{V}.$
If $\lambda=0$, then $\lambda f+\kappa=\kappa\in \mathcal{V}.$
Otherwise, $\lambda>0$. Fix an $\varepsilon>0$ and an $x_{_0}\in X$.
Since $f$ is $\tau_{_1}$-lower semicontinuous there exists a 
$\tau_{_1}$-open neighborhood $U_{_{x_{_0}}}$ of $x_{_0}$ such that for each 
$y\in U_{_{x_{_0}}}$, $f(y)<f(x_{_0})+\frac{\varepsilon}{\lambda}$ holds.
Let
\begin{center}
$A_{_y}=\{y\in X\vert \lambda f(y)+\kappa<\lambda f(x_{_0})+\kappa+\varepsilon\}$.
\end{center}
Then, $A_{_y}=\{y\in X\vert \lambda f(y)<\lambda f(x_{_0})+\varepsilon\}=
\{y\in X\vert f(y)<f(x_{_0})+\frac{\varepsilon}{\lambda}\}=U_{_{x_{_0}}}$.
It follows that $\lambda f+ \kappa$ is $\tau_{_1}$-lower semicontinuous. Similarly, we prove that $\lambda f+ \kappa$ is $\tau_{_2}$-upper semicontinuous.
\par\smallskip\par\noindent

Conversely, let $\mathcal{U}$ be a
separating and closed (with respect to uniform norm topology) 
generalized cone lattice of $\mathfrak{D}$. The constant functions are in $\mathcal{U}$, so this isn't empty.
Define a order $\precsim$ on $\mathfrak{D}$ as follows:
\begin{center}
$x\sim y$ if $f(x)=f(y)$ and $x\prec y$ if $f(x)<f(y)$
for all $x, y\in X$ and $f\in \mathcal{V}$.
\end{center}
We prove that $\precsim$ is a $\tau_{_1}\times \tau_{_2}$-closed subset of $X\times X$
and thus $\mathfrak{X}=(X,\tau_1,\tau_2,\precsim)$
is a bitopological ordered space. 
Indeed,
Let $a, b\in X$ with $a\not\precsim b$. 
Since $\mathcal{V}$ is separating, there exists at least one 
$f\in\mathcal{U}$ such that $f(a)\neq f(b)$ ($f(a)< f(b)$ or $f(b)< f(a)$).
Let $\mathfrak{F}_{_{a,b}}=\{f\in \mathcal{V}\vert f(a)\neq f(b)\}$.
Then, there exists $f^{\ast}\in \mathfrak{F}_{_{a,b}}$
such that $f^{\ast}(b)<f^{\ast}(a)$, because otherwise, for all $f$ in $\mathcal{V}$ 
there holds $f(a)\leq f(b)$. But then, $a\precsim b$ which is a contradiction to
$a\not\precsim b$.
By the density of reals,
there exists
$r>0$ for which 
$f^{\ast}(b)<r<f^{\ast}(a)$ holds. Let $U_{_b}=(f^{\ast})^{-1}(]-\infty,r[)$ and $U_{_a}=(f^{\ast})^{-1}(]r,\infty[)$. Then,
$U_{_b}$ is a  $\tau_{_2}$-decreasing neighborhood of $b$ and 
$U_{_a}$ is a $\tau_{_1}$-increasing neighborhood of $a$ such that
$U_{_a}\cap U_{_b}=\emptyset$. Hence, by Theorem \ref{a3} we conclude that
$\precsim$ is a $\tau_{_1}\times \tau_{_2}$-closed subset of $X\times X$ which implies that 
$\mathfrak{X}$ is a bitopological ordered space. It remains to prove that $\mathcal{U}$
is the set of all Richter-Peleg utility representations of $\precsim$. 
By definition, all $f\in \mathcal{U}$ are 
Richter-Peleg utility representations of $\precsim$.

Let $\varphi$ be a Richter-Peleg utility representation of $\precsim$ on $\mathfrak{X}$
and let $x, y\in X$. We prove that there exists $h\in \mathcal{U}$
such that $h(x)=\varphi(x)$ and $h(y)=\varphi(y)$, and thus, since  
$\mathcal{U}$ is closed with respect to uniform norm topology, by Proposition \ref{a812} we have that $\varphi\in \mathcal{U}$.
If $\varphi(x)=\varphi(y)=\mu$, then since $\mathcal{U}$ is non-empty,
by taking $\lambda=0$ we observe that the function $h=0.f+\mu$ with $f\in \mathcal{U}$ belongs to $\mathcal{U}$ and $h(x)=\varphi(x)$
and $h(y)=\varphi(y)$.
Now, assume that $\varphi(x)<\varphi(y)$. Then, since $\varphi$ is order preserving we conclude that $y\not\precsim x$. By corollary \ref{a112} there exists an increasing $\tau_{_1}$-lower semicontinuous and $\tau_{_2}$-upper semicontinuous function $f$
such that $f(x)<f(y)$.
Choose the real numbers $\lambda, \kappa$ such that 
$\lambda f(x)+\kappa=\varphi(x)$ and 
$\lambda f(y)+\kappa=\varphi(y)$. Since 
$\lambda=\frac{\varphi(y)-\varphi(x)}{f(y)-f(x)}>0$
and $\mathcal{U}$ is closed with respect to uniform norm topology, 
by Proposition \ref{a812}, we conclude that
$\varphi\in \mathcal{U}$.
It follows that 
$\mathcal{U}$ is the set of all 
$\tau_{_1}$-lower and $\tau_{_2}$-upper semicontinous
Richter-Peleg utility representations of $\precsim$.
\end{proof}

\section{The case of preordered sets}
Let $\mathcal{P}=(X,\precsim)$ be a preordered set.
A subset $D$ of $\mathcal{P}$ is {\it directed} provided it is nonempty, and every finite subset of $D$ has an upper bound in $D$. 
We use $\bigvee^{\mathcal{P}} D$ (resp. $\bigwedge^{\mathcal{P}} D$)
to represent the supremum (resp. infimum) of $D$ if $D$ is a directed set and the supremum (resp. infimum) exists in the preordered set.
A {\it directed
complete preordered set} is a preordered set such that each of
its directed subsets has a supremum. Note that {\it dcpo} normally stands for directed complete partially ordered set (poset).
A directed complete poset which is a lattice is called a {\it directed complete lattice}.
(ii) A {\it complete lattice} $\mathcal{L}$ is a poset in which every subset has a supremum and an infimum.
An {\it ideal} of $\mathcal{P}$ is a directed lower set.
The topology on $\mathcal{P}$ generated by $\{X\setminus \displaystyle\uparrow \{x\}\vert
x\in X\}$
is called the {\it lower topology} and denoted by $\omega(\mathcal{P})=\omega$. 
Let $x, y$ be elements of $X$. 
We say that $x$ is {\it way-below} $y$, written $x\ll y$, if for any directed subset $D$ with
$\bigvee^{\mathcal{P}} D$ exists and $y\precsim \bigvee^{\mathcal{P}} D$, implies that $x\precsim d$ for some $d\in D$.
If  $x\ll y$ then $x\precsim y$ (consider the directed set $D=\{y\}$).
A poset is a {\it continuous poset}  if every element is the join of a directed set of those elements which are way-below it.
A lattice $\mathcal{L}$ is called a {\it continuous lattice} if it is a complete lattice and if every element is the join of those elements which are way-below it.
A subset $A$ of $X$ is called {\it Scott-open} if $A$ is an upper set ($A=\displaystyle\uparrow A$) and for a directed set $D$ with 
$\bigvee^{\mathcal{P}} D\in A$, we have $d\in A$ for some
$d\in D$.
The Scott-open sets satisfy the axioms of a
topology, which we call the {\it Scott topology} and we denote by
$\sigma(\mathcal{P})=\sigma$.
If $\mathcal{P}$ is a continuous poset, then
all sets of the form
$\turnnw{\twoheadleftarrow}{x}=\{y\vert x\ll y\}$
are Scott-open sets. 
The scott topology $\sigma$ on a poset $\precsim$ is always compatible, that is, $\sigma=\precsim_{\sigma}$.
The Scott
topology of the natural order of the real line is the topology of
lower semicontinuity (the topology for which the non-trivial open
sets are the intervals  $(a,\infty)$).
The common refinement $\sigma(\mathcal{P})\vee \omega(\mathcal{P})$ of the Scott and lower topologies
is called the {\it Lawson topology} and is denoted
by $\lambda(\mathcal{P})$(cf. \cite[Definition III-11.5]{GH}).
Given two preordered sets $\mathcal{P}=(X,\precsim)$ and $\mathcal{Q}=(Y,\sqsubseteq)$,
a function $f: \mathcal{P}\to \mathcal{Q}$
is an {\it order embedding} if
for all $x, y\in X$, one has $x\precsim y$ if and only if $f(x) \sqsubseteq f(y)$.
In case of a poset $\mathcal{P}$, this condition forces $f$ to be one-to-one. 
An order embedding is a type of monotone function in order theory that allows one preordered set to be included in another. 
An {\it extension} of a preordered set $\mathcal{P}$
is a pair $(f,\mathcal{Q})$, where $\mathcal{Q}$ is a preordered set and 
$f: \mathcal{P}\to \mathcal{Q}$
is an order embedding. 
A completion of a preordered set $\mathcal{P}$ is an extension 
$(f,\mathcal{Q})$
of
$\mathcal{P}$
such that $\mathcal{Q}$ is a complete lattice.
A subset $I$ of $X$ is an {\it ideal}, if the following conditions hold:
(i) $I$ is non-empty; (ii)
for every $x\in I$, any $y\in X$ and $y\precsim x$ implies that $y\in X$ ($I$ is a lower set), and
for every $x, y\in I$, there is some element $z\in I$, such that 
$x \precsim z$ and $y\precsim z$ ($I$ is a directed set).
The smallest ideal that contains an element $x\in X$ is called a {\it principal ideal}.
This is denoted by
${\displaystyle \downarrow x}=\{y\in X\vert y\precsim x\}$.
Given a subset $A$ of a preordered set $\mathcal{P}=(X,\precsim)$, we denote by 
$A^{\uparrow}$ and $A^{\downarrow}$ the sets of all upper and lower bounds of $A$, respectively.
There are various definitions of a cut of a preordered set $X$ yielding a completion of $X$.
MacNeille \cite{MN} has introduced the famous ``completion by cuts'' for arbitrary preordered sets.
A {\it cut} is a pair $(A,B)$ such that $A=B^{\downarrow}$ and $B=A^{\uparrow}$.
The collection of all cuts, ordered by $(A,B)\leq (C,D)$ if and only if $A\subseteq C$ and $B\subseteq D$
is a complete lattice, called the {\it Dedekind-MacNeille completion} of $\mathcal{P}$.
Any isomorphic copy of Dedekind-MacNeille completion is referred to as the normal completion of $\mathcal{P}$.
Normal completions can be characterized in a number of ways. Because any part of a cut determines the other, 
working with lower cuts is generally more convenient. The normal completion of $\mathcal{P}$ by lower cuts is defined as follows:
Define $A^{\delta}=(A^{\uparrow})^{\downarrow}$,
then the {\it lower cuts completion} of $X$ consists of all subsets $A$ for which
$A^{\delta}=A$.
If $\delta(X)=\{A^\delta\vert A\subseteq X\}$, then
$A^{\delta}\precsim^{\delta}B^{\delta}$ in completion if and only if $A^{\delta}\subseteq B^{\delta}$ as sets. 
The lower cuts completion $(\delta(X), \precsim^{\delta})$ of $\mathcal{P}$
is a complete lattice (\cite[Lemma 1]{abi}).
If $(X,\precsim)$ is a qoset,
then each element $x\in X$  corresponds to its principal ideal
$\downarrow x$
into the lower cuts completion $\delta(X)$. 
The lower cut completion shall be referred to as the {\it normal completion of} $\mathcal{P}$ in the following paragraphs.

Ern\'{e} \cite{ern2} introduced a new way-below relation and 
on the basis of this notion, he defines the concept of precontinuous preordered sets for arbitrary preordered sets. 
The work of Ern\'{e} was influenced by the use of Frink ideals \cite{fri} instead of directed lower sets.

Formally, this notion is defined as follows.
\begin{definition}\label{pas}{\rm (\cite{fri}) A subset $I$ of a preordered set $\mathcal{P}=(X,\precsim)$ is called a {\it Frink ideal} in $X$ if $\delta (Z)\subseteq I$ for all finite
subsets $Z\subseteq I$.
In what follows, $Fid(X,\precsim)$
denote the set of all Frink ideals of $\mathcal{P}$.
We say that
$\mathcal{P}$ is {\it precontinuous} if and only if for each $x\in X$ there is a smallest $I\in Fid(X,\precsim)$
such that  $x\in I^{\uparrow\downarrow}$. }
\end{definition}

\begin{remark}\label{opa2}{\rm In a continuous lattice $\mathcal{L}$, $\turnnw{\twoheadrightarrow}{x}$ 
is automatically directed and thus we may write the axiom of approximation, which defines continuous lattices, as 
\begin{center}
$x=sup \turnnw{\twoheadrightarrow}{x}=sup\{u\in \mathcal{L}\vert u\ll x\}$ 
\end{center}
or as
\begin{center}
whenever $x\not\precsim y$, then there is a $u\ll x$ with $u\not\precsim y$.
\end{center}

}
\end{remark}

\begin{proposition}\label{vjh}{\rm Let $\mathcal{L}=(X,\precsim)$ be a continuous lattice. Then, 
$\mathfrak{X}=(X,\sigma(\mathcal{L}),$
$\omega(\mathcal{L}),\precsim)$ is a joincompact bitopological ordered space.}
\end{proposition}
\begin{proof} By \cite[Theorems III-1.9 and III-1.10]{GH}, $\sigma(\mathcal{L})\vee
\omega(\mathcal{L})=\lambda(\mathcal{L})$ is a compact topology.
To show that $\precsim$ is closed,
suppose that $x\not\precsim y$ for some $x, y\in X$. 
By the remark \ref{opa2} (see also \cite[Definition I-1.6]{GH}), there exists $z\in X$
such that $z\ll x$ and $z\not\precsim y$. Then, $X\setminus \uparrow{z}$ is a decreasing
$\omega$-open neighborhood of $y$ and 
$\turnnw{\twoheadleftarrow}{z}$
is an increasing Scott-open neighborhood of $x$ such that 
$\turnnw{\twoheadleftarrow}{z}\cap (X\setminus \uparrow{z})=\emptyset$.
It follows that $\precsim$ is $\sigma\times\omega$-closed subset of $X\times X$ and 
so $(X,\sigma,\omega,\precsim)$ is a bitopological ordered space. 
\end{proof}

\begin{remark}\label{opa1}{\rm
In $[0,1]$, $x\ll y$ if and only if $x<y$ or $x=y=0$. Thus this ordered set is a continuous lattice (see \cite[Example 3.2]{HKMS}). 
The {\it lower topology} on $[0,1]$ is $\mathfrak{L}=\{[0,a)\vert 0<a\precsim 1\}\cup \{\emptyset, [0,1]\}$ and the {\it upper topology} is 
$\mathfrak{U}=\{(a,1]\vert 0\precsim a<1\}\cup \{\emptyset, [0,1]\}$.
A function $f$ from a topological space $X$ into the real unit interval $[0,1]$ is Scott-continuous if and only if it is lower semicontinuous in the sense of classical analysis (see \cite[Example 1-5.9]{GW}). On the other hand, 
if $X$ is a topological space, then 
$f: X\to [0,1]$ is 
upper semicontinuous if and only if it is a continuous map of 
$f: X\to [0,1]$
with the lower topology (see \cite[Page 295]{HKMS}).
Hence, if
$(X,\sigma,\omega)=(X,\sigma(\mathcal{L}),\omega(\mathcal{L}))$ is the bitopological space generated by a continuous lattice
$\mathcal{L}=(X,\precsim)$, 
without loss of generality, by Proposition \ref{a2},
we may assume that
the notion of $\sigma$-lower  semicontinuous function
coincides with the notion of Scott-continuous function and the notion of
$\omega$-upper semicontinuous function
coincides with the notion of $\omega$-continuous function.
The following notation will be used in the sequel: 
$\mathfrak{I}=([0,1], \sigma, \omega)$.
}
\end{remark}

\begin{theorem}\label{kjh}{\rm Let $\mathcal{L}=(X,\precsim)$ be a continuous lattice. 
Then,  the set of all Scott and $\omega$-continuous functions is a 
Richter-Peleg 
multi-utility representation of $\precsim$.}
\end{theorem}
\begin{proof} Suppose that $\mathcal{L}=(X,\precsim)$ is a continuous lattice. 
By Proposition \ref{vjh},  $\mathfrak{X}=(X,\sigma(\mathcal{L}),\omega(\mathcal{L}),\precsim)$
is a joincompact bitopological ordered space.
Let $\mathcal{V}$ be the set of all 
Scott and $\omega$-continuous
Richter-Peleg utility representation of $\precsim$.
By Proposition \ref{a111} we have $\mathcal{V}\neq \emptyset$.
Therefore,  $\mathcal{V}$ is a Scott and $\omega$-continuous
Richter-Peleg 
multi-utility representation of $\precsim$.
\end{proof}

A function $p: X\to Y$ is a {\it quotient map} if satisfy the following conditions: 
(i) $p$ is surjective; (ii) $p$ is continuous (i.e. $U$ is open in $Y$ implies that
$p^{-1}(U)$ is open in $X$), and (iii) $V\subseteq Y$ and $p^{-1}(V)$ open in $X$ implies that $V$ open in $Y$.
In this case we say the map $f$ is a {\it quotient map}.
In fact, if $f$ is a quotient map, then,
\begin{center}
$U$ is open in $Y$ if and only if $p^{-1}(U)$ is open in $X$.
\end{center}
The topology produced by $p$ is known as {\it quotient topology}.
If $\approx$ is an equivalence relation on a topological space $(X,\tau)$,
then the quotient set by this equivalence relation $\approx$ will be denoted by 
$X_{_\approx}$, and its elements (equivalence classes) by $[x]$. Let the {\it projection map }
$\pi: X\to X_{_\approx}$
which carries each point of $X$ to the element of $X_{_\approx}$
that contains it. The projection map is a quotient map.
The quotient topology defined by $\pi$ is denoted by $\tau_{_\approx}$ and
the topological space $(X_{_\approx},\tau_{_\approx})$ is called the {\it quotient space}
of $X$ determined by $\approx$.
Thus, the typical open set in
$X_{_\approx}$
is a collection of equivalence classes whose union is an open set in $X$.
A continuous map $f: X\to Y$ {\it respects the equivalence relation} $R$ if equivalent points have identical
images, that is if $xRy$ implies $f(x)=f(y)$.

Let $\precsim$ be a preorder on a set $X$.
Define an equivalence class $\approx$ on $X$ by 
$x \approx y$ if and only if $x\precsim y$ and $y\precsim x$. Then, $[x]=\{y\in X\vert $
$x\precsim y$ and $y\precsim x$ $\}$.
We can define a partial order $\sqsubseteq$ on $X_{_\approx}$ by:
\begin{center}
$[x] \sqsubseteq [y]$ if and only if $[x]=[y]$ or there are $x^{\prime}\in [x], y^{\prime}\in [y]$ such that $x^{\prime}\precsim y^{\prime}$.
\end{center}
Clearly,
\begin{center}
$x\precsim y$ if and only if $[x]\sqsubseteq [y]$.
\end{center}

In the following, the symbols $X_{_\approx}$ and $\sqsubseteq$,
will always stand for the notions
that have
been defined just above.

\begin{remark}\label{olymb}{\rm Let $\mathcal{P}=(X,\precsim)$ be a preordered set and 
$\mathcal{Q}=(X_{_\approx},\sqsubseteq)$ be the poset, as it has been defined just
before. Let $Fid(X_{_\approx},\sqsubseteq)$ denote the set of all Frink ideals in 
$\mathcal{Q}$.
If $x, y\in X$, we write $[x] \ll^{\ast}_{_e} [y]$, if 
for each $[I]\in Fid(X_{_\approx},\sqsubseteq)$, $[y]\in [I]^{\delta}$ implies $[x]\in [I]$. 
}
\end{remark}

The following three propositions are useful in the proof of main theorem.

\begin{proposition}\label{tame}{\rm Let $\mathfrak{X}=(X,\tau_1,\tau_2,\precsim)$
be a topological preordered space and let
$\mathfrak{X}_{_\approx}=(X_{_\approx},\mathfrak{t}_{_1},\mathfrak{t}_{_2},\sqsubseteq)$ be the quotient space of $\mathfrak{X}$.
Then, $\mathfrak{X}$ is a topological preordered space if and only if 
$\mathfrak{X}$ is a topological ordered space.
}
\end{proposition}
\begin{proof} By Theorem \ref{a3}($\mathfrak{a}$) and by definition of the
quotient topology, it is easy to show that:\
$\precsim$ is $\tau_1\times \tau_2$-closed if and only if
$\sqsubseteq$ is
$\mathfrak{t}_{_1}\times \mathfrak{t}_{_2}$-closed. The rest is obvious.
\end{proof}

\begin{proposition}\label{tfte}{\rm Let $\precsim$ be a preorder on a set $X$.
Then, $(X,\precsim)$ is precontinuous if and only if $(X_{_\approx},\sqsubseteq)$ is precontinuous.
}
\end{proposition}
\begin{proof}
It is an immediate consequence of
Definition \ref{pas} and Remark \ref{olymb}.
\end{proof}

Scott's continuous functions are monotonic, so $x\approx y$ implies $f(x)=f(y)$ 
for a Scott continuous function $f$, that is, $f$ respects $\approx$.

The following proposition is an immediate consequence of \cite[Theorem 2.82]{mol}

\begin{proposition}\label{a1me}{\rm Let $\mathcal{P}=(X,\precsim)$ be a preordered set and let
$p$ be the quotient map of $X$ to $X_{_\approx}$.
Let $f: X\to [0,1]$
be a Scott and $\omega$-continuous Richter-Peleg representation of $\precsim$.
Then, $\sqsubseteq$ induces a 
Scott and $\omega$-continuous Richter-Peleg representation 
$\widetilde{f}: X_{_\approx}\to [0,1]$
such that $\widetilde{f}\circ \pi =f$.
Conversely, if $\widetilde{f}: X_{_\approx}\to [0,1]$ is 
a Scott and $\omega$-continuous Peleg-Richter representation 
of $\sqsubseteq$, then $\precsim $ induces a
Scott and $\omega$-continuous Peleg-Richter representation 
$f: X\to [0,1]$
such that $\widetilde{f}\circ \pi=f$.
}
\end{proposition}

The following corollary is an immediate consequence of Theorem \ref{kjh},
Remark \ref{olymb} and Propositions \ref{tame}, \ref{tfte} and \ref{a1me}.

\begin{corollary}\label{lant}{\rm Let $\mathfrak{X}=(X,\tau_1,\tau_2,\precsim)$
be a topological preordered space and let
$\mathfrak{X}_{_\approx}=(X_{_\approx},\mathfrak{t}_{_1},\mathfrak{t}_{_2},\sqsubseteq)$ be the quotient space of $\mathfrak{X}$.
Then, the set of 
Scott and $\omega$-continuous Peleg-Richter representation of $\precsim$
in
$\mathfrak{X}$ is isomorphic to the
set of 
Scott and $\omega$-continuous Peleg-Richter representation of $\sqsubseteq$
in
$\mathfrak{X}_{_\approx}$.
}
\end{corollary}

The following lemma is a result of Ern\'{e} in \cite{ern2}.

\begin{lemma}\label{aeks}{\rm Let $\mathcal{P}=(X,\precsim)$ be a poset. 
Then, the normal completion of $\mathcal{P}$ is a continuous lattice
if and only if $\mathcal{P}$
is precontinuous.
}
\end{lemma}

The following lemma is Theorem 1-3.12 in \cite{GW}.

\begin{lemma}\label{GW}{\rm A complete lattice $\mathcal{L}=(X,\precsim)$ is continuous if and only if there is an injection of $\mathcal{L}$ into some power $[0,1]^I$ of the unit interval preserving arbitrary meets and directed joins.}
\end{lemma}

\begin{definition}\label{cgs}{\rm 
Let $\{f_{_a}\vert a\in A\}$ be a collection of functions.
Then, the {\it evaluation map} $e: X \longrightarrow \displaystyle\prod_{a\in A}f_{_a}$
induced by the collection $\{f_{_a}\vert a\in A\}$
is defined as follows: for each $x\in X$, $e(x)=(f(x))_{_{a}}$.
That is, for each $x\in X$, $e(x)$ is the point in $\displaystyle\prod_{a\in A}f_{_a}$
whose ath coordinate is $f_{_a}(x)$ for each $a\in A$.}
\end{definition}

\begin{theorem}\label{gjh}{\rm Let $\mathcal{L}=(X,\precsim)$ be a preordered set.
Then,  $\precsim$ has the set of Scott and $\omega$-continuous
functions from $X$ into the real unit interval $[0,1]$ 
as a Richter-Peleg 
multi-utility representation if and only if $\precsim$ is precontinuous.}
\end{theorem}
\begin{proof} 
Let $X_{_\approx}$ be the quotient space of $X$. Since 
$(X,\precsim)$ is precontinuous,  by Proposition \ref{tfte}, we have that 
$(X_{_\approx},\sqsubseteq)$ is precontinuous.
Let
$(\delta(X_{_\approx}), \sqsubseteq^{\delta})$
be the MacNeille completion of $(X_{_\approx},\sqsubseteq)$. 
By Lemma \ref{aeks},
$(\delta(X_{_\approx}), \sqsubseteq^{\delta})$ is a continuous lattice. 
Let 
\begin{center}
$\mathcal{V}=\{g: \delta(X_{_\approx})\to [0,1]\vert $
$g$
is a
Scott and $\omega$-continuous Richter-Peleg representation of
$\sqsubseteq^{\delta}\}$.
\end{center}
By Theorem \ref{kjh},
we have that $\mathcal{V}$ is a Richter-Peleg multi-utility representation of
$\sqsubseteq^{\delta}$. Let $\varphi^{\delta}$ be the usual embedding of  
$(X_{_\approx},\sqsubseteq)$ to
$(\delta(X_{_\approx}), \sqsubseteq^{\delta})$.
Clearly, $\varphi^{\delta}$ preserves all existing joins and meets (see \cite[Theorem 1]{abi}).
 Then, for each $g\in \mathcal{V}$, $g\circ \varphi^{\delta}$
 is a Scott and 
 $\omega$-continuous Richter-Peleg function of $X_{_\approx}$ to $[0,1]$.
 Let
 \begin{center}
 $\widetilde{\mathcal{V}}=\{\widetilde{f}=g\circ \varphi^{\delta}\vert g\in \mathcal{V}\}$.
  \end{center}
It follows that $\widetilde{\mathcal{V}}$ is a Richter-Peleg multi-utility representation of
$\sqsubseteq$. 
Define,
 \begin{center}
 $\widehat{\mathcal{V}}=\{h=\widetilde{f}\circ \pi\vert \widetilde{f}\in
 \widetilde{\mathcal{V}}\}$, where $\pi$ is the projection map of $X$ to $X_{_\approx}$\}
 \end{center}
If $x, y\in X$, then
\begin{center}
$x\sim y\Rightarrow h(x)=h(y)$ and $x\prec y\Rightarrow h(x)<h(y)$.
\end{center}
 Then, by Proposition \ref{a1me}, we have that $h$
is a 
Scott and 
$\omega$-continuous function of $(X,\precsim)$ to $[0, 1]$.

Hence, $h$ is a Richter-Peleg utility representation of
$(X,\precsim)$.
It follows that
$\widehat{\mathcal{V}}$ is a Richter-Peleg multi-utility representation of
$(X,\precsim)$.

Conversely, suppose that $\mathcal{V}$ is a Richter-Peleg multi-utility representation of $\precsim$ such that each $f\in \mathcal{V}$ is a Scott and $\omega$-open
continuous function
from $X$ into the real unit interval $[0,1]$. 
For every $f\in \mathcal{V}$ we define the  preorder $\precsim_{_f}$ such that
\begin{center}
$\precsim_{_f}=\{(x,y)\in X\times X\vert f(x)\leq f(y)\}$.
\end{center}
Then,
the monotonicity and Scott and $\omega$-continuity
of $f$ implies that $\precsim_{_f}$ is a complete binary relation and
for every point $x\in X$, the set $\{y\in X\vert  y\precsim_{_f} x\}$ is $\sigma$-closed and
the set $\{z\in X\vert  x\precsim_{_f} z\}$
is $\omega$-closed in $X$.
It, thus, follows that $\precsim=\displaystyle\bigcap_{f\in \mathcal{V}}\precsim_{_f}$
is $\sigma\times\omega$-closed preorder and therefore $(X,\sigma,\omega,\precsim)$ 
is a bitopological ordered space.
As it is known, in $[0, 1]$, $x\ll y$, if and only if $x<y$ or $x=y=0$. Thus, $([0, 1],\ll)$ is a continuous lattice and thus
$\ll$ is the specialization order of its Scott topology.
Therefore,
$\mathbb{I}=([0,1],\mathfrak{L},\mathfrak{U})$, where $\mathfrak{L}, \mathfrak{U}$ are the lower and upper
topologies corresponds to the unit interval for bitopological ordered spaces (see Example \ref{ex0}).
For each $x\in X$, let
$h(x)=\{y\vert y\in X\ {\rm and}\ y\precsim x\}$. Then, the intersection of any $\{h(x)\vert x\in X\}$ is a lower cut in $X$
(see \cite[Definition 1]{abi}).
By \cite[Theorem 1]{abi}, $h(x)$ is a one -to- one mapping from $(X,\precsim)$
into $(\delta(X), \precsim^{\delta})$ which
extends $\precsim$ and preserves arbitrary meets and directed  joins.
 Define
\begin{center}
$\widehat{f}: \delta(X)\longrightarrow \mathbb{I}$ by $\widehat{f}(x)=\displaystyle\inf_{y\in x^{^\uparrow}}f(y)$.
\end{center}
Since $\precsim^{\delta}$ extends $\precsim$, by monotonicity of $f$ (it is Scott continuous), for each $y\in x^{^\uparrow}$
we have that $x\precsim y$. Therefore,
\begin{center}
$x \precsim^{\delta} y$ if and only if $\widehat{f}(x)\leq \widehat{f}(y)$.
\end{center}
The subset $[0,1]$ of reals is complete with respect to upper and lower topologies, thus by remark \ref{opa1}
we concludes that
$h$
is an
injection of $\delta(X)$ into $\mathbb{I}$ which preserves arbitrary meets and directed  joins.
On the other hand, if $x\in X$, then $\widehat{f}(x^{^\uparrow})=f(x)$. It follows that $\widehat{f}\circ h=f$.
For the sake of clarity, we'll assume that the above-mentioned family $\mathcal{V}$ satisfies the relation
$\mathcal{V}=\{f_{_a}\vert a\in A\}$
in the next section of the proof.
Let
\begin{center}
$\widehat{\mathcal{V}}=\{\widehat{f}_{_a}\in \delta(X)\vert \widehat{f}_{_a}\circ h=f_{_a}\}$. 
\end{center}
We correspond with
each $\widehat{f}_{_a}\in \widehat{\mathcal{V}}$ the set $\widehat{f}_{_a}(X)=I_{_{\widehat{f}_{_a}}} \subseteq [0,1]$ and denote
$\mathbb{I}_{_{f_{_a}}}=(I_{_{\widehat{f}_{_a}}},\mathfrak{L},\mathfrak{U})=(I_{_{\widehat{f}_{_a}}},\sigma,\omega)$ (remark \ref{opa1}).
Define 
the evaluation map 
\begin{center}
$e: \delta(X) \longrightarrow \displaystyle\prod_{\widehat{f}_{_a}\in \widehat{V}}\mathbb{I}_{_{f_{_a}}}$
\end{center}
induced by the collection $\{\widehat{f}_{_a}\vert a\in A\}$.
Then, for each $x\in X$, $e(x)$ is the point in 
$\displaystyle\prod_{\widehat{f}_{_a}\in \widehat{V}}\mathbb{I}_{_{f_{_a}}}$
whose ath coordinate is $\widehat{f}_{_a}(x)$ for each $a\in A$.
Since each $\widehat{f}_{_a}$ has a range contained in some closed and bounded interval $I_{_{\widehat{f}_{_a}}}$, by Theorem 8.16 in \cite{wil}, $\delta(X)$ is homeomorphic to a subspace of the cube $\displaystyle\prod_{\widehat{f}_{_a}\in \widetilde{V}}I_{_{\widehat{f}_{_a}}}$. Since $\precsim$ is $\sigma\times \omega$-closed, it separates points from closed sets and 
since the
Lawson topology $\lambda=\sigma\vee \omega$ in $\delta(X)$ is $T_1$, by Theorem 8.16 in \cite{wil}, we have that the evaluation map
$e: \delta(X) \longrightarrow \displaystyle\prod_{\widehat{f}_{_a}\in \widehat{V}}I_{_{\widehat{f}_{_a}}}$ is an injection.
On the other hand, since each $\widehat{f}_{_a}$ preserves arbitrary meets and directed  joins so does $e$.
By Lemma \ref{GW} we conclude that the normal completion $(\delta(X),\precsim^{\delta})$ is a continuous lattice.
 Therefore, by Lemma
\ref{aeks} we conclude that $(X,\precsim)$ is precontinuous.
\end{proof}

\par\bigskip\smallskip\par\noindent

\par\noindent
{\it Address}: {\tt {Athanasios Andrikopoulos} \\ {Department of Computer Engineering \& Informatics\\ University of Patras\\ Greece}
\par\noindent
{\it E-mail address}:{\tt aandriko@ceid.upatras.gr}

\end{document}